\documentclass[11pt]{amsart}
\usepackage{texdraw}

\theoremstyle{plain}
\newtheorem{thm}{\bf Theorem}[section]
\newtheorem{df}[thm]{\bf Definition}
\newtheorem{prop}[thm]{\bf Proposition}
\newtheorem{cor}[thm]{\bf Corollary}
\newtheorem{lem}[thm]{\bf Lemma}

\newtheorem{ex}[thm]{\bf Example}
\numberwithin{equation}{section}

\newcommand{\nc}{\newcommand}
\nc{\pf}{{\bfseries Proof. }}
\nc{\tris}{ \bsegment \move(0 0)\lvec(10 0)\lvec(10 10)\lvec(0
0)\ifill f:0.7 \esegment }

\nc{\recs}{ \bsegment \move(0 0)\lvec(10 0)\lvec(10 5)\lvec(0
5)\lvec(0 0)\ifill f:0.7 \esegment }

\newbox{\tmppic}
\newbox{\tmpdraw}
\newbox{\tmpfig}
\newbox{\tmpa}
\newbox{\tmpb}

\begin{document}

\title[Fock space representation for $U_q(C_2^{(1)})$]
      {Fock space representations for the quantum affine algebra $U_q(C_2^{(1)})$}
\author[S.-J. Kang and J.-H. Kwon]{Seok-Jin Kang$^{*}$ and Jae-Hoon Kwon$^{\dagger}$}
\address{School of Mathematics \\
         Korea Institute for Advanced Study \\
         Seoul 130-012, Korea}
\thanks{ $^{*}$This research was supported by KOSEF Grant
\# 98-0701-01-5-L and the Young Scientist Award, Korean Academy of
Science and Technology. \\ $^{\dagger}$This research was supported
by KOSEF Grant \# 98-0701-01-05-L} \email{sjkang@kias.re.kr,
jhkwon@kias.re.kr}

\maketitle

\begin{abstract}
We construct the Fock space representations for the quantum affine
algebra of type $C_2^{(1)}$ in terms of Young walls. Using this
construction, we give a generalized Lascoux-Leclerc-Thibon
algorithm for computing the global bases of the basic
representations.
\end{abstract}

\section{Introduction}
The theory of {\it crystal bases} for the quantum group
$U_q(\frak{g})$ associated with a symmetrizable Kac-Moody algebra
$\frak{g}$ was developed by Kashiwara \cite{Kash,Kash93,Kash94}.
The crystal bases are bases of $U_q(\frak{g})$-modules at $q=0$
and they contain a lot of important information on
$U_q(\frak{g})$-modules. For example, they have oriented graph
structures, called {\it crystal graphs}, which behave very nicely
under the tensor product. As a consequence, many problems in
representation theory are reduced to those in combinatorics. It is
one of the most important problems in the theory of crystal bases
to give various realizations of crystals. There have been many
works on this problem (see for example, \cite{KMN2, KashNakash,
KashSaito, Lit, Nakash1, Nakash2}). From crystal bases, Kashiwara
also recovered true bases of integrable $U_q(\frak{g})$-modules,
called the {\it global bases}, in a canonical way \cite{Kash}.
These bases were proved to be equal to the {\it canonical bases}
constructed by Lustzig in a geometric way \cite{GL,Lusz}.
Recently, using the global bases of the basic representations of
the quantum affine algebra of type $A_n^{(1)}$, Lascoux, Leclerc
and Thibon discovered that there is a deep connection between the
representation theory of the quantum affine algebras and the Hecke
algebras \cite{LLT} (see also \cite{Ariki, BK}).

In this paper, we focus on the quantum affine algebra of type
$C_2^{(1)}$. For a dominant integral weight $\Lambda$ of level 1,
let $B(\Lambda)$ be the crystal of the basic representation
$V(\Lambda)$. In \cite{HK}, Hong and Kang gave a realization of
$B(\Lambda)$ in terms of {\it Young walls} associated with
$U_q(C_2^{(1)})$. These are made by building colored blocks on a
{\it ground state wall $Y_{\Lambda}$} following certain patterns
and rules. Let ${\mathcal Z}(\Lambda)$ be the set of {\it proper
Young walls} on $Y_{\Lambda}$, and ${\mathcal Y}(\Lambda)$ the set
of {\it reduced proper Young walls} on $Y_{\Lambda}$. They gave an
affine crystal structure on ${\mathcal Z}(\Lambda)$, and then
showed that $B(\Lambda)$ is isomorphic to the subcrystal
${\mathcal Y}(\Lambda)$ of ${\mathcal Z}(\Lambda)$.

Let
\begin{equation}
{\mathcal F}(\Lambda)=\bigoplus_{Y\in {\mathcal
Z}(\Lambda)}\mathbb{Q}(q)Y
\end{equation}
be the $\mathbb{Q}(q)$-vector space with a basis ${\mathcal
Z}(\Lambda)$. We give an integrable $U_q(C_2^{(1)})$-module
structure on ${\mathcal F}(\Lambda)$, the {\it Fock space
representation}. We show that the crystal of ${\mathcal
F}(\Lambda)$ is isomorphic to the abstract affine crystal
${\mathcal Z}(\Lambda)$ given in \cite{HK}. Then, by finding all
the maximal vectors in ${\mathcal Z}(\Lambda)$, we obtain a
decomposition of ${\mathcal F}(\Lambda)$ as follows
\begin{equation}
{\mathcal
F}(\Lambda)=\bigoplus_{m=0}^{\infty}V(\Lambda-m\delta)^{\oplus
p(m)}.
\end{equation}

From the embedding of $V(\Lambda)$ into ${\mathcal F}(\Lambda)$,
we show that each global basis element $G(Y)$ ($Y\in {\mathcal
Y}(\Lambda)$) can be written as a $\mathbb{Z}[q]$-linear
combination of proper Young walls which are smaller than or equal
to $Y$ with respect to a certain ordering; that is,
\begin{equation}
G(Y)=Y+\sum_{\substack{Z\in{\mathcal Z}(\Lambda) \\
|Y|\rhd|Z^R|}}G_{Y,Z}(q)Z,
\end{equation}
where  $G_{Y,Z}(q)\in q\mathbb{Z}[q]$ for $Y\neq Z$. We also
discuss an algorithm for computing the coefficient polynomials
$G_{Y,Z}(q)$ in $G(Y)$. This kind of algorithm known as {\it
Lascoux-Leclerc-Thibon algorithm}, was introduced by Lascoux,
Leclerc, and Thibon in case of the quantum affine algebra of type
$A_n^{(1)}$ \cite{LLT}. There are several variants of
Lascoux-Leclerc-Thibon algorithm (see \cite{LT} for classical type
$A_n$, \cite{L} for classical type $C_n$, and \cite{KK} for affine
types $A_{2n-1}^{(1)}$, $A_{2n}^{(2)}$, $B_n^{(1)}$, $D_n^{(1)}$,
$D_{n+1}^{(2)}$). Our results in this paper are based on the work
\cite{KK}.

\section{Quantum affine algebra $U_q(C_2^{(1)})$}
Let $I=\{0,1,2\}$ be the index set. The generalized Cartan matrix
$A=(a_{ij})_{i,j \in I}$ of affine type $C_2^{(1)}$ and its Dynkin
diagram are given by {\savebox{\tmppic}{\begin{texdraw} \drawdim
mm \fontsize{7}{7}\selectfont \setunitscale 0.5 \textref h:C v:C
\move(0 0) \bsegment \move(0 0)\lcir r:2 \move(1.7 1)\lvec(13.3 1)
\move(1.7 -1)\lvec(13.3 -1) \move(15 0)\lcir r:2 \move(16.7
1)\lvec(28.3 1) \move(16.7 -1)\lvec(28.3 -1) \move(30 0)\lcir r:2
\move(13 0)\lvec(10.3 2.7)\move(13 0)\lvec(10.3 -2.7) \move(17
0)\lvec(19.7 2.7)\move(17 0)\lvec(19.7 -2.7) \esegment \htext(0
-5){$0$} \htext(15 -5){$1$} \htext(30 -5){$2$}
\end{texdraw}}%
\begin{equation*}
A =
\begin{pmatrix}
2 & -1 & 0\\
-2 & 2 & -2\\
0 & -1 & 2
\end{pmatrix}
\quad\text{and}\quad \raisebox{-0.4\height}{\usebox{\tmppic}}\,.
\end{equation*}%
}

Let
$P^{\vee}=\mathbb{Z}h_0\oplus\mathbb{Z}h_1\oplus\mathbb{Z}h_2\oplus\mathbb{Z}d$
be a free abelian group, called {\it dual weight lattice} and set
$\frak{h}=\mathbb{Q}\otimes_{\mathbb{Z}}P^{\vee}$. For $i\in I$,
we define $\alpha_i$ and $\Lambda_i\in\frak{h}^*$ by
\begin{equation*}
\begin{aligned}
\alpha_i (h_j) & = a_{ji}, \quad \alpha_i (d) = \delta_{0, i}, \\
\Lambda_i(h_j) & = \delta_{ij}, \quad \Lambda_i(d)=0 \qquad (i,j
\in I).
\end{aligned}
\end{equation*}
The $\alpha_i$ are called the {\it simple roots} and the
$\Lambda_i$ are called the {\it fundamental weights}.

Let $c=h_0 + h_1 + h_2$ and $\delta = \alpha_0 + 2\alpha_1 +
\alpha_2$. Then we have $\alpha_i(c)=0$, $\delta(h_i)=0$ for all
$i\in I$ and $\delta(d)=1$. We call $c$ (resp. $\delta$) the {\it
canonical central element} (resp. {\it null root}). The free
abelian group $P=\mathbb{Z} \Lambda_0 \oplus \mathbb{Z} \Lambda_1
\oplus \mathbb{Z} \Lambda_2 \oplus \mathbb{Z} \delta$ is called
the {\it weight lattice}.

Let $q$ be an indeterminate. We denote by $q^h$ $(h\in P^{\vee})$
the basis elements of the group algebra $\mathbb{Q}(q)[P^{\vee}]$
with the multiplication $q^h q^{h'} = q^{h+h'}$ $(h,h'\in
P^{\vee})$. Set $q_0=q_2=q^2$, $q_1=q$ and $K_i=q_i^{h_i}$ ($i\in
I$). We will also use the following notations.
\begin{equation*}
[k]_i = \frac{q_i^k - q_i^{-k}}{q_i - q_i^{-1}}, \quad [n]_i ! =
\prod_{k=1}^{n} [k]_i, \quad \text{and} \quad e_i^{(n)} =
e_i^n/[n]_i!, f_i^{(n)} = f_i^n/[n]_i!.
\end{equation*}

\begin{df}{\rm
The {\em quantum affine algebra $U_q(C_2^{(1)})$} is the
associative algebra with 1 over $\mathbb{Q}(q)$ generated by the
elements $e_i$, $f_i$ $(i\in I)$ and $q^h$ $(h\in P^{\vee})$
subject to the following defining relations:}
\end{df}
\begin{equation*}
\begin{aligned}
\ & q^0 = 1, \ \ q^h q^{h'} = q^{h+ h'} \quad (h, h'\in P^{\vee}),\\
\ & q^h e_i q^{-h} = q^{\alpha_i(h)} e_i, \quad
q^h f_i q^{-h} = q^{-\alpha_i(h)} f_i \quad (h\in P^{\vee}, i\in I), \\
\ & e_i f_j - f_j e_i = \delta_{i,j} \frac{K_i - K_i^{-1}}{q_i -
q_i^{-1}}
\quad (i,j \in I), \\
\ & e_0^2e_1 - (q^2 + q^{-2})e_0 e_1 e_0 + e_1e_0^2 = 0, \\
\ & f_0^2f_1 - (q^2 + q^{-2})f_0 f_1 f_0 + f_1f_0^2 = 0, \\
\ & e_1^3e_0 - (q^2 + 1 + q^{-2})e_1^2e_0e_1
      + (q^2 + 1 + q^{-2})e_1e_0e_1^2 - e_0e_1^3 = 0, \\
\ & f_1^3f_0 - (q^2 + 1 + q^{-2})f_1^2f_0f_1
      + (q^2 + 1 + q^{-2})f_1f_0f_1^2 - f_0f_1^3 = 0, \\
\ & e_1^3e_2 - (q^2 + 1 + q^{-2})e_1^2e_2e_1
      + (q^2 + 1 + q^{-2})e_1e_2e_1^2 - e_2e_1^3 = 0, \\
\ & f_1^3f_2 - (q^2 + 1 + q^{-2})f_1^2f_2f_1
      + (q^2 + 1 + q^{-2})f_1f_2f_1^2 - f_2f_1^3 = 0, \\
\ & e_2^2e_1 - (q^2 + q^{-2})e_2 e_1 e_2 + e_1e_2^2 = 0, \\
\ & f_2^2f_1 - (q^2 + q^{-2})f_2 f_1 f_2 + f_1f_2^2 = 0, \\
\ & e_0e_2 = e_2e_0, \quad f_0f_2 = f_2f_0.
\end{aligned}
\end{equation*}
It is also called the {\it quantum affine algebra of type
$C_2^{(1)}$.}

\section{Crystal bases}
In this section, we review the crystal basis theory for the
quantum affine algebra $U_q(C_2^{(1)})$. All the statements and
the results in this section hold for a quantized enveloping
algebra associated with a symmetrizable Kac-Moody algebra (see
\cite{Kash}). A $U_q(C_2^{(1)})$-module $M$ is called {\it
integrable} if

\begin{itemize}
\item[(i)] $M=\bigoplus_{\lambda \in P} M_{\lambda}$ where
$M_{\lambda} = \{ v\in M \mid q^h v = q^{\lambda(h)} v \text{ for
all }h\in P^{\vee}\}$,

\item[(ii)] $M$ is a direct sum of finite
dimensional irreducible $U_i$-modules, where $U_i$ ($i\in I$) is
the subalgebra generated by $e_i$, $f_i$, $K_i^{\pm1}$.
\end{itemize}

Fix $i\in I$. An element $v\in M_{\lambda}$ may be written
uniquely as
\begin{equation*}
v = \sum_{k\geq0} f_i^{(k)} v_k,
\end{equation*}
where $v_k \in \ker e_i \cap M_{\lambda+k\alpha_i}$. We define the
endomorphisms $\tilde{e}_i$ and $\tilde{f}_i$ on $M$, called the
{\it Kashiwara operators},  by
\begin{equation*}
\tilde{e}_i v = \sum_{k\geq1} f_i^{(k-1)} v_k, \qquad \tilde{f}_i
v = \sum_{k\geq0} f_i^{(k+1)} v_k.
\end{equation*}

Let $\mathbb{A}_0=\{\,f/g\in\mathbb{Q}(q)\,|\,f,g\in\mathbb{Q}[q],
g(0)\neq 0 \,\}$ be the localization of $\mathbb{Q}[q]$ at $q=0$.
\begin{df}\label {crystal basis}{\rm
A {\it crystal basis} of $M$ is a pair $(L,B)$, where
\begin{itemize}
\item[(i)] $L$ is a free $\mathbb{A}_0$-submodule of $M$ such that
$M\cong\mathbb{Q}(q)\otimes_{\mathbb{A}_0}L$,

\item[(ii)] $B$ is a $\mathbb{Q}$-basis of $L/qL$,

\item[(iii)] $L=\bigoplus_{\lambda\in P}L_{\lambda}$, where
$L_{\lambda}=L\cap M_{\lambda}$,

\item[(iv)] $B=\bigsqcup_{\lambda\in P}B_{\lambda}$, where
$B_{\lambda}=B\cap(L_{\lambda}/qL_{\lambda})$,

\item[(v)] $\tilde{e}_iL\subset L$, $\tilde{f}_iL\subset L$ for all
$i\in I$,

\item[(vi)] $\tilde{e}_iB\subset B\cup\{0\}$, $\tilde{f}_iB\subset
B\cup\{0\}$ for all $i\in I$,

\item[(vii)] for $b,b'\in B$, $\tilde{f}_ib=b'$ if and
only if $b= \tilde{e}_ib'$.
\end{itemize}}
\end{df}
The set $B$ becomes a colored oriented graph, called the {\it
crystal graph}, where the arrows are defined by
$b\stackrel{i}{\rightarrow}b'$ if and only if $\tilde{f}_i b=b'$,
for $b,b'\in B$.

For each $b\in B$ and $i\in I$, we define
$\varepsilon_i(b)=\max\{\,k\geq 0\,|\,\tilde{e}_i^kb\in B\,\}$,
$\varphi_i(b)=\max\{\,k\geq 0\,|\,\tilde{f}_i^kb\in B\,\}$. Then
we have
\begin{equation}
\begin{split}
&\varphi_i(b)=\varepsilon_i(b)+\langle h_i,{\rm wt }(b)\rangle \\
&{\rm wt}(\tilde{e}_ib)={\rm wt}(\tilde{e}_ib)+\alpha_i,\quad
{\rm wt}(\tilde{f}_ib)={\rm wt}(\tilde{f}_ib)-\alpha_i,\\
&\varepsilon_i(\tilde{e}_ib)=\varepsilon_i(b)-1,\quad
\varphi_i(\tilde{e}_ib)=\varphi_i(b)+1\quad
\text{ if $\tilde{e}_ib\in B$}, \\
&\varepsilon_i(\tilde{f}_ib)=\varepsilon_i(b)+1,\quad
\varphi_i(\tilde{f}_ib)=\varphi_i(b)-1\quad \text{ if
$\tilde{f}_ib\in B$}.
\end{split}
\end{equation}

Set $P^{+}=\{\,\lambda \in
\frak{h}^*\,|\,\text{$\lambda(h_i)\in\mathbb{Z}_{\geq 0
}$,}\,\,\,i\in I\,\}$. For $\lambda\in P^+$, let $V(\lambda)$ be
the irreducible highest weight $U_q(C_2^{(1)})$-module with
highest weight $\lambda$ and highest weight vector $u_{\lambda}$.

\begin{thm}\label{crystal basis}{\rm \cite{Kash}}
Let $L(\lambda)$ be the free $\mathbb{A}_0$-submodule of
$V(\lambda)$ spanned by the vectors of the form
$\tilde{f}_{i_1}\cdots\tilde{f}_{i_r}u_{\lambda}$ {\rm (}$i_k\in
I$, $r\in\mathbb{Z}_{\geq0}${\rm )} and set $
B(\lambda)=\{\,\tilde{f}_{i_1}\cdots\tilde{f}_{i_r}u_{\lambda}+qL(\lambda)\in
L(\lambda)/qL(\lambda)\,\}\backslash\{0\}$. Then
$(L(\lambda),B(\lambda))$ is a crystal basis of $V(\lambda)$, and
every crystal basis of $V(\lambda)$ is isomorphic to
$(L(\lambda),B(\lambda))$.\qed
\end{thm}

There exists an involution $-$ of $U_q(C_2^{(1)})$ as a
$\mathbb{Q}$-algebra defined  by
\begin{gather}
\overline{e}_i=e_i,\quad \overline{f}_i=f_i,\quad
\overline{q^h}=q^{-h},\quad \overline{q}=q^{-1}
\end{gather}
for $i\in I$ and $h\in P^{\vee}$. Set
$\mathbb{A}=\mathbb{Q}[q,q^{-1}]$. We denote by
$U^-_{\mathbb{A}}(C_2^{(1)})$ the $\mathbb{A}$-subalgebra of
$U_q(C_2^{(1)})$ generated by $f^{(n)}_i$ {\rm ($i\in I$,
$n\in\mathbb{Z}_{\geq 0}$)}. Set
$V(\lambda)^{\mathbb{A}}=U_{\mathbb{A}}^{-}(C_2^{(1)})u_{\lambda}$.

\begin{thm}\label {global basis} {\rm \cite{Kash}}
There exists a unique $\mathbb{A}$-basis
$G(\lambda)=\{\,G(b)\,|\,b\in B(\lambda)\,\}$ of
$V(\lambda)^{\mathbb{A}}$ such that
\begin{equation*}
G(b)\equiv b \mod {qL(\lambda)}\quad \text{and} \quad
\overline{G(b)}=G(b)
\end{equation*}
for all $b\in B(\lambda)$.\qed
\end{thm}
The basis $G(\lambda)$ is called the {\it global basis} or {\it
canonical basis} of $V(\lambda)$ associated with the crystal graph
$B(\lambda)$.

By extracting the properties of crystal graphs, we can define the
notion of abstract {\it crystals} \cite{Kash93,Kash94}.
\begin{df}{\rm
An {\it affine crystal} is a set $B$ together with the maps ${\rm
wt} : B\rightarrow P$, $\varepsilon_i, \varphi_i : B\rightarrow
\mathbb{Z}\cup\{-\infty\}$, $\tilde{e}_i,\tilde{f}_i :
B\rightarrow B\cup\{0\}$ ($i\in I$) such that for $i\in I$ and
$b\in B$,
\begin{itemize}
\item[(i)] $\varphi_i(b)=\varepsilon_i(b)+\langle h_i,{\rm wt
}(b)\rangle $,

\item[(ii)] ${\rm wt}(\tilde{e}_ib)={\rm wt}(b)+\alpha_i$, ${\rm
wt}(\tilde{f}_ib)={\rm wt}(b)-\alpha_i$,

\item[(iii)] $\varepsilon_i(\tilde{e}_ib)=\varepsilon_i(b)-1$,
$\varphi_i(\tilde{e}_ib)=\varphi_i(b)+1$ if $\tilde{e}_ib\in B$,

\item[(iv)] $\varepsilon_i(\tilde{f}_ib)=\varepsilon_i(b)+1$,
$\varphi_i(\tilde{f}_ib)=\varphi_i(b)-1$ if $\tilde{f}_ib\in B$,

\item[(v)] $\tilde{f}_ib=b'$ if and
only if $\tilde{e}_ib'=b$ for $b,b'\in B$,

\item[(vi)] $\tilde{e}_ib=\tilde{f}_ib=0$ if $\varepsilon_i(b)=-\infty$.
\end{itemize}}
\end{df}
The crystal $B(\lambda)$ of $V(\lambda)$ ($\lambda\in P^+$)
satisfies the above conditions and it is an affine crystal.

 A {\it morphism} $\psi:B_1 \rightarrow B_2$ of
crystals is a map $\psi:B_1 \cup \{0\} \rightarrow B_2 \cup \{0\}$
satisfying the conditions:
\begin{itemize}
\item[(i)] $\psi(0)=0$,

\item[(ii)] ${\rm wt}(\psi(b))= {\rm wt}(b)$,
$\varepsilon_i(\psi(b))=\varepsilon_i(b)$, and $\varphi_i(\psi(b))
= \varphi_i(b)$ if $b\in B_1$ and $\psi(b) \in B_2$,

\item[(iii)] $\tilde{f}_i \psi(b) = \psi(b')$ if $b,
b'\in B_1$, $\psi(b), \psi(b') \in B_2$ and $\tilde{f}_i b =b'$.
\end{itemize}

\section{Young walls}
In this section, we will give a brief review of the results in
\cite{HK}. The {\it Young walls} will be built of two kinds of
\emph{blocks}; \vskip 3mm
\renewcommand{\arraystretch}{2}
\begin{center}
\begin{tabular}{c|c|c|c|c|c}
type & shape & width & thickness & height & volume \\
\hline I & \raisebox{-0.4\height}{
\begin{texdraw}
\drawdim em \setunitscale 0.1 \linewd 0.5 \move(0 0)\lvec(10
0)\lvec(10 5)\lvec(0 5)\lvec(0 0) \move(10 0)\lvec(15 5)\lvec(15
10)\lvec(5 10)\lvec(0 5) \move(10 5)\lvec(15 10)
\end{texdraw}
} $=$ \raisebox{-0.4\height}{
\begin{texdraw}
\drawdim em \setunitscale 0.1 \linewd 0.5 \textref h:C v:C \move(0
0)\lvec(10 0)\lvec(10 5)\lvec(0 5)\lvec(0 0)
\end{texdraw}
}
& 1 & 1 & $\frac{1}{2}$ & $\frac{1}{2}$ \\
II &

\raisebox{-0.4\height}{
\begin{texdraw}
\drawdim em \setunitscale 0.1 \linewd 0.5 \move(-10 0)\lvec(0
0)\lvec(0 10)\lvec(-10 10)\lvec(-10 0) \move(0 0)\lvec(2.5
2.5)\lvec(2.5 12.5)\lvec(-7.5 12.5)\lvec(-10 10) \move(0
10)\lvec(2.5 12.5) \lpatt(0.3 1) \move(0 0)\lvec(-2.5
-2.5)\lvec(-12.5 -2.5)\lvec(-10 0)
\end{texdraw}
} $=$ \raisebox{-0.4\height}{
\begin{texdraw}
\drawdim em \setunitscale 0.1 \linewd 0.5 \move(-10 0)\lvec(0
0)\lvec(0 10)\lvec(-10 0) \lpatt(0.3 1)\move(-10 0)\lvec(-10
10)\lvec(0 10)
\end{texdraw}}
, \hskip 2mm \raisebox{-0.4\height}{
\begin{texdraw}
\drawdim em \setunitscale 0.1 \linewd 0.5 \move(-10 0)\lvec(0
0)\lvec(0 10)\lvec(-10 10)\lvec(-10 0) \move(0 0)\lvec(2.5
2.5)\lvec(2.5 12.5)\lvec(-7.5 12.5)\lvec(-10 10) \move(0
10)\lvec(2.5 12.5) \lpatt(0.3 1) \move(2.5 2.5)\lvec(5 5)\lvec(2.5
5)
\end{texdraw}
} $=$ \raisebox{-0.4\height}{
\begin{texdraw}
\drawdim em \setunitscale 0.1 \linewd 0.5 \move(-10 0)\lvec(-10
10)\lvec(0 10)\lvec(-10 0)\lpatt(0.3 1)\move(-10 0)\lvec(0
0)\lvec(0 10)
\end{texdraw}
}

& 1 & $\frac{1}{2}$ & 1 & $\frac{1}{2}$ \\

\end{tabular}
\end{center}
\renewcommand{\arraystretch}{1}
\vskip 3mm We also give a coloring of blocks as follows;

\vskip 3mm
\begin{center}
\raisebox{-0.33\height}[0.69\height][0.35\height]{%
\begin{texdraw}
\drawdim em \setunitscale 0.1 \linewd 0.5 \move(0 0)\lvec(10
0)\lvec(10 10)\lvec(0 10)\lvec(0 0) \move(10 0)\lvec(12.5
2.5)\lvec(12.5 12.5)\lvec(2.5 12.5)\lvec(0 10) \move(10
10)\lvec(12.5 12.5) \htext(3 3){$_0$}
\end{texdraw}%
} \hskip 5mm
\raisebox{-0.33\height}[0.69\height][0.35\height]{%
\begin{texdraw}
\drawdim em \setunitscale 0.1 \linewd 0.5 \move(0 0)\lvec(10
0)\lvec(10 10)\lvec(0 10)\lvec(0 0) \move(10 0)\lvec(12.5
2.5)\lvec(12.5 12.5)\lvec(2.5 12.5)\lvec(0 10) \move(10
10)\lvec(12.5 12.5) \htext(3 3){$_2$}
\end{texdraw}%
} \hskip 5mm
\raisebox{-0.33\height}[0.69\height][0.35\height]{%
\begin{texdraw}
\drawdim em \setunitscale 0.1 \linewd 0.5 \move(0 0)\lvec(10
0)\lvec(10 5)\lvec(0 5)\lvec(0 0) \move(10 0)\lvec(15 5)\lvec(15
10)\lvec(5 10)\lvec(0 5) \move(10 5)\lvec(15 10) \htext(3 1){\tiny
$1$}
\end{texdraw}%
}
\end{center}
\vskip 3mm

Given a dominant integral weight $\Lambda=\Lambda_i$ ($i\in I$) of
level 1, that is, $\Lambda(c)=1$, we fix a frame $Y_{\Lambda}$
called the {\it ground state wall} of weight $\Lambda$ as follows

\vskip 3mm
\begin{center}
\begin{tabular}{rcl}
$Y_{\Lambda_0}$ & $=$ & \raisebox{-0.4\height}{
\begin{texdraw}
\drawdim mm \setunitscale 0.5 \fontsize{7}{7}\selectfont \textref
h:C v:C \move(0 0)\lvec(-60 0)\move(0 10)\lvec(-60 10) \move(-57.5
12.5)\lvec(2.5 12.5)\lvec(2.5 2.5)\lvec(0 0) \move(0 0) \bsegment
\move(0 0)\lvec(0 10)\lvec(2.5 12.5) \esegment \move(-10 0)
\bsegment \move(0 0)\lvec(0 10)\lvec(2.5 12.5) \esegment \move(-20
0) \bsegment \move(0 0)\lvec(0 10)\lvec(2.5 12.5) \esegment
\move(-30 0) \bsegment \move(0 0)\lvec(0 10)\lvec(2.5 12.5)
\esegment \move(-40 0) \bsegment \move(0 0)\lvec(0 10)\lvec(2.5
12.5) \esegment \move(-50 0) \bsegment \move(0 0)\lvec(0
10)\lvec(2.5 12.5) \esegment \move(-60 0) \bsegment \move(0
0)\lvec(0 10)\lvec(2.5 12.5) \esegment \move(0 0) \bsegment
\lpatt(0.3 1) \move(0 0)\lvec(-2.5 -2.5)\lvec(-67.5 -2.5)
\move(-60 0)\rlvec(-5 0)\move(-60 10)\rlvec(-5 0)\move(-57.5
12.5)\rlvec(-5 0) \move(0 0)\rlvec(-2.5 -2.5) \move(-10
0)\rlvec(-2.5 -2.5) \move(-20 0)\rlvec(-2.5 -2.5) \move(-30
0)\rlvec(-2.5 -2.5) \move(-40 0)\rlvec(-2.5 -2.5) \move(-50
0)\rlvec(-2.5 -2.5) \move(-60 0)\rlvec(-2.5 -2.5) \esegment
\move(0 0) \bsegment \htext(-5 5){$2$} \htext(-15 5){$0$}
\htext(-25 5){$2$} \htext(-35 5){$0$} \htext(-45 5){$2$}
\htext(-55 5){$0$} \esegment
\end{texdraw}
}
\\[4mm]
$Y_{\Lambda_1}$ & $=$ & \raisebox{-0.4\height}{
\begin{texdraw}
\drawdim mm \setunitscale 0.5 \fontsize{7}{7}\selectfont \textref
h:C v:C \move(-60 0)\lvec(0 0)\move(0 5)\lvec(-60 5) \move(0
0)\lvec(5 5)\lvec(5 10)\lvec(-55 10) \move(0 0)\bsegment\move(0
0)\lvec(0 5)\lvec(5 10)\esegment \move(-10 0)\bsegment\move(0
0)\lvec(0 5)\lvec(5 10)\esegment \move(-20 0)\bsegment\move(0
0)\lvec(0 5)\lvec(5 10)\esegment \move(-30 0)\bsegment\move(0
0)\lvec(0 5)\lvec(5 10)\esegment \move(-40 0)\bsegment\move(0
0)\lvec(0 5)\lvec(5 10)\esegment \move(-50 0)\bsegment\move(0
0)\lvec(0 5)\lvec(5 10)\esegment \move(-60 0)\bsegment\move(0
0)\lvec(0 5)\lvec(5 10)\esegment \move(0 0) \bsegment \lpatt(0.3
1) \move(-60 0)\rlvec(-5 0)\move(-60 5)\rlvec(-5 0)\move(-55
10)\rlvec(-5 0) \esegment \move(0 0) \bsegment \htext(-5 2.5){$1$}
\htext(-15 2.5){$1$} \htext(-25 2.5){$1$} \htext(-35 2.5){$1$}
\htext(-45 2.5){$1$} \htext(-55 2.5){$1$} \esegment
\end{texdraw}
}
\\[4mm]
$Y_{\Lambda_2}$ & $=$ & \raisebox{-0.4\height}{
\begin{texdraw}
\drawdim mm \setunitscale 0.5 \fontsize{7}{7}\selectfont \textref
h:C v:C \move(0 0)\lvec(-60 0)\move(0 10)\lvec(-60 10) \move(-57.5
12.5)\lvec(2.5 12.5)\lvec(2.5 2.5)\lvec(0 0) \move(0 0) \bsegment
\move(0 0)\lvec(0 10)\lvec(2.5 12.5) \esegment \move(-10 0)
\bsegment \move(0 0)\lvec(0 10)\lvec(2.5 12.5) \esegment \move(-20
0) \bsegment \move(0 0)\lvec(0 10)\lvec(2.5 12.5) \esegment
\move(-30 0) \bsegment \move(0 0)\lvec(0 10)\lvec(2.5 12.5)
\esegment \move(-40 0) \bsegment \move(0 0)\lvec(0 10)\lvec(2.5
12.5) \esegment \move(-50 0) \bsegment \move(0 0)\lvec(0
10)\lvec(2.5 12.5) \esegment \move(-60 0) \bsegment \move(0
0)\lvec(0 10)\lvec(2.5 12.5) \esegment \move(0 0) \bsegment
\lpatt(0.3 1) \move(0 0)\lvec(-2.5 -2.5)\lvec(-67.5 -2.5)
\move(-60 0)\rlvec(-5 0)\move(-60 10)\rlvec(-5 0)\move(-57.5
12.5)\rlvec(-5 0) \move(0 0)\rlvec(-2.5 -2.5) \move(-10
0)\rlvec(-2.5 -2.5) \move(-20 0)\rlvec(-2.5 -2.5) \move(-30
0)\rlvec(-2.5 -2.5) \move(-40 0)\rlvec(-2.5 -2.5) \move(-50
0)\rlvec(-2.5 -2.5) \move(-60 0)\rlvec(-2.5 -2.5) \esegment
\move(0 0) \bsegment \htext(-5 5){$0$} \htext(-15 5){$2$}
\htext(-25 5){$0$} \htext(-35 5){$2$} \htext(-45 5){$0$}
\htext(-55 5){$2$} \esegment
\end{texdraw}
}
\end{tabular}
\end{center}
\vskip 3mm

On this frame, we build a wall of thickness less than or equal to
one unit. The rules for building walls are as follows:
\begin{itemize}
\item[(1)] The colored blocks should be stacked in columns. No block
can be placed on top of a column of half-unit thickness.

\item[(2)] Except for the right-most column, there should be no free
space to the right of any block.

\item[(3)] The colored blocks should be stacked in a specified
pattern which is determined as follows;\vskip 5mm

\begin{center}
$\Lambda_0$: \raisebox{-1\height}{\begin{texdraw} \drawdim mm
\setunitscale 0.5 \fontsize{7}{7}\selectfont \textref h:C v:C

\move(0 0)\tris \move(10 0)\tris \move(-20 0)\tris \move(-10
0)\tris

\move(0 0) \bsegment \move(0 0)\rlvec(20 0) \move(0 10)\rlvec(20
0) \move(0 15)\rlvec(20 0) \move(0 20)\rlvec(20 0) \move(0
30)\rlvec(20 0) \move(0 35)\rlvec(20 0) \move(0 40)\rlvec(20 0)
\move(0 0)\rlvec(0 40) \move(10 0)\rlvec(0 40) \move(20 0)\rlvec(0
40) \move(0 0)\rlvec(10 10) \move(10 0)\rlvec(10 10) \move(0
20)\rlvec(10 10) \move(10 20)\rlvec(10 10)

\htext(2.5 7.5){$2$} \htext(7.5 2.5){$0$} \htext(12.5 7.5){$0$}
\htext(17.5 2.5){$2$} \htext(2.5 27.5){$2$} \htext(7.5 22.5){$0$}
\htext(12.5 27.5){$0$} \htext(17.5 22.5){$2$} \htext(5 12.5){$1$}
\htext(15 12.5){$1$} \htext(5 17.5){$1$} \htext(15 17.5){$1$}
\htext(5 32.5){$1$} \htext(15 32.5){$1$} \htext(5 37.5){$1$}
\htext(15 37.5){$1$} \esegment

\move (0 40) \bsegment \move(0 0)\rlvec(20 0) \move(0 10)\rlvec(20
0) \move(0 15)\rlvec(20 0) \move(0 20)\rlvec(20 0) \move(0
30)\rlvec(20 0) \move(0 35)\rlvec(20 0) \move(0 40)\rlvec(20 0)
\move(0 0)\rlvec(0 40) \move(10 0)\rlvec(0 40) \move(20 0)\rlvec(0
40) \move(0 0)\rlvec(10 10) \move(10 0)\rlvec(10 10) \move(0
20)\rlvec(10 10) \move(10 20)\rlvec(10 10)

\htext(2.5 7.5){$2$} \htext(7.5 2.5){$0$} \htext(12.5 7.5){$0$}
\htext(17.5 2.5){$2$} \htext(2.5 27.5){$2$} \htext(7.5 22.5){$0$}
\htext(12.5 27.5){$0$} \htext(17.5 22.5){$2$} \htext(5 12.5){$1$}
\htext(15 12.5){$1$} \htext(5 17.5){$1$} \htext(15 17.5){$1$}
\htext(5 32.5){$1$} \htext(15 32.5){$1$} \htext(5 37.5){$1$}
\htext(15 37.5){$1$} \esegment

\move(-20 0) \bsegment \move(0 0)\rlvec(20 0) \move(0 10)\rlvec(20
0) \move(0 15)\rlvec(20 0) \move(0 20)\rlvec(20 0) \move(0
30)\rlvec(20 0) \move(0 35)\rlvec(20 0) \move(0 40)\rlvec(20 0)
\move(0 0)\rlvec(0 40) \move(10 0)\rlvec(0 40) \move(20 0)\rlvec(0
40) \move(0 0)\rlvec(10 10) \move(10 0)\rlvec(10 10) \move(0
20)\rlvec(10 10) \move(10 20)\rlvec(10 10) \htext(2.5 7.5){$2$}
\htext(7.5 2.5){$0$} \htext(12.5 7.5){$0$} \htext(17.5 2.5){$2$}
\htext(2.5 27.5){$2$} \htext(7.5 22.5){$0$} \htext(12.5 27.5){$0$}
\htext(17.5 22.5){$2$} \htext(5 12.5){$1$} \htext(15 12.5){$1$}
\htext(5 17.5){$1$} \htext(15 17.5){$1$} \htext(5 32.5){$1$}
\htext(15 32.5){$1$} \htext(5 37.5){$1$} \htext(15 37.5){$1$}
\esegment

\move (-20 40) \bsegment \move(0 0)\rlvec(20 0) \move(0
10)\rlvec(20 0) \move(0 15)\rlvec(20 0) \move(0 20)\rlvec(20 0)
\move(0 30)\rlvec(20 0) \move(0 35)\rlvec(20 0) \move(0
40)\rlvec(20 0) \move(0 0)\rlvec(0 40) \move(10 0)\rlvec(0 40)
\move(20 0)\rlvec(0 40) \move(0 0)\rlvec(10 10) \move(10
0)\rlvec(10 10) \move(0 20)\rlvec(10 10) \move(10 20)\rlvec(10 10)
\htext(2.5 7.5){$2$} \htext(7.5 2.5){$0$} \htext(12.5 7.5){$0$}
\htext(17.5 2.5){$2$} \htext(2.5 27.5){$2$} \htext(7.5 22.5){$0$}
\htext(12.5 27.5){$0$} \htext(17.5 22.5){$2$} \htext(5 12.5){$1$}
\htext(15 12.5){$1$} \htext(5 17.5){$1$} \htext(15 17.5){$1$}
\htext(5 32.5){$1$} \htext(15 32.5){$1$} \htext(5 37.5){$1$}
\htext(15 37.5){$1$}\esegment

\move(20 0)\lvec(-25 0) \move(20 10)\lvec(-25 10) \move(20
20)\lvec(-25 20) \move(20 30)\lvec(-25 30) \move(20 40)\lvec(-25
40) \move(20 50)\lvec(-25 50) \move(20 60)\lvec(-25 60)
 \move(20 70)\lvec(-25 70) \move(20 80)\lvec(-25 80)

\move(20 0)\lvec(20 85) \move(10 0)\lvec(10 85)\move(0 0)\lvec(0
85)\move(-10 0)\lvec(-10 85)\move(-20 0)\lvec(-20 85)

\end{texdraw}}
\hskip 2cm
$\Lambda_2$: \raisebox{-1\height}{\begin{texdraw} \drawdim mm
\setunitscale 0.5 \fontsize{7}{7}\selectfont \textref h:C v:C

\move(0 0)\tris \move(10 0)\tris \move(-20 0)\tris \move(-10
0)\tris

\move(0 0) \bsegment \move(0 0)\rlvec(20 0) \move(0 10)\rlvec(20
0) \move(0 15)\rlvec(20 0) \move(0 20)\rlvec(20 0) \move(0
30)\rlvec(20 0) \move(0 35)\rlvec(20 0) \move(0 40)\rlvec(20 0)
\move(0 0)\rlvec(0 40) \move(10 0)\rlvec(0 40) \move(20 0)\rlvec(0
40) \move(0 0)\rlvec(10 10) \move(10 0)\rlvec(10 10) \move(0
20)\rlvec(10 10) \move(10 20)\rlvec(10 10) \htext(2.5 7.5){$0$}
\htext(7.5 2.5){$2$} \htext(12.5 7.5){$2$} \htext(17.5 2.5){$0$}
\htext(2.5 27.5){$0$} \htext(7.5 22.5){$2$} \htext(12.5 27.5){$2$}
\htext(17.5 22.5){$0$} \htext(5 12.5){$1$} \htext(15 12.5){$1$}
\htext(5 17.5){$1$} \htext(15 17.5){$1$} \htext(5 32.5){$1$}
\htext(15 32.5){$1$} \htext(5 37.5){$1$} \htext(15 37.5){$1$}
\esegment

\move (0 40) \bsegment \move(0 0)\rlvec(20 0) \move(0 10)\rlvec(20
0) \move(0 15)\rlvec(20 0) \move(0 20)\rlvec(20 0) \move(0
30)\rlvec(20 0) \move(0 35)\rlvec(20 0) \move(0 40)\rlvec(20 0)
\move(0 0)\rlvec(0 40) \move(10 0)\rlvec(0 40) \move(20 0)\rlvec(0
40) \move(0 0)\rlvec(10 10) \move(10 0)\rlvec(10 10) \move(0
20)\rlvec(10 10) \move(10 20)\rlvec(10 10) \htext(2.5 7.5){$0$}
\htext(7.5 2.5){$2$} \htext(12.5 7.5){$2$} \htext(17.5 2.5){$0$}
\htext(2.5 27.5){$0$} \htext(7.5 22.5){$2$} \htext(12.5 27.5){$2$}
\htext(17.5 22.5){$0$} \htext(5 12.5){$1$} \htext(15 12.5){$1$}
\htext(5 17.5){$1$} \htext(15 17.5){$1$} \htext(5 32.5){$1$}
\htext(15 32.5){$1$} \htext(5 37.5){$1$} \htext(15 37.5){$1$}
\esegment

\move(-20 0) \bsegment \move(0 0)\rlvec(20 0) \move(0 10)\rlvec(20
0) \move(0 15)\rlvec(20 0) \move(0 20)\rlvec(20 0) \move(0
30)\rlvec(20 0) \move(0 35)\rlvec(20 0) \move(0 40)\rlvec(20 0)
\move(0 0)\rlvec(0 40) \move(10 0)\rlvec(0 40) \move(20 0)\rlvec(0
40) \move(0 0)\rlvec(10 10) \move(10 0)\rlvec(10 10) \move(0
20)\rlvec(10 10) \move(10 20)\rlvec(10 10) \htext(2.5 7.5){$0$}
\htext(7.5 2.5){$2$} \htext(12.5 7.5){$2$} \htext(17.5 2.5){$0$}
\htext(2.5 27.5){$0$} \htext(7.5 22.5){$2$} \htext(12.5 27.5){$2$}
\htext(17.5 22.5){$0$} \htext(5 12.5){$1$} \htext(15 12.5){$1$}
\htext(5 17.5){$1$} \htext(15 17.5){$1$} \htext(5 32.5){$1$}
\htext(15 32.5){$1$} \htext(5 37.5){$1$} \htext(15 37.5){$1$}
\esegment

\move (-20 40) \bsegment \move(0 0)\rlvec(20 0) \move(0
10)\rlvec(20 0) \move(0 15)\rlvec(20 0) \move(0 20)\rlvec(20 0)
\move(0 30)\rlvec(20 0) \move(0 35)\rlvec(20 0) \move(0
40)\rlvec(20 0) \move(0 0)\rlvec(0 40) \move(10 0)\rlvec(0 40)
\move(20 0)\rlvec(0 40) \move(0 0)\rlvec(10 10) \move(10
0)\rlvec(10 10) \move(0 20)\rlvec(10 10) \move(10 20)\rlvec(10 10)
\htext(2.5 7.5){$0$} \htext(7.5 2.5){$2$} \htext(12.5 7.5){$2$}
\htext(17.5 2.5){$0$} \htext(2.5 27.5){$0$} \htext(7.5 22.5){$2$}
\htext(12.5 27.5){$2$} \htext(17.5 22.5){$0$} \htext(5 12.5){$1$}
\htext(15 12.5){$1$} \htext(5 17.5){$1$} \htext(15 17.5){$1$}
\htext(5 32.5){$1$} \htext(15 32.5){$1$} \htext(5 37.5){$1$}
\htext(15 37.5){$1$} \esegment

\move(20 0)\lvec(-25 0) \move(20 10)\lvec(-25 10) \move(20
20)\lvec(-25 20) \move(20 30)\lvec(-25 30) \move(20 40)\lvec(-25
40) \move(20 50)\lvec(-25 50) \move(20 60)\lvec(-25 60)
 \move(20 70)\lvec(-25 70) \move(20 80)\lvec(-25 80)

\move(20 0)\lvec(20 85) \move(10 0)\lvec(10 85)\move(0 0)\lvec(0
85)\move(-10 0)\lvec(-10 85)\move(-20 0)\lvec(-20 85)
\end{texdraw}}

\end{center}
\vskip 5mm
\begin{center}
$\Lambda_1$: \raisebox{-1\height}{\begin{texdraw} \drawdim mm
\setunitscale 0.5 \fontsize{7}{7}\selectfont \textref h:C v:C

\move(0 0)\recs  \move(10 0)\recs \move(-20 0)\recs \move(-10
0)\recs

\move(0 0)\bsegment

\move(0 0)\lvec(10 0)\lvec(10 10)\lvec(0 10)\lvec(0 0) \move(0
5)\lvec(10 5)\htext(5 2.5){$1$}\htext(5 7.5){$1$}

\move(0 10)\lvec(10 10)\lvec(10 20)\lvec(0 20)\lvec(0 10)\move(0
10)\lvec(10 20)\htext(2.5 17.5){$2$}\htext(7.5 12.5){$0$}

\move(0 20)\lvec(10 20)\lvec(10 30)\lvec(0 30)\lvec(0 20)\move(0
25)\lvec(10 25)\htext(5 22.5){$1$}\htext(5 27.5){$1$}

\move(0 30)\lvec(10 30)\lvec(10 40)\lvec(0 40)\lvec(0 30)\move(0
30)\lvec(10 40)\htext(2.5 37.5){$2$}\htext(7.5 32.5){$0$}

\move(0 40)\lvec(10 40)\lvec(10 50)\lvec(0 50)\lvec(0 40)\move(0
45)\lvec(10 45)\htext(5 42.5){$1$}\htext(5 47.5){$1$}

\move(0 50)\lvec(10 50)\lvec(10 60)\lvec(0 60)\lvec(0 50)\move(0
50)\lvec(10 60)\htext(2.5 57.5){$2$}\htext(7.5 52.5){$0$}

\move(0 60)\lvec(10 60)\lvec(10 70)\lvec(0 70)\lvec(0 60)\move(0
65)\lvec(10 65)\htext(5 62.5){$1$}\htext(5 67.5){$1$}

\move(0 70)\lvec(10 70)\lvec(10 80)\lvec(0 80)\lvec(0 70)\move(0
70)\lvec(10 80)\htext(2.5 77.5){$2$}\htext(7.5 72.5){$0$}

\esegment

\move(10 0)\bsegment

\move(0 0)\lvec(10 0)\lvec(10 10)\lvec(0 10)\lvec(0 0) \move(0
5)\lvec(10 5)\htext(5 2.5){$1$}\htext(5 7.5){$1$}

\move(0 10)\lvec(10 10)\lvec(10 20)\lvec(0 20)\lvec(0 10)\move(0
10)\lvec(10 20)\htext(2.5 17.5){$0$}\htext(7.5 12.5){$2$}

\move(0 20)\lvec(10 20)\lvec(10 30)\lvec(0 30)\lvec(0 20)\move(0
25)\lvec(10 25)\htext(5 22.5){$1$}\htext(5 27.5){$1$}

\move(0 30)\lvec(10 30)\lvec(10 40)\lvec(0 40)\lvec(0 30)\move(0
30)\lvec(10 40)\htext(2.5 37.5){$0$}\htext(7.5 32.5){$2$}

\move(0 40)\lvec(10 40)\lvec(10 50)\lvec(0 50)\lvec(0 40)\move(0
45)\lvec(10 45)\htext(5 42.5){$1$}\htext(5 47.5){$1$}

\move(0 50)\lvec(10 50)\lvec(10 60)\lvec(0 60)\lvec(0 50)\move(0
50)\lvec(10 60)\htext(2.5 57.5){$0$}\htext(7.5 52.5){$2$}

\move(0 60)\lvec(10 60)\lvec(10 70)\lvec(0 70)\lvec(0 60)\move(0
65)\lvec(10 65)\htext(5 62.5){$1$}\htext(5 67.5){$1$}

\move(0 70)\lvec(10 70)\lvec(10 80)\lvec(0 80)\lvec(0 70)\move(0
70)\lvec(10 80)\htext(2.5 77.5){$0$}\htext(7.5 72.5){$2$}

\esegment

\move(-20 0)\bsegment

\move(0 0)\lvec(10 0)\lvec(10 10)\lvec(0 10)\lvec(0 0) \move(0
5)\lvec(10 5)\htext(5 2.5){$1$}\htext(5 7.5){$1$}

\move(0 10)\lvec(10 10)\lvec(10 20)\lvec(0 20)\lvec(0 10)\move(0
10)\lvec(10 20)\htext(2.5 17.5){$2$}\htext(7.5 12.5){$0$}

\move(0 20)\lvec(10 20)\lvec(10 30)\lvec(0 30)\lvec(0 20)\move(0
25)\lvec(10 25)\htext(5 22.5){$1$}\htext(5 27.5){$1$}

\move(0 30)\lvec(10 30)\lvec(10 40)\lvec(0 40)\lvec(0 30)\move(0
30)\lvec(10 40)\htext(2.5 37.5){$2$}\htext(7.5 32.5){$0$}

\move(0 40)\lvec(10 40)\lvec(10 50)\lvec(0 50)\lvec(0 40)\move(0
45)\lvec(10 45)\htext(5 42.5){$1$}\htext(5 47.5){$1$}

\move(0 50)\lvec(10 50)\lvec(10 60)\lvec(0 60)\lvec(0 50)\move(0
50)\lvec(10 60)\htext(2.5 57.5){$2$}\htext(7.5 52.5){$0$}

\move(0 60)\lvec(10 60)\lvec(10 70)\lvec(0 70)\lvec(0 60)\move(0
65)\lvec(10 65)\htext(5 62.5){$1$}\htext(5 67.5){$1$}

\move(0 70)\lvec(10 70)\lvec(10 80)\lvec(0 80)\lvec(0 70)\move(0
70)\lvec(10 80)\htext(2.5 77.5){$2$}\htext(7.5 72.5){$0$}

\esegment

\move(-10 0)\bsegment

\move(0 0)\lvec(10 0)\lvec(10 10)\lvec(0 10)\lvec(0 0) \move(0
5)\lvec(10 5)\htext(5 2.5){$1$}\htext(5 7.5){$1$}

\move(0 10)\lvec(10 10)\lvec(10 20)\lvec(0 20)\lvec(0 10)\move(0
10)\lvec(10 20)\htext(2.5 17.5){$0$}\htext(7.5 12.5){$2$}

\move(0 20)\lvec(10 20)\lvec(10 30)\lvec(0 30)\lvec(0 20)\move(0
25)\lvec(10 25)\htext(5 22.5){$1$}\htext(5 27.5){$1$}

\move(0 30)\lvec(10 30)\lvec(10 40)\lvec(0 40)\lvec(0 30)\move(0
30)\lvec(10 40)\htext(2.5 37.5){$0$}\htext(7.5 32.5){$2$}

\move(0 40)\lvec(10 40)\lvec(10 50)\lvec(0 50)\lvec(0 40)\move(0
45)\lvec(10 45)\htext(5 42.5){$1$}\htext(5 47.5){$1$}

\move(0 50)\lvec(10 50)\lvec(10 60)\lvec(0 60)\lvec(0 50)\move(0
50)\lvec(10 60)\htext(2.5 57.5){$0$}\htext(7.5 52.5){$2$}

\move(0 60)\lvec(10 60)\lvec(10 70)\lvec(0 70)\lvec(0 60)\move(0
65)\lvec(10 65)\htext(5 62.5){$1$}\htext(5 67.5){$1$}

\move(0 70)\lvec(10 70)\lvec(10 80)\lvec(0 80)\lvec(0 70)\move(0
70)\lvec(10 80)\htext(2.5 77.5){$0$}\htext(7.5 72.5){$2$}

\esegment

\move(20 0)\lvec(-25 0) \move(20 10)\lvec(-25 10) \move(20
20)\lvec(-25 20) \move(20 30)\lvec(-25 30) \move(20 40)\lvec(-25
40) \move(20 50)\lvec(-25 50) \move(20 60)\lvec(-25 60)
 \move(20 70)\lvec(-25 70) \move(20 80)\lvec(-25 80)

\move(20 0)\lvec(20 85) \move(10 0)\lvec(10 85)\move(0 0)\lvec(0
85)\move(-10 0)\lvec(-10 85)\move(-20 0)\lvec(-20 85)
\end{texdraw}
}
\end{center}\vskip 3mm
Here the shaded blocks in the above patterns are the ones in the
ground state walls.
\end{itemize}
\vskip 5mm

A wall built on $Y_{\Lambda}$ following the above rules is called
a {\it Young wall on $Y_{\Lambda}$}, for the heights of its
columns are weakly decreasing as we proceed from right to left.

\begin{df}{\rm \mbox{}

\begin{itemize}
\item[(1)] A column of a Young wall is called a {\it full
column} if its volume is of an integral value.

\item[(2)] A Young wall is said to be {\it proper} if none of the full columns
have the same heights.
\end{itemize}
}
\end{df}
We denote by ${\mathcal Z}(\Lambda)$ the set of all proper Young
walls on $Y_{\Lambda}$. For $Y\in{\mathcal Z}(\Lambda)$, we often
write $Y=(y_k)_{k=0}^{\infty}=(\cdots,y_2,y_1,y_0)$ as an infinite
sequence of its columns. Let $|y_k|$ be the number of blocks in
$y_k$ added to $Y_{\Lambda}$. Then the {\it associated partition}
is defined to be $|Y|=(|y_k|)_{k=0}^{\infty}$.

\begin{ex}\label{YW}{\rm
We illustrate several examples of proper Young walls. For
convenience, we omit the columns of the ground state wall on which
no block has been added.

\vskip 3mm
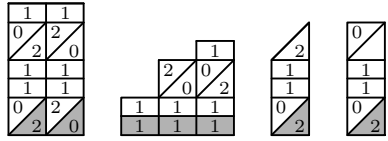
\begin{center}
\begin{texdraw}
\drawdim mm \setunitscale 0.5 \fontsize{7}{7}\selectfont \textref
h:C v:C \move(0 0)\tris \move(10 0)\tris \move(30 0)\recs \move(40
0)\recs \move(50 0)\recs \move(70 0)\tris \move(90 0)\tris

\move(0 0) \bsegment \move(0 0)\rlvec(20 0) \move(0 10)\rlvec(20
0) \move(0 15)\rlvec(20 0) \move(0 20)\rlvec(20 0) \move(0
30)\rlvec(20 0) \move(0 35)\rlvec(20 0) \move(0 35)\rlvec(20 0)
\move(0 0)\rlvec(0 35) \move(10 0)\rlvec(0 35) \move(20 0)\rlvec(0
35) \move(0 0)\rlvec(10 10) \move(10 0)\rlvec(10 10) \move(0
20)\rlvec(10 10) \move(10 20)\rlvec(10 10) \htext(2.5 7.5){$0$}
\htext(7.5 2.5){$2$} \htext(12.5 7.5){$2$} \htext(17.5 2.5){$0$}
\htext(2.5 27.5){$0$} \htext(7.5 22.5){$2$} \htext(12.5 27.5){$2$}
\htext(17.5 22.5){$0$} \htext(5 12.5){$1$} \htext(15 12.5){$1$}
\htext(5 17.5){$1$} \htext(15 17.5){$1$} \htext(5 32.5){$1$}
\htext(15 32.5){$1$} \esegment

\move(60 0) \bsegment \move(0 0)\rlvec(-30 0) \move(0 5)\rlvec(-30
0) \move(0 10)\rlvec(-20 0) \move(0 20)\rlvec(-20 0) \move(0
25)\rlvec(-10 0) \move(0 0)\rlvec(0 25) \move(-10 0)\rlvec(0 25)
\move(-20 0)\rlvec(0 20) \move(-30 0)\rlvec(0 5) \move(-10
10)\rlvec(10 10) \move(-20 10)\rlvec(10 10) \move(-30 5)\lvec(-30
10)\lvec(-20 10) \htext(-5 2.5){$1$} \htext(-15 2.5){$1$}
\htext(-5 7.5){$1$} \htext(-15 7.5){$1$} \htext(-25 2.5){$1$}
\htext(-25 7.5){$1$} \htext(-5 22.5){$1$} \htext(-7.5 17.5){$0$}
\htext(-2.5 12.5){$2$} \htext(-17.5 17.5){$2$} \htext(-12.5
12.5){$0$} \esegment

\move(80 0) \bsegment \move(0 0)\lvec(0 30)\lvec(-10 20)\lvec(-10
0)\lvec(0 0) \move(-10 0)\lvec(0 10)\lvec(-10 10) \move(0
15)\lvec(-10 15) \move(0 20)\lvec(-10 20)\lvec(0 30) \htext(-2.5
2.5){$2$} \htext(-7.5 7.5){$0$} \htext(-5 12.5){$1$} \htext(-5
17.5){$1$} \htext(-2.5 22.5){$2$} \esegment

\move(100 0) \bsegment \move(0 0)\lvec(0 30)\lvec(-10 30)\lvec(-10
0)\lvec(0 0) \move(-10 0)\lvec(0 10)\lvec(-10 10) \move(0
15)\lvec(-10 15) \move(0 20)\lvec(-10 20)\lvec(0 30) \htext(-2.5
2.5){$2$} \htext(-7.5 7.5){$0$} \htext(-5 12.5){$1$} \htext(-5
17.5){$1$} \htext(-7.5 27.5){$0$} \esegment
\end{texdraw}

\end{center}
\vskip 3mm }
\end{ex}

\begin{df}\label {admrmv}{\rm Let $Y$ be a proper Young wall on $Y_{\Lambda}$.

\begin{enumerate}
\item[(1)] A block of color $i$ (in short, an $i$-block) in $Y$ is called a {\it removable
$i$-block} if $Y$ remains a proper Young wall after removing the
block.
\item[(2)] A place in $Y$ is called an {\it admissible $i$-slot} if one
may add an $i$-block to obtain another proper Young wall.
\item[(3)] A column in $Y$ is said to be {\it $i$-removable {\rm(}resp.
$i$-admissible{\rm )}} if there is a removable $i$-block {\rm
(}resp. an admissible $i$-slot{\rm )} in that column.
\end{enumerate}
}
\end{df}
We now define the {\it abstract Kashiwara operators}
$\tilde{E}_i$, $\tilde{F}_i$ on ${\mathcal Z}(\Lambda)$ as
follows. Fix $i\in I$ and let $Y=(y_k)_{k=0}^{\infty}$ be a proper
Young wall on $Y_{\Lambda}$.
\begin{itemize}
\item[(1)] To each column $y_k$ of $Y$, we assign
\begin{equation*}
\begin{cases}
-- & \text{if $y_k$ is twice $i$-removable,} \\
- & \text{if $y_k$ is once $i$-removable but not $i$-admissible,} \\
-+ & \text{if $y_k$ is once $i$-removable and once
$i$-admissible,} \\
+ & \text{if $y_k$ is once $i$-admissible but not $i$-removable,}
\\
++ & \text{if $y_k$ is twice $i$-admissible,} \\
\ \ \cdot & \text{otherwise}.
\end{cases}
\end{equation*}

\item[(2)] From this infinite sequence of $\pm$'s and $\cdot$'s, we cancel
out every $(+,-)$-pair to obtain a finite sequence of $-$'s
followed by $+$'s, reading from left to right. This finite
sequence $(-\cdots-,+\cdots+)$ is called the {\it $i$-signature}
of $Y$.

\item[(3)] We define $\tilde{E}_i Y$ to be the proper Young wall
obtained from $Y$ by removing the $i$-block corresponding to the
right-most $-$ in the $i$-signature of $Y$. We define $\tilde{E}_i
Y=0$ if there is no $-$ in the $i$-signature of $Y$.

\item[(4)] We define $\tilde{F}_i Y$ to be the proper Young wall
obtained from $Y$ by adding an $i$-block to the column
corresponding to the left-most $+$ in the $i$-signature of $Y$. We
define $\tilde{F}_i Y=0$ if there is no $+$ in the $i$-signature
of $Y$.
\end{itemize}
Next, we define
\begin{equation*}
\begin{split}
{\rm wt}(Y)&=\Lambda-\sum_{i\in I}k_i\alpha_i \in P, \\
\varepsilon_i(Y)&=\text{the number of $-$'s in the $i$-signature of $Y$}, \\
\varphi_i(Y)&=\text{the number of $+$'s in the $i$-signature of
$Y$},
\end{split}
\end{equation*}
where $k_i$ denotes the number of $i$-blocks in $Y$ which have
been added to $Y_{\Lambda}$.
\begin{prop}\label {HK}{\rm \cite{HK}}
The set ${\mathcal Z}(\Lambda)$ together with the maps ${\rm wt}$,
$\varepsilon_i$, $\varphi_i$, $\tilde{E}_i$ and $\tilde{F}_i$ {\rm
($i\in I$)} becomes an affine crystal.\qed
\end{prop}
The part of a column with $a_i$-many $i$-blocks for each $i\in I$
($a_0=a_2=1,a_1=2$) in some cyclic order is called a {\it
$\delta$-column}. A $\delta$-column in a proper Young wall is {\it
removable} if it can be removed to yield another proper Young
wall.
\begin{df}{\rm
A proper Young wall $Y$ is said to be {\it reduced} if none of its
columns contain a removable $\delta$-column.}
\end{df}

\begin{ex}{\rm
Among the proper Young walls given in Example \ref{YW}, the second
and the fourth ones are reduced, but the others are not.

}
\end{ex}

Let ${\mathcal Y}(\Lambda)\subset{\mathcal Z}(\Lambda) $ be the
set of all reduced proper Young walls on $Y_{\Lambda}$. Then we
have
\begin{thm}{\rm \cite{HK}}
The set ${\mathcal Y}(\Lambda)$ is an affine crystal. Moreover,
there exists an affine crystal isomorphism ${\mathcal Y}(\Lambda)
\stackrel{\sim}{\longrightarrow} B(\Lambda)$, where $B(\Lambda)$
is the crystal of the basic representation $V(\Lambda)$. \qed
\end{thm}

Let $Y=(y_k)_{k=0}^{\infty}$ be a proper Young wall in ${\mathcal
Z}(\Lambda)$. Let $S$ be an interval in $\mathbb{Z}_{\geq 0}$
which is finite or infinite; i.e. $S=\{\,k\,|\,s\leq k < t\,\}$
for some $0\leq s < t\leq\infty$. We call $Y'=(y_k)_{k\in S}$ a
{\it part of $Y$}. If $S$ is infinite; that is, if
$Y'=(y_k)_{k=s}^{\infty}$, then $Y'$ is itself a proper Young wall
in ${\mathcal Y}(\Lambda')$ for some $\Lambda'$. If $S$ is finite;
that is, if $Y'=(y_{t-1},\cdots,y_s)$, then $Y'$ is not a proper
Young wall, but a finite collection of successive columns in $Y$.
Also, by adding or removing blocks only in columns $y_s$ ($s\in
S$) of $Y'$, we can extend the notions of admissible $i$-slots,
removable $i$-blocks, the $i$-signature, $\varepsilon_i$, and
$\varphi_i$ of a part $Y'$ (however, we define {\rm wt} only for
proper Young walls). The notion of parts will be used when we
define the action of $U_q(C_2^{(1)})$ on ${\mathcal Z}(\Lambda)$.

\begin{ex}\label{graph}{\rm
\samepage Crystal graph $B(\Lambda_1)$.}
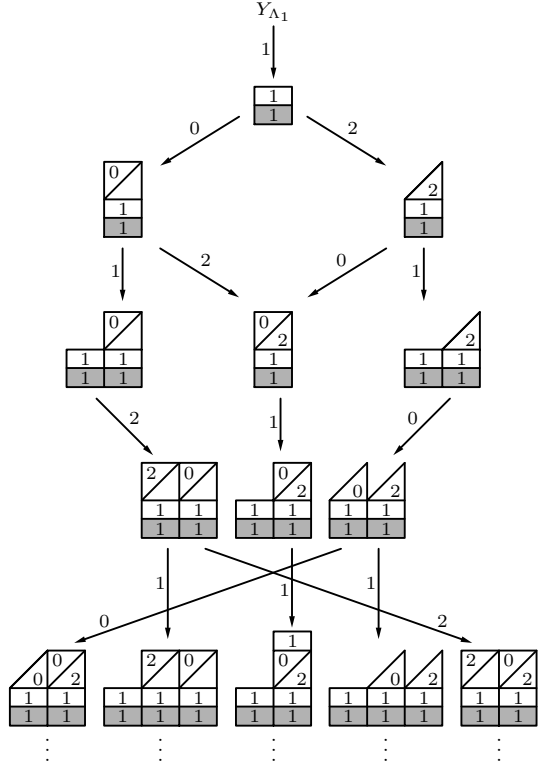
\begin{center}
\begin{texdraw}
\drawdim mm \fontsize{7}{7}\selectfont \textref h:C v:C
\setunitscale 0.5 \htext(5 30){$Y_{\Lambda_1}$} \move(0 33)
\move(0 0)\recs \move(0 0) \bsegment \move(0 0)\lvec(10 0)\lvec(10
10)\lvec(0 10)\lvec(0 0) \move(10 5)\lvec(0 5) \htext(5
2.5){$1$}\htext(5 7.5){$1$} \esegment
\move(-40 -30)\recs \move(-40 -30) \bsegment \move(0 0)\lvec(10
0)\lvec(10 20)\lvec(0 20)\lvec(0 0) \move(0 5)\lvec(10 5) \move(10
10)\lvec(0 10)\lvec(10 20) \htext(5 2.5){$1$}\htext(5
7.5){$1$}\htext(2.5 17.5){$0$} \esegment
\move(40 -30)\recs \move(40 -30) \bsegment \move(0 0)\lvec(10
0)\lvec(10 20)\lvec(0 10)\lvec(0 0) \move(0 5)\lvec(10 5) \move(10
10)\lvec(0 10)\lvec(10 20) \htext(5 2.5){$1$}\htext(5
7.5){$1$}\htext(7.5 12.5){$2$} \esegment
\move(-40 -70)\recs \move(-50 -70)\recs \move(-40 -70) \bsegment
\move(10 10)\lvec(-10 10)\lvec(-10 0)\lvec(10 0)\lvec(10
20)\lvec(0 20) \lvec(0 0) \move(-10 5)\lvec(10 5) \move(0
10)\lvec(10 20) \htext(5 2.5){$1$}\htext(5 7.5){$1$}\htext(2.5
17.5){$0$} \htext(-5 2.5){$1$}\htext(-5 7.5){$1$} \esegment
\move(0 -70)\recs \move(0 -70) \bsegment \move(0 0)\lvec(10
0)\lvec(10 20)\lvec(0 20)\lvec(0 0) \move(0 5)\lvec(10 5) \move(10
10)\lvec(0 10)\lvec(10 20) \htext(5 2.5){$1$}\htext(5
7.5){$1$}\htext(7.5 12.5){$2$}\htext(2.5 17.5){$0$} \esegment
\move(50 -70)\recs \move(40 -70)\recs \move(50 -70) \bsegment
\move(10 10)\lvec(-10 10)\lvec(-10 0)\lvec(10 0)\lvec(10
20)\lvec(0 10) \lvec(0 0) \move(-10 5)\lvec(10 5) \move(0
10)\lvec(10 20) \htext(5 2.5){$1$}\htext(5 7.5){$1$}\htext(7.5
12.5){$2$} \htext(-5 2.5){$1$}\htext(-5 7.5){$1$} \esegment
\move(-20 -110)\recs \move(-30 -110)\recs \move(-20 -110)
\bsegment \move(0 20)\lvec(0 0)\lvec(10 0)\lvec(10 20)\lvec(-10
20)\lvec(-10 0)\lvec(0 0) \move(-10 5)\lvec(10 5) \move(10
10)\lvec(0 10)\lvec(10 20) \move(0 10)\lvec(-10 10)\lvec(0 20)
\htext(5 2.5){$1$}\htext(5 7.5){$1$} \htext(-5 2.5){$1$}\htext(-5
7.5){$1$} \htext(-7.5 17.5){$2$}\htext(2.5 17.5){$0$} \esegment
\move(5 -110)\recs \move(-5 -110)\recs \move(5 -110) \bsegment
\move(10 10)\lvec(-10 10)\lvec(-10 0)\lvec(10 0)\lvec(10
20)\lvec(0 20) \lvec(0 0) \move(-10 5)\lvec(10 5) \move(0
10)\lvec(10 20) \htext(5 2.5){$1$}\htext(5 7.5){$1$}\htext(2.5
17.5){$0$}\htext(7.5 12.5){$2$} \htext(-5 2.5){$1$}\htext(-5
7.5){$1$} \esegment
\move(30 -110)\recs \move(20 -110)\recs \move(30 -110) \bsegment
\move(0 20)\lvec(0 0)\lvec(10 0)\lvec(10 20)\move(-10 10)\lvec(-10
0)\lvec(0 0) \move(-10 5)\lvec(10 5) \move(10 10)\lvec(0
10)\lvec(10 20) \move(0 10)\lvec(-10 10)\lvec(0 20) \htext(5
2.5){$1$}\htext(5 7.5){$1$} \htext(-5 2.5){$1$}\htext(-5 7.5){$1$}
\htext(-2.5 12.5){$0$}\htext(7.5 12.5){$2$} \esegment
\move(-55 -160)\recs \move(-65 -160)\recs \move(-55 -160)
\bsegment \move(0 20)\lvec(0 0)\lvec(10 0)\lvec(10 20)\lvec(0
20)\lvec(-10 10)\lvec(-10 0)\lvec(0 0) \move(-10 5)\lvec(10 5)
\move(10 10)\lvec(0 10)\lvec(10 20) \move(0 10)\lvec(-10
10)\lvec(0 20) \htext(5 2.5){$1$}\htext(5 7.5){$1$} \htext(-5
2.5){$1$}\htext(-5 7.5){$1$} \htext(-2.5 12.5){$0$}\htext(7.5
12.5){$2$}\htext(2.5 17.5){$0$} \vtext(0 -6){$\cdots$} \esegment
\move(-20 -160)\recs \move(-30 -160)\recs \move(-40 -160)\recs
\move(-20 -160) \bsegment \move(0 20)\lvec(0 0)\lvec(10 0)\lvec(10
20)\lvec(-10 20)\lvec(-10 0)\lvec(0 0) \move(-10 5)\lvec(10 5)
\move(10 10)\lvec(0 10)\lvec(10 20) \move(0 10)\lvec(-10
10)\lvec(0 20) \htext(5 2.5){$1$}\htext(5 7.5){$1$} \htext(-5
2.5){$1$}\htext(-5 7.5){$1$} \htext(-7.5 17.5){$2$}\htext(2.5
17.5){$0$} \move(-10 0)\lvec(-20 0)\lvec(-20 10)\lvec(-10 10)
\move(-10 5)\lvec(-20 5) \htext(-15 2.5){$1$} \htext(-15 7.5){$1$}
\vtext(-5 -6){$\cdots$} \esegment
\move(5 -160)\recs \move(-5 -160)\recs \move(5 -160) \bsegment
\move(10 10)\lvec(-10 10)\lvec(-10 0)\lvec(10 0)\lvec(10
20)\lvec(0 20) \lvec(0 0) \move(-10 5)\lvec(10 5) \move(0
10)\lvec(10 20) \htext(5 2.5){$1$}\htext(5 7.5){$1$}\htext(2.5
17.5){$0$}\htext(7.5 12.5){$2$} \htext(-5 2.5){$1$}\htext(-5
7.5){$1$} \move(10 20)\lvec(10 25)\lvec(0 25)\lvec(0 20) \htext(5
22.5){$1$} \vtext(0 -6){$\cdots$} \esegment
\move(40 -160)\recs \move(30 -160)\recs \move(20 -160)\recs
\move(40 -160) \bsegment \move(0 20)\lvec(0 0)\lvec(10 0)\lvec(10
20)\move(-10 10)\lvec(-10 0)\lvec(0 0) \move(-10 5)\lvec(10 5)
\move(10 10)\lvec(0 10)\lvec(10 20) \move(0 10)\lvec(-10
10)\lvec(0 20) \htext(5 2.5){$1$}\htext(5 7.5){$1$} \htext(-5
2.5){$1$}\htext(-5 7.5){$1$} \htext(-2.5 12.5){$0$}\htext(7.5
12.5){$2$} \move(-10 0)\lvec(-20 0)\lvec(-20 10)\lvec(-10 10)
\move(-10 5)\lvec(-20 5) \htext(-15 2.5){$1$} \htext(-15 7.5){$1$}
\vtext(-5 -6){$\cdots$} \esegment
\move(65 -160)\recs \move(55 -160)\recs \move(65 -160) \bsegment
\move(0 20)\lvec(0 0)\lvec(10 0)\lvec(10 20)\lvec(-10 20)\lvec(-10
0)\lvec(0 0) \move(-10 5)\lvec(10 5) \move(10 10)\lvec(0
10)\lvec(10 20) \move(0 10)\lvec(-10 10)\lvec(0 20) \htext(5
2.5){$1$}\htext(5 7.5){$1$} \htext(-5 2.5){$1$}\htext(-5 7.5){$1$}
\htext(-7.5 17.5){$2$}\htext(2.5 17.5){$0$} \htext(7.5 12.5){$2$}
\vtext(0 -6){$\cdots$} \esegment
\move(0 0) \bsegment \linewd 0.45 \arrowheadsize l:3 w:1.5
\arrowheadtype t:F \move(5 26)\ravec(0 -14)\htext(3 20){$1$}
\move(-4 2)\ravec(-21 -13)\htext(-16 -1){$0$} \move(14 2)\ravec(21
-13)\htext(26 -1){$2$} \move(-25 -33)\ravec(21 -13)\htext(-13
-36){$2$} \move(35 -33)\ravec(-21 -13)\htext(23 -36){$0$}
\move(-35 -33)\ravec(0 -14)\htext(-37 -39){$1$} \move(45
-33)\ravec(0 -14)\htext(43 -39){$1$} \move(7 -73)\ravec(0
-14)\htext(5 -79){$1$} \move(52 -73)\ravec(-15 -15)\htext(42
-78){$0$} \move(-42 -73)\ravec(15 -15)\htext(-32 -78){$2$}
\move(23 -113)\ravec(-71 -25)\htext(-40 -132){$0$} \move(-13
-113)\ravec(71 -25)\htext(50 -132){$2$} \move(33 -113)\ravec(0
-24)\htext(31 -122){$1$} \move(-23 -113)\ravec(0 -24)\htext(-25
-122){$1$} \move(10 -113)\ravec(0 -21)\htext(8 -124){$1$}
\esegment \move(0 -170)
\end{texdraw}
\end{center}
\end{ex}

\section{Fock space representation}
Given a dominant integral weight $\Lambda=\Lambda_i$ ($i\in I$),
we define ${\mathcal F}(\Lambda)=\bigoplus_{Y\in {\mathcal
Z}(\Lambda)}\mathbb{Q}(q)Y$ to be the $\mathbb{Q}(q)$-vector space
with a basis ${\mathcal Z}(\Lambda)$. In this section, we will
define a $U_q(C_2^{(1)})$-module structure on ${\mathcal
F}(\Lambda)$, the {\it Fock space representation of
$U_q(C_2^{(1)})$}. Then we will show that the abstract affine
crystal ${\mathcal Z}(\Lambda)$ is isomorphic to the crystal of
${\mathcal F}(\Lambda)$.

Let us define the action of $U_q(C_2^{(1)})$ on ${\mathcal
Z}(\Lambda)$. For $Y=(y_k)_{k=0}^{\infty}\in {\mathcal
Z}(\Lambda)$ and $q^h$ ($h\in P^{\vee}$), we define
\begin{equation}
q^h Y=q^{\langle h,{\rm wt}(Y) \rangle}Y.
\end{equation}
The actions of $e_i$ and $f_i$ ($i\in I$) are given according to
the type of the $i$-block.\vskip 3mm

{\bf Case 1.} Suppose that $i=1$ ($q=q_1$).

Let $b$ be a removable $1$-block in $y_k$ of $Y$. If the
$1$-signature of $y_k$ is $--$, or if the $1$-signature of $y_k$
is $-$ and there is another $1$-block below $b$, we define
$Y\nearrow b$ to be the Young wall obtained by removing $b$ from
$Y$. If the $1$-signature of $y_k$ is $-+$, or if the
$1$-signature of $y_k$ is $-$ and there is no $1$-block below $b$,
we define
\begin{equation*}
Y\nearrow b=q^{-1}(1-(-q^2)^{l(b)+1})Z,
\end{equation*}
where $Z$ is the Young wall obtained by removing $b$ from $Y$ and
$l(b)$ is the number of $y_l$'s with $l<k$ such that
$|y_l|=|y_k|$. That is,if

\vskip 5mm \hskip 4cm $Y=$\raisebox{-0.5\height}{
\begin{texdraw}\fontsize{9}{9}
\drawdim em \setunitscale 0.1 \linewd 0.5

\move(-20 -20)\lvec(-20 0)\lvec(0 0)\lvec(0 10)\lvec(50
10)\lvec(50 30)\lvec(60 30)\lvec(60 40)\lvec(80 40)\lvec(80
-20)\lvec(-20 -20)

\move(50 10)\lvec(50 -20)

\htext(57 10){$Y_R(b)$}

\move(0 10) \lvec(0 15) \lvec(10 15)\lvec(10 10)\lvec(0 10)\lfill
f:0.8

\move(10 10)\lvec(10 15)\lvec(50 15)\lvec(50 10)\lvec(10 10)
\lfill f:0.8

\move(20 10) \lvec(20 15) \move(40 10) \lvec(40 15) \htext(25
11){$\cdots$}

\move(10 0) \arrowheadtype t:F \arrowheadsize l:4 w:2 \avec(3 12)

\htext(12 -1){$_b$}

\move(10 15)\clvec(12 20)(23 20)(25 20)

\move(50 15)\clvec(48 20)(37 20)(35 20)

\htext(26 18){$_{l(b)}$}
\end{texdraw}}\hskip 5mm,  \vskip 5mm

\hskip .9cm then \hskip 5mm$Y\nearrow b=$\raisebox{-0.5\height}{
\begin{texdraw}\fontsize{9}{9}
\drawdim em \setunitscale 0.1 \linewd 0.5

\move(-20 -20)\lvec(-20 0)\lvec(0 0)\lvec(0 10)\lvec(50
10)\lvec(50 30)\lvec(60 30)\lvec(60 40)\lvec(80 40)\lvec(80
-20)\lvec(-20 -20)

\move(10 10)\lvec(10 15)\lvec(50 15)\lvec(50 10)\lvec(10 10)
\lfill f:0.8

\move(20 10) \lvec(20 15) \move(40 10) \lvec(40 15) \htext(25
11){$\cdots$}

\htext(-73 12){$\dfrac{(1-(-q^2)^{l(b)+1})}{q}\times$}
\end{texdraw}}\hskip 5mm .
\vskip 3mm

In either case, if $k\geq 1$, let $Y_R(b)=(y_{l},\cdots,y_{0})$ be
the part of $Y$ with finite columns such that $l$ is the integer
satisfying $|y_k|=|y_{k-1}|=\cdots=|y_{l+1}|<|y_l|$. Set
$R_1(b;Y)=\varphi_1(Y_R(b))-\varepsilon_1(Y_R(b))$ if $k\geq 1$
and $0$ if $k=0$. Then we define
\begin{equation}
e_1\,Y=\sum_{b}q^{-R_1(b;Y)}(Y\nearrow b),
\end{equation}
\vskip -2mm \noindent where $b$ runs over all removable $1$-blocks
in $Y$.

On the other hand, suppose that $b$ is an admissible $1$-slot in
$y_k$ of $Y$. If the $1$-signature of $y_k$ is $++$, or if the
$1$-signature of $y_k$ is $+$ and there is no $1$-block below $b$,
we define $Y\swarrow b$ to be the Young wall obtained by adding a
$1$-block at $b$. If the $1$-signature of $y_k$ is $-+$, or if the
$1$-signature of $y_k$ is $+$ and there is another $1$-block below
$b$, then we define
\begin{equation*}
Y\swarrow b=q^{-1}(1-(-q^2)^{l(b)+1})Z,
\end{equation*}
where $Z$ is the Young wall obtained by adding a $1$-block at $b$
and $l(b)$ is the number of $y_l$'s with $l>k$ such that
$|y_l|=|y_k|$. That is, if \vskip 5mm \hskip 4cm
$Y=$\raisebox{-0.5\height}{
\begin{texdraw}\fontsize{9}{9}
\drawdim em \setunitscale 0.1 \linewd 0.5

\move(-30 -20)\lvec(-30 0)\lvec(-10 0)\lvec(-10 10)\lvec(50
10)\lvec(50 30)\lvec(65 30)\lvec(65 -20)\lvec(-30 -20)

\move(0 10)\lvec(0 -20)\htext(-24 -14){$Y_L(b)$}

\move(40 15) \lvec(40 20)\lvec(50 20)

\move(0 10) \lvec(0 15) \lvec(10 15)\lvec(10 10)\lvec(0 10)\lfill
f:0.8

\move(10 10)\lvec(10 15)\lvec(50 15)\lvec(50 10)\lvec(10 10)
\lfill f:0.8

\move(30 10) \lvec(30 15) \move(40 10) \lvec(40 15) \htext(15
11){$\cdots$}

\move(55 0) \arrowheadtype t:F \arrowheadsize l:4 w:2 \avec(43 17)

\htext(60 -3){$_b$}

\move(0 15)\clvec(2 20)(13 20)(15 20)

\move(40 15)\clvec(38 20)(27 20)(25 20)

\htext(16 18){$_{l(b)}$}
\end{texdraw}}\hskip 5mm,
\vskip 0.5cm

\hskip .7cm then \hskip 5mm$Y\swarrow b=$\raisebox{-0.5\height}{
\begin{texdraw}\fontsize{9}{9}
\drawdim em \setunitscale 0.1 \linewd 0.5

\move(-30 -20)\lvec(-30 0)\lvec(-10 0)\lvec(-10 10)\lvec(50
10)\lvec(50 30)\lvec(65 30)\lvec(65 -20)\lvec(-30 -20)

\move(40 15) \lvec(40 20)\lvec(50 20)\lvec(50 15)\lfill f:0.8

\move(0 10) \lvec(0 15) \lvec(10 15)\lvec(10 10)\lvec(0 10)\lfill
f:0.8

\move(10 10)\lvec(10 15)\lvec(50 15)\lvec(50 10)\lvec(10 10)
\lfill f:0.8

\move(30 10) \lvec(30 15) \move(40 10) \lvec(40 15) \htext(15
11){$\cdots$}

\htext(-90 0){$\dfrac{(1-(-q^2)^{l(b)+1})}{q}\times$}
\end{texdraw}}\hskip 5mm .\vskip 5mm

In either case, let $Y_L(b)=(\cdots,y_{l+2},y_{l+1})$, where $l$
is the integer such that $|y_{l+1}|<|y_l|=|y_{l-1}|=\cdots=|y_k|$,
and set $L_1(b;Y)=\varphi_1(Y_L(b))-\varepsilon_1(Y_L(b))$. Then
we define
\begin{equation}
f_1\,Y=\sum_{b}q^{L_1(b;Y)}(Y\swarrow b),
\end{equation}
\vskip -2mm \noindent where $b$ runs over all admissible $1$-slots
in $Y$.\vskip 3mm

{\bf Case 2.} Suppose that $i=0,2$ ($q^2=q_i$).

If $b$ is a removable $i$-block in $y_k$ of $Y$, then we define
$Y\nearrow b$ to be the Young wall obtained by removing $b$ from
$Y$. Consider the following $i$-block $b$ in $y_k$ of $Y$, called
a {\it virtually removable $i$-block}.\vskip 3mm

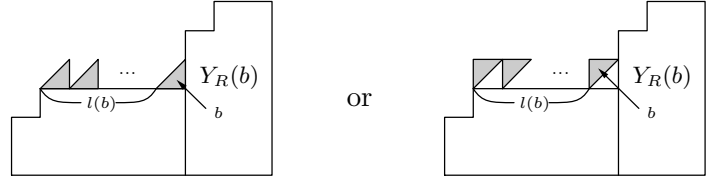
\begin{center}
\raisebox{-0.4\height}{\begin{texdraw}\fontsize{9}{9} \drawdim em
\setunitscale 0.1 \linewd 0.5

\move(0 -20)\lvec(0 0)\lvec(10 0)\lvec(10 10)\lvec(60 10)\lvec(60
30)\lvec(70 30)\lvec(70 40)\lvec(90 40)\lvec(90 -20)\lvec(0 -20)

\move(60 10)\lvec(60 -20)

\move(10 10)\lvec(20 20)\lvec(20 10)\lvec(10 10)\lfill f:0.8

\move(20 10)\lvec(30 20)\lvec(30 10)\lvec(20 10)\lfill f:0.8

\move(50 10)\lvec(60 20)\lvec(60 10)\lvec(50 10)\lfill f:0.8

\htext(65 10){$Y_R(b)$}

\htext(37 14){$_{\cdots}$}

\move(67 3) \arrowheadtype t:F \arrowheadsize l:4 w:2 \avec(57 13)

\htext(70 0){$_b$}

\move(24 5)\clvec(16 5)(13 5)(10 10)

\move(36 5)\clvec(44 5)(47 5)(50 10)

\htext(26 2){$_{l(b)}$}
\end{texdraw}}\hskip 1cm or
\hskip .7cm\raisebox{-0.4\height}{
\begin{texdraw}\fontsize{9}{9}
\drawdim em \setunitscale 0.1 \linewd 0.5

\move(0 -20)\lvec(0 0)\lvec(10 0)\lvec(10 10)\lvec(60 10)\lvec(60
30)\lvec(70 30)\lvec(70 40)\lvec(90 40)\lvec(90 -20)\lvec(0 -20)

\move(60 10)\lvec(60 -20)

\htext(65 10){$Y_R(b)$}

\move(10 10)\lvec(20 20)\lvec(10 20)\lvec(10 10)\lfill f:0.8

\move(20 10)\lvec(30 20)\lvec(20 20)\lvec(20 10)\lfill f:0.8

\move(50 10)\lvec(60 20)\lvec(50 20)\lvec(50 10)\lfill f:0.8


\htext(37 14){$_{\cdots}$}

\move(67 3) \arrowheadtype t:F \arrowheadsize l:4 w:2 \avec(53 17)

\htext(70 0){$_b$}

\move(24 5)\clvec(16 5)(13 5)(10 10)

\move(36 5)\clvec(44 5)(47 5)(50 10)

\htext(26 2){$_{l(b)}$}
\end{texdraw}}
\end{center}\vskip 3mm

\noindent In this case, we define $Y\nearrow b$ to be \vskip 3mm

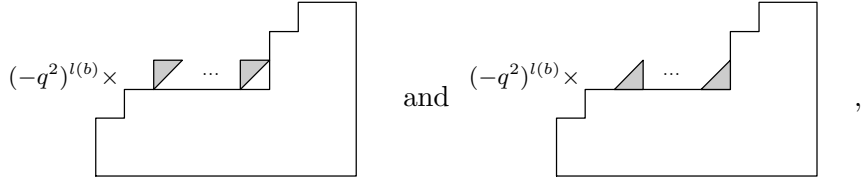
\begin{center}
\hskip -10mm\raisebox{-0.4\height}{
\begin{texdraw}\fontsize{9}{9}
\drawdim em \setunitscale 0.1 \linewd 0.5

\htext(-30 10){${(-q^2)^{l(b)}\times}$}

\move(0 -20)\lvec(0 0)\lvec(10 0)\lvec(10 10)\lvec(60 10)\lvec(60
30)\lvec(70 30)\lvec(70 40)\lvec(90 40)\lvec(90 -20)\lvec(0 -20)

\move(20 10)\lvec(30 20)\lvec(20 20)\lvec(20 10)\lfill f:0.8

\move(50 10)\lvec(60 20)\lvec(50 20)\lvec(50 10)\lfill f:0.8


\htext(37 14){$_{\cdots}$}
\end{texdraw}} \hskip 5mm and
\raisebox{-0.4\height}{
\begin{texdraw}\fontsize{9}{9}
\drawdim em \setunitscale 0.1 \linewd 0.5

\htext(-30 10){${(-q^2)^{l(b)}\times}$}

\move(0 -20)\lvec(0 0)\lvec(10 0)\lvec(10 10)\lvec(60 10)\lvec(60
30)\lvec(70 30)\lvec(70 40)\lvec(90 40)\lvec(90 -20)\lvec(0 -20)

\move(20 10)\lvec(30 20)\lvec(30 10)\lvec(20 10)\lfill f:0.8

\move(50 10)\lvec(60 20)\lvec(60 10)\lvec(50 10)\lfill f:0.8


\htext(37 14){$_{\cdots}$}

\end{texdraw}}\hskip 5mm ,
\end{center}\vskip 3mm

\noindent respectively, where $l(b)\geq 1$ is given in the above
figure. In either case, if $k\geq 1$, let
$Y_R(b)=(y_{k-1},\cdots,y_{0})$. Set
$R_i(b;Y)=\varphi_i(Y_R(b))-\varepsilon_i(Y_R(b))$ if $k\geq 1$,
and $0$ if $k=0$. Then we define
\begin{equation}
e_i\,Y=\sum_{b}q^{-2R_i(b;Y)}(Y\nearrow b),
\end{equation}
\vskip -2mm \noindent where $b$ runs over all removable and
virtually removable $i$-blocks in $Y$.

On the other hand, if $b$ is an admissible $i$-slot in $y_k$ of
$Y$, then we define $Y\swarrow b$ to be the Young wall obtained by
adding an $i$-block at $b$. Consider the following $i$-slot $b$ in
$y_k$ of $Y$, called a {\it virtually admissible $i$-slot}.\vskip
3mm

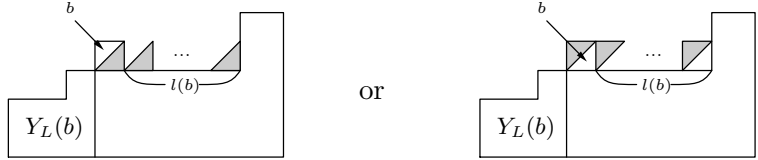
\begin{center}
\raisebox{-0.4\height}{
\begin{texdraw}\fontsize{9}{9}
\drawdim em \setunitscale 0.1 \linewd 0.5

\move(-20 -20)\lvec(-20 0)\lvec(0 0)\lvec(0 10)\lvec(60
10)\lvec(60 30)\lvec(75 30)\lvec(75 -20)\lvec(-20 -20)

\move(10 10)\lvec(10 -20)\htext(-14 -14){$Y_L(b)$}

\move(10 10)\lvec(10 20)\lvec(20 20)

\move(10 10)\lvec(20 20)\lvec(20 10)\lvec(10 10)\lfill f:0.8

\move(20 10)\lvec(30 20)\lvec(30 10)\lvec(20 10)\lfill f:0.8

\move(50 10)\lvec(60 20)\lvec(60 10)\lvec(50 10)\lfill f:0.8


\htext(37 14){$_{\cdots}$}

\move(3 27) \arrowheadtype t:F \arrowheadsize l:4 w:2 \avec(13 17)

\htext(0 30){$_b$}

\move(34 5)\clvec(26 5)(23 5)(20 10)

\move(46 5)\clvec(54 5)(57 5)(60 10)

\htext(36 2){$_{l(b)}$}
\end{texdraw}}\hskip 1cm or
\hskip 1cm\raisebox{-0.4\height}{
\begin{texdraw}\fontsize{9}{9}
\drawdim em \setunitscale 0.1 \linewd 0.5

\move(-20 -20)\lvec(-20 0)\lvec(0 0)\lvec(0 10)\lvec(60
10)\lvec(60 30)\lvec(75 30)\lvec(75 -20)\lvec(-20 -20)

\move(10 10)\lvec(10 -20)\htext(-14 -14){$Y_L(b)$}

\move(10 10)\lvec(20 20)\lvec(10 20)\lvec(10 10)\lfill f:0.8

\move(20 10)\lvec(30 20)\lvec(20 20)\lvec(20 10)\lfill f:0.8

\move(50 10)\lvec(60 20)\lvec(50 20)\lvec(50 10)\lfill f:0.8

\htext(37 14){$_{\cdots}$}

\move(3 27) \arrowheadtype t:F \arrowheadsize l:4 w:2 \avec(17 13)

\htext(0 30){$_b$}

\move(34 5)\clvec(26 5)(23 5)(20 10)

\move(46 5)\clvec(54 5)(57 5)(60 10)

\htext(36 2){$_{l(b)}$}
\end{texdraw}}
\end{center}\vskip 3mm

\noindent In this case, we define $Y\swarrow b$ to be \vskip 3mm

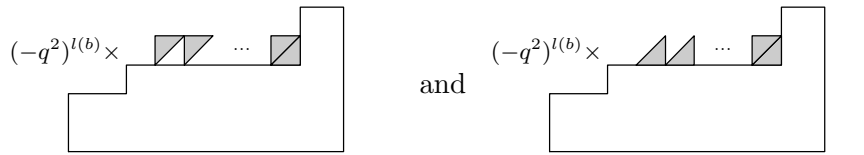
\begin{center}
\hskip -6mm\raisebox{-0.4\height}{
\begin{texdraw}\fontsize{9}{9}
\drawdim em \setunitscale 0.1 \linewd 0.5

\htext(-40 10){${(-q^2)^{l(b)}\times}$}

\move(-20 -20)\lvec(-20 0)\lvec(0 0)\lvec(0 10)\lvec(60
10)\lvec(60 30)\lvec(75 30)\lvec(75 -20)\lvec(-20 -20)

\move(10 10)\lvec(20 20)\lvec(10 20)\lvec(10 10)\lfill f:0.8

\move(20 10)\lvec(30 20)\lvec(20 20)\lvec(20 10)\lfill f:0.8

\move(50 10)\lvec(60 20)\lvec(50 20)\lvec(50 10)\lfill f:0.8

\move(50 10)\lvec(60 20)\lvec(60 10)\lvec(50 10)\lfill f:0.8


\htext(37 14){$_{\cdots}$}

\end{texdraw}}\hskip 10mm and
\hskip 2mm \raisebox{-0.4\height}{\begin{texdraw}\fontsize{9}{9}
\drawdim em \setunitscale 0.1 \linewd 0.5

\htext(-40 10){${(-q^2)^{l(b)}\times}$}

\move(-20 -20)\lvec(-20 0)\lvec(0 0)\lvec(0 10)\lvec(60
10)\lvec(60 30)\lvec(75 30)\lvec(75 -20)\lvec(-20 -20)

\move(10 10)\lvec(20 20)\lvec(20 10)\lvec(10 10)\lfill f:0.8

\move(20 10)\lvec(30 20)\lvec(30 10)\lvec(20 10)\lfill f:0.8

\move(50 10)\lvec(60 20)\lvec(60 10)\lvec(50 10)\lfill f:0.8

\move(50 10)\lvec(60 20)\lvec(50 20)\lvec(50 10)\lfill f:0.8


\htext(37 14){$_{\cdots}$}

\end{texdraw}}\hskip 5mm ,
\end{center}\vskip 3mm

\noindent respectively, where $l(b)\geq 1$ is given in the above
figure. In either case, let $Y_L(b)=(\cdots,y_{k+2},y_{k+1})$ and
set $L_i(b;Y)=\varphi_i(Y_L(b))-\varepsilon_i(Y_L(b))$. Then we
define
\begin{equation}
f_i\,Y=\sum_{b}q^{2L_i(b;Y)}(Y\swarrow b),
\end{equation}
\vskip -2mm \noindent where $b$ runs over all admissible and
virtually admissible $i$-slots in $Y$.

\begin{ex}{\rm \mbox{}\vskip 5mm

(1) \hskip 5mm $e_1$ \raisebox{-0.4\height}{\begin{texdraw}
\drawdim mm \setunitscale 0.5 \fontsize{7}{7}\selectfont \textref
h:C v:C

\move(0 0)\recs \move(10 0)\recs \move(-20 0)\recs \move(-10
0)\recs \move(-30 0)\recs

\move(-30 0)\bsegment

\move(0 0)\lvec(10 0)\lvec(10 10)\lvec(0 10)\lvec(0 0) \move(0
5)\lvec(10 5)\htext(5 2.5){$1$}\htext(5 7.5){$1$}

\esegment
\move(-20 0)\bsegment

\move(0 0)\lvec(10 0)\lvec(10 10)\lvec(0 10)\lvec(0 0) \move(0
5)\lvec(10 5)\htext(5 2.5){$1$}\htext(5 7.5){$1$}

\move(0 10)\lvec(10 10)\lvec(10 20)\lvec(0 20)\lvec(0 10)\move(0
10)\lvec(10 20)\htext(2.5 17.5){$2$}\htext(7.5 12.5){$0$}

\move(0 20)\lvec(10 20)\lvec(10 25)\lvec(0 25)\lvec(0 20)\move(0
25)\lvec(10 25)\htext(5 22.5){$1$}

\esegment
\move(-10 0)\bsegment

\move(0 0)\lvec(10 0)\lvec(10 10)\lvec(0 10)\lvec(0 0) \move(0
5)\lvec(10 5)\htext(5 2.5){$1$}\htext(5 7.5){$1$}

\move(0 10)\lvec(10 10)\lvec(10 20)\lvec(0 20)\lvec(0 10)\move(0
10)\lvec(10 20)\htext(2.5 17.5){$0$}\htext(7.5 12.5){$2$}

\move(0 20)\lvec(10 20)\lvec(10 25)\lvec(0 25)\lvec(0 20)\move(0
25)\lvec(10 25)\htext(5 22.5){$1$} \esegment

\move(0 0)\bsegment

\move(0 0)\lvec(10 0)\lvec(10 10)\lvec(0 10)\lvec(0 0) \move(0
5)\lvec(10 5)\htext(5 2.5){$1$}\htext(5 7.5){$1$}

\move(0 10)\lvec(10 10)\lvec(10 20)\lvec(0 20)\lvec(0 10)\move(0
10)\lvec(10 20)\htext(2.5 17.5){$2$}\htext(7.5 12.5){$0$}

\move(0 20)\lvec(10 20)\lvec(10 25)\lvec(0 25)\lvec(0 20)\move(0
25)\lvec(10 25)\htext(5 22.5){$1$} \esegment
\move(10 0)\bsegment

\move(0 0)\lvec(10 0)\lvec(10 10)\lvec(0 10)\lvec(0 0) \move(0
5)\lvec(10 5)\htext(5 2.5){$1$}\htext(5 7.5){$1$}

\move(0 10)\lvec(10 10)\lvec(10 20)\lvec(0 20)\lvec(0 10)\move(0
10)\lvec(10 20)\htext(2.5 17.5){$0$}\htext(7.5 12.5){$2$}

\move(0 20)\lvec(10 20)\lvec(10 30)\lvec(0 30)\lvec(0 20)\move(0
25)\lvec(10 25)\htext(5 22.5){$1$}\htext(5 27.5){$1$}

\esegment
\end{texdraw}}\hskip 3mm$=$\hskip 5mm
$q^2$~ \raisebox{-0.4\height}{\begin{texdraw} \drawdim mm
\setunitscale 0.5 \fontsize{7}{7}\selectfont \textref h:C v:C

\move(0 0)\recs \move(10 0)\recs \move(-20 0)\recs \move(-10
0)\recs

\move(-20 0)\bsegment

\move(0 0)\lvec(10 0)\lvec(10 10)\lvec(0 10)\lvec(0 0) \move(0
5)\lvec(10 5)\htext(5 2.5){$1$}\htext(5 7.5){$1$}

\move(0 10)\lvec(10 10)\lvec(10 20)\lvec(0 20)\lvec(0 10)\move(0
10)\lvec(10 20)\htext(2.5 17.5){$2$}\htext(7.5 12.5){$0$}

\move(0 20)\lvec(10 20)\lvec(10 25)\lvec(0 25)\lvec(0 20)\move(0
25)\lvec(10 25)\htext(5 22.5){$1$}

\esegment
\move(-10 0)\bsegment

\move(0 0)\lvec(10 0)\lvec(10 10)\lvec(0 10)\lvec(0 0) \move(0
5)\lvec(10 5)\htext(5 2.5){$1$}\htext(5 7.5){$1$}

\move(0 10)\lvec(10 10)\lvec(10 20)\lvec(0 20)\lvec(0 10)\move(0
10)\lvec(10 20)\htext(2.5 17.5){$0$}\htext(7.5 12.5){$2$}

\move(0 20)\lvec(10 20)\lvec(10 25)\lvec(0 25)\lvec(0 20)\move(0
25)\lvec(10 25)\htext(5 22.5){$1$} \esegment

\move(0 0)\bsegment

\move(0 0)\lvec(10 0)\lvec(10 10)\lvec(0 10)\lvec(0 0) \move(0
5)\lvec(10 5)\htext(5 2.5){$1$}\htext(5 7.5){$1$}

\move(0 10)\lvec(10 10)\lvec(10 20)\lvec(0 20)\lvec(0 10)\move(0
10)\lvec(10 20)\htext(2.5 17.5){$2$}\htext(7.5 12.5){$0$}

\move(0 20)\lvec(10 20)\lvec(10 25)\lvec(0 25)\lvec(0 20)\move(0
25)\lvec(10 25)\htext(5 22.5){$1$} \esegment
\move(10 0)\bsegment

\move(0 0)\lvec(10 0)\lvec(10 10)\lvec(0 10)\lvec(0 0) \move(0
5)\lvec(10 5)\htext(5 2.5){$1$}\htext(5 7.5){$1$}

\move(0 10)\lvec(10 10)\lvec(10 20)\lvec(0 20)\lvec(0 10)\move(0
10)\lvec(10 20)\htext(2.5 17.5){$0$}\htext(7.5 12.5){$2$}

\move(0 20)\lvec(10 20)\lvec(10 30)\lvec(0 30)\lvec(0 20)\move(0
25)\lvec(10 25)\htext(5 22.5){$1$}\htext(5 27.5){$1$}

\esegment

\end{texdraw}} \vskip 5mm

\hskip 3cm $+$\hskip 3mm $q(1+q^6)$
\raisebox{-0.4\height}{\begin{texdraw} \drawdim mm \setunitscale
0.5 \fontsize{7}{7}\selectfont \textref h:C v:C

\move(0 0)\recs \move(10 0)\recs \move(-20 0)\recs \move(-10
0)\recs \move(-30 0)\recs

\move(-30 0)\bsegment

\move(0 0)\lvec(10 0)\lvec(10 10)\lvec(0 10)\lvec(0 0) \move(0
5)\lvec(10 5)\htext(5 2.5){$1$}\htext(5 7.5){$1$}

\esegment
\move(-20 0)\bsegment

\move(0 0)\lvec(10 0)\lvec(10 10)\lvec(0 10)\lvec(0 0) \move(0
5)\lvec(10 5)\htext(5 2.5){$1$}\htext(5 7.5){$1$}

\move(0 10)\lvec(10 10)\lvec(10 20)\lvec(0 20)\lvec(0 10)\move(0
10)\lvec(10 20)\htext(2.5 17.5){$2$}\htext(7.5 12.5){$0$}
\esegment
\move(-10 0)\bsegment

\move(0 0)\lvec(10 0)\lvec(10 10)\lvec(0 10)\lvec(0 0) \move(0
5)\lvec(10 5)\htext(5 2.5){$1$}\htext(5 7.5){$1$}

\move(0 10)\lvec(10 10)\lvec(10 20)\lvec(0 20)\lvec(0 10)\move(0
10)\lvec(10 20)\htext(2.5 17.5){$0$}\htext(7.5 12.5){$2$}

\move(0 20)\lvec(10 20)\lvec(10 25)\lvec(0 25)\lvec(0 20)\move(0
25)\lvec(10 25)\htext(5 22.5){$1$} \esegment

\move(0 0)\bsegment

\move(0 0)\lvec(10 0)\lvec(10 10)\lvec(0 10)\lvec(0 0) \move(0
5)\lvec(10 5)\htext(5 2.5){$1$}\htext(5 7.5){$1$}

\move(0 10)\lvec(10 10)\lvec(10 20)\lvec(0 20)\lvec(0 10)\move(0
10)\lvec(10 20)\htext(2.5 17.5){$2$}\htext(7.5 12.5){$0$}

\move(0 20)\lvec(10 20)\lvec(10 25)\lvec(0 25)\lvec(0 20)\move(0
25)\lvec(10 25)\htext(5 22.5){$1$} \esegment
\move(10 0)\bsegment

\move(0 0)\lvec(10 0)\lvec(10 10)\lvec(0 10)\lvec(0 0) \move(0
5)\lvec(10 5)\htext(5 2.5){$1$}\htext(5 7.5){$1$}

\move(0 10)\lvec(10 10)\lvec(10 20)\lvec(0 20)\lvec(0 10)\move(0
10)\lvec(10 20)\htext(2.5 17.5){$0$}\htext(7.5 12.5){$2$}

\move(0 20)\lvec(10 20)\lvec(10 30)\lvec(0 30)\lvec(0 20)\move(0
25)\lvec(10 25)\htext(5 22.5){$1$}\htext(5 27.5){$1$}

\esegment
\end{texdraw}
}\hskip 3mm $+$\hskip 3mm \raisebox{-0.4\height}{\begin{texdraw}
\drawdim mm \setunitscale 0.5 \fontsize{7}{7}\selectfont \textref
h:C v:C

\move(0 0)\recs \move(10 0)\recs \move(-20 0)\recs \move(-10
0)\recs \move(-30 0)\recs

\move(-30 0)\bsegment

\move(0 0)\lvec(10 0)\lvec(10 10)\lvec(0 10)\lvec(0 0) \move(0
5)\lvec(10 5)\htext(5 2.5){$1$}\htext(5 7.5){$1$}

\esegment
\move(-20 0)\bsegment

\move(0 0)\lvec(10 0)\lvec(10 10)\lvec(0 10)\lvec(0 0) \move(0
5)\lvec(10 5)\htext(5 2.5){$1$}\htext(5 7.5){$1$}

\move(0 10)\lvec(10 10)\lvec(10 20)\lvec(0 20)\lvec(0 10)\move(0
10)\lvec(10 20)\htext(2.5 17.5){$2$}\htext(7.5 12.5){$0$}

\move(0 20)\lvec(10 20)\lvec(10 25)\lvec(0 25)\lvec(0 20)\move(0
25)\lvec(10 25)\htext(5 22.5){$1$}

\esegment
\move(-10 0)\bsegment

\move(0 0)\lvec(10 0)\lvec(10 10)\lvec(0 10)\lvec(0 0) \move(0
5)\lvec(10 5)\htext(5 2.5){$1$}\htext(5 7.5){$1$}

\move(0 10)\lvec(10 10)\lvec(10 20)\lvec(0 20)\lvec(0 10)\move(0
10)\lvec(10 20)\htext(2.5 17.5){$0$}\htext(7.5 12.5){$2$}

\move(0 20)\lvec(10 20)\lvec(10 25)\lvec(0 25)\lvec(0 20)\move(0
25)\lvec(10 25)\htext(5 22.5){$1$} \esegment

\move(0 0)\bsegment

\move(0 0)\lvec(10 0)\lvec(10 10)\lvec(0 10)\lvec(0 0) \move(0
5)\lvec(10 5)\htext(5 2.5){$1$}\htext(5 7.5){$1$}

\move(0 10)\lvec(10 10)\lvec(10 20)\lvec(0 20)\lvec(0 10)\move(0
10)\lvec(10 20)\htext(2.5 17.5){$2$}\htext(7.5 12.5){$0$}

\move(0 20)\lvec(10 20)\lvec(10 25)\lvec(0 25)\lvec(0 20)\move(0
25)\lvec(10 25)\htext(5 22.5){$1$} \esegment
\move(10 0)\bsegment

\move(0 0)\lvec(10 0)\lvec(10 10)\lvec(0 10)\lvec(0 0) \move(0
5)\lvec(10 5)\htext(5 2.5){$1$}\htext(5 7.5){$1$}

\move(0 10)\lvec(10 10)\lvec(10 20)\lvec(0 20)\lvec(0 10)\move(0
10)\lvec(10 20)\htext(2.5 17.5){$0$}\htext(7.5 12.5){$2$}

\move(0 20)\lvec(10 20)\lvec(10 25)\lvec(0 25)\lvec(0 20)\htext(5
22.5){$1$}

\esegment
\end{texdraw}
}\vskip 5mm

(2) \hskip 5mm $f_2$ \raisebox{-0.6\height}{\begin{texdraw}
\drawdim mm \setunitscale 0.5 \fontsize{7}{7}\selectfont \textref
h:C v:C

\move(0 10)\tris  \move(10 10)\tris \move(-20 10)\tris \move(-10
10)\tris

\move(-20 0)\bsegment

\move(0 10)\lvec(10 10)\lvec(10 20)\lvec(0 20)\lvec(0 10)\move(0
10)\lvec(10 20)\htext(2.5 17.5){$0$}\htext(7.5 12.5){$2$}

\move(0 20)\lvec(10 20)\lvec(10 30)\lvec(0 30)\lvec(0 20)\move(0
25)\lvec(10 25)\htext(5 22.5){$1$}\htext(5 27.5){$1$}

\esegment

\move(-10 0)\bsegment

\move(0 10)\lvec(10 10)\lvec(10 20)\lvec(0 20)\lvec(0 10)\move(0
10)\lvec(10 20)\htext(2.5 17.5){$2$}\htext(7.5 12.5){$0$}

\move(0 20)\lvec(10 20)\lvec(10 30)\lvec(0 30)\lvec(0 20)\move(0
25)\lvec(10 25)\htext(5 22.5){$1$}\htext(5 27.5){$1$}

\move(0 30)\lvec(10 30)\lvec(10 40)\lvec(0 30)\htext(7.5
32.5){$0$}

\esegment

\move(0 0)\bsegment

\move(0 10)\lvec(10 10)\lvec(10 20)\lvec(0 20)\lvec(0 10)\move(0
10)\lvec(10 20)\htext(2.5 17.5){$0$}\htext(7.5 12.5){$2$}

\move(0 20)\lvec(10 20)\lvec(10 30)\lvec(0 30)\lvec(0 20)\move(0
25)\lvec(10 25)\htext(5 22.5){$1$}\htext(5 27.5){$1$}

\move(0 30)\lvec(10 30)\lvec(10 40)\lvec(0 30)\htext(7.5
32.5){$2$}

\esegment

\move(10 0)\bsegment

\move(0 10)\lvec(10 10)\lvec(10 20)\lvec(0 20)\lvec(0 10)\move(0
10)\lvec(10 20)\htext(2.5 17.5){$2$}\htext(7.5 12.5){$0$}

\move(0 20)\lvec(10 20)\lvec(10 30)\lvec(0 30)\lvec(0 20)\move(0
25)\lvec(10 25)\htext(5 22.5){$1$}\htext(5 27.5){$1$}

\move(0 30)\lvec(10 30)\lvec(10 40)\lvec(0 30)\htext(7.5
32.5){$0$}

\esegment
\end{texdraw}}\hskip 3mm$=$\hskip 5mm
$q^2$ \raisebox{-0.6\height}{\begin{texdraw} \drawdim mm
\setunitscale 0.5 \fontsize{7}{7}\selectfont \textref h:C v:C

\move(0 10)\tris  \move(10 10)\tris \move(-20 10)\tris \move(-10
10)\tris

\move(-20 0)\bsegment

\move(0 10)\lvec(10 10)\lvec(10 20)\lvec(0 20)\lvec(0 10)\move(0
10)\lvec(10 20)\htext(2.5 17.5){$0$}\htext(7.5 12.5){$2$}

\move(0 20)\lvec(10 20)\lvec(10 30)\lvec(0 30)\lvec(0 20)\move(0
25)\lvec(10 25)\htext(5 22.5){$1$}\htext(5 27.5){$1$}

\move(0 30)\lvec(10 30)\lvec(10 40)\lvec(0 30)\htext(7.5
32.5){$2$}

\esegment

\move(-10 0)\bsegment

\move(0 10)\lvec(10 10)\lvec(10 20)\lvec(0 20)\lvec(0 10)\move(0
10)\lvec(10 20)\htext(2.5 17.5){$2$}\htext(7.5 12.5){$0$}

\move(0 20)\lvec(10 20)\lvec(10 30)\lvec(0 30)\lvec(0 20)\move(0
25)\lvec(10 25)\htext(5 22.5){$1$}\htext(5 27.5){$1$}

\move(0 30)\lvec(10 30)\lvec(10 40)\lvec(0 30)\htext(7.5
32.5){$0$}

\esegment

\move(0 0)\bsegment

\move(0 10)\lvec(10 10)\lvec(10 20)\lvec(0 20)\lvec(0 10)\move(0
10)\lvec(10 20)\htext(2.5 17.5){$0$}\htext(7.5 12.5){$2$}

\move(0 20)\lvec(10 20)\lvec(10 30)\lvec(0 30)\lvec(0 20)\move(0
25)\lvec(10 25)\htext(5 22.5){$1$}\htext(5 27.5){$1$}

\move(0 30)\lvec(10 30)\lvec(10 40)\lvec(0 30)\htext(7.5
32.5){$2$}

\esegment

\move(10 0)\bsegment

\move(0 10)\lvec(10 10)\lvec(10 20)\lvec(0 20)\lvec(0 10)\move(0
10)\lvec(10 20)\htext(2.5 17.5){$2$}\htext(7.5 12.5){$0$}

\move(0 20)\lvec(10 20)\lvec(10 30)\lvec(0 30)\lvec(0 20)\move(0
25)\lvec(10 25)\htext(5 22.5){$1$}\htext(5 27.5){$1$}

\move(0 30)\lvec(10 30)\lvec(10 40)\lvec(0 30)\htext(7.5
32.5){$0$}

\esegment
\end{texdraw}} \vskip 5mm

\hskip 3cm$+$\hskip 3mm $q^8$
 \raisebox{-0.6\height}{\begin{texdraw}
\drawdim mm \setunitscale 0.5 \fontsize{7}{7}\selectfont \textref
h:C v:C

\move(0 10)\tris  \move(10 10)\tris \move(-20 10)\tris \move(-10
10)\tris

\move(-20 0)\bsegment

\move(0 10)\lvec(10 10)\lvec(10 20)\lvec(0 20)\lvec(0 10)\move(0
10)\lvec(10 20)\htext(2.5 17.5){$0$}\htext(7.5 12.5){$2$}

\move(0 20)\lvec(10 20)\lvec(10 30)\lvec(0 30)\lvec(0 20)\move(0
25)\lvec(10 25)\htext(5 22.5){$1$}\htext(5 27.5){$1$}

\esegment

\move(-10 0)\bsegment

\move(0 10)\lvec(10 10)\lvec(10 20)\lvec(0 20)\lvec(0 10)\move(0
10)\lvec(10 20)\htext(2.5 17.5){$2$}\htext(7.5 12.5){$0$}

\move(0 20)\lvec(10 20)\lvec(10 30)\lvec(0 30)\lvec(0 20)\move(0
25)\lvec(10 25)\htext(5 22.5){$1$}\htext(5 27.5){$1$}

\move(0 30)\lvec(0 40)\lvec(10 40)\lvec(0 30)\htext(2.5 37.5){$2$}

\esegment

\move(0 0)\bsegment

\move(0 10)\lvec(10 10)\lvec(10 20)\lvec(0 20)\lvec(0 10)\move(0
10)\lvec(10 20)\htext(2.5 17.5){$0$}\htext(7.5 12.5){$2$}

\move(0 20)\lvec(10 20)\lvec(10 30)\lvec(0 30)\lvec(0 20)\move(0
25)\lvec(10 25)\htext(5 22.5){$1$}\htext(5 27.5){$1$}

\move(0 30)\lvec(0 40)\lvec(10 40)\lvec(0 30)\htext(2.5 37.5){$0$}

\esegment

\move(10 0)\bsegment

\move(0 10)\lvec(10 10)\lvec(10 20)\lvec(0 20)\lvec(0 10)\move(0
10)\lvec(10 20)\htext(2.5 17.5){$2$}\htext(7.5 12.5){$0$}

\move(0 20)\lvec(10 20)\lvec(10 30)\lvec(0 30)\lvec(0 20)\move(0
25)\lvec(10 25)\htext(5 22.5){$1$}\htext(5 27.5){$1$}

\move(0 30)\lvec(10 30)\lvec(10 40)\lvec(0 30)\lvec(0 40)\lvec(10
40)\htext(2.5 37.5){$2$}\htext(7.5 32.5){$0$}

\esegment
\end{texdraw}}
\hskip 3mm $+$ \hskip 3mm $q^4$
\raisebox{-0.6\height}{\begin{texdraw} \drawdim mm \setunitscale
0.5 \fontsize{7}{7}\selectfont \textref h:C v:C

\move(0 10)\tris  \move(10 10)\tris \move(-20 10)\tris \move(-10
10)\tris

\move(-20 0)\bsegment

\move(0 10)\lvec(10 10)\lvec(10 20)\lvec(0 20)\lvec(0 10)\move(0
10)\lvec(10 20)\htext(2.5 17.5){$0$}\htext(7.5 12.5){$2$}

\move(0 20)\lvec(10 20)\lvec(10 30)\lvec(0 30)\lvec(0 20)\move(0
25)\lvec(10 25)\htext(5 22.5){$1$}\htext(5 27.5){$1$}

\esegment

\move(-10 0)\bsegment

\move(0 10)\lvec(10 10)\lvec(10 20)\lvec(0 20)\lvec(0 10)\move(0
10)\lvec(10 20)\htext(2.5 17.5){$2$}\htext(7.5 12.5){$0$}

\move(0 20)\lvec(10 20)\lvec(10 30)\lvec(0 30)\lvec(0 20)\move(0
25)\lvec(10 25)\htext(5 22.5){$1$}\htext(5 27.5){$1$}

\move(0 30)\lvec(10 30)\lvec(10 40)\lvec(0 30)\htext(7.5
32.5){$0$}

\esegment

\move(0 0)\bsegment

\move(0 10)\lvec(10 10)\lvec(10 20)\lvec(0 20)\lvec(0 10)\move(0
10)\lvec(10 20)\htext(2.5 17.5){$0$}\htext(7.5 12.5){$2$}

\move(0 20)\lvec(10 20)\lvec(10 30)\lvec(0 30)\lvec(0 20)\move(0
25)\lvec(10 25)\htext(5 22.5){$1$}\htext(5 27.5){$1$}

\move(0 30)\lvec(10 30)\lvec(10 40)\lvec(0 30)\htext(7.5
32.5){$2$}

\esegment

\move(10 0)\bsegment

\move(0 10)\lvec(10 10)\lvec(10 20)\lvec(0 20)\lvec(0 10)\move(0
10)\lvec(10 20)\htext(2.5 17.5){$2$}\htext(7.5 12.5){$0$}

\move(0 20)\lvec(10 20)\lvec(10 30)\lvec(0 30)\lvec(0 20)\move(0
25)\lvec(10 25)\htext(5 22.5){$1$}\htext(5 27.5){$1$}

\move(0 30)\lvec(10 30)\lvec(10 40)\lvec(0 30)\lvec(0 40)\lvec(10
40)\htext(2.5 37.5){$2$}\htext(7.5 32.5){$0$}

\esegment
\end{texdraw}}

}
\end{ex}

\begin{thm}
${\mathcal F}(\Lambda)$ is an integrable $U_q(C_2^{(1)})$-module.
\end{thm}
\pf First, it follows directly from the definition of the actions
of $U_q(C_2^{(1)})$ that
\begin{equation}\label {rel1}
\begin{split}
q^hq^{h'}Y&=q^{h+h'}Y, \\
q^h e_i q^{-h}Y&=q^{\langle h,\alpha_i \rangle}e_i Y, \\
q^h f_i q^{-h}Y&=q^{-\langle h,\alpha_i \rangle}f_i Y
\end{split}
\end{equation}
for $Y\in {\mathcal Z}(\Lambda)$, $i\in I$ and $h,h'\in P^{\vee}$.
Since $e_i$ and $f_i$ ($i\in I$) act locally nilpotently on
${\mathcal F}(\Lambda)$, if we show that
\begin{equation}\label {rel2}
\begin{split}
[e_i,f_j]Y=\delta_{ij}\frac{K_i-K_i^{-1}}{q_i-q_i^{-1}}Y
\end{split}
\end{equation}
for $Y\in {\mathcal Z}(\Lambda)$ and $i,j\in I$, then the Serre
relations will follow from Proposition B.1 in \cite{KMPY}. The
verification of \eqref{rel2} can be obtained by modifying the
arguments in \cite{KK}.\qed\vskip 3mm

Let ${\mathcal L}(\Lambda)=\bigoplus_{Y\in{\mathcal
Z}(\Lambda)}\mathbb{A}_0 Y$. Then we also obtain

\begin{thm}\label {crystal basis for Fock space}
The pair $({\mathcal L}(\Lambda),{\mathcal Z}(\Lambda))$ is a
crystal basis of ${\mathcal F}(\Lambda)$ and the crystal of
${\mathcal F}(\Lambda)$ is isomorphic to the affine crystal
${\mathcal Z}(\Lambda)$ defined in Section 4.
\end{thm}
\pf Let us give a sketch of the proof (see \cite{KK} for a
detailed argument). It is clear that $({\mathcal
L}(\Lambda),{\mathcal Z}(\Lambda))$ satisfies the first four
conditions in Definition \ref{crystal basis}. Fix $i\in I$. For
each $Y\in{\mathcal Z}(\Lambda)$, we can find a subset $B_{Y,i}$
of ${\mathcal Z}(\Lambda)$ containing $Y$ such that
\begin{itemize}
\item[(i)] $V_{Y,i}=\bigoplus_{Z\in B_{Y,i}}\mathbb{Q}(q)Z$ is a $U_{(i)}$-submodule
of ${\mathcal F}(\Lambda)$,

\item[(ii)] $(L_{Y,i},B_{Y,i})$ is a crystal basis of $V_{Y,i}$ where
$L_{Y,i}=\bigoplus_{Z\in B_{Y,i}}\mathbb{A}_0 Z$,

\item[(iii)] $\tilde{e}_iZ=\tilde{E}_iZ$ and $\tilde{f}_iZ=\tilde{F}_iZ$
for all $Z\in B_{Y,i}$.
\end{itemize}
From (i) and (ii), it follows that $({\mathcal
L}(\Lambda),{\mathcal Z}(\Lambda))$ satisfies the rest three
conditions in Definition \ref{crystal basis}. The condition (iii)
implies that the Kashiwara operators on ${\mathcal Z}(\Lambda)$
induced by the $U_q(C_2^{(1)})$ action on ${\mathcal F}(\Lambda)$
coincide with the abstract Kashiwara operators on ${\mathcal
Z}(\Lambda)$. Therefore, the crystal of ${\mathcal F}(\Lambda)$ is
isomorphic to the abstract affine crystal ${\mathcal Z}(\Lambda)$
defined in Section 4. \qed

\begin{cor}
\begin{equation}
{\mathcal
F}(\Lambda)=\bigoplus_{m=0}^{\infty}V(\Lambda-m\delta)^{\oplus
p(m)},
\end{equation}
where $p(m)$ is the number of partitions of $m$.
\end{cor}
\pf We will show that
\begin{itemize}
\item[(i)] the weight of each maximal vector in ${\mathcal Z}(\Lambda)$
is of the form $\Lambda-m\delta$ for some $m\geq 0$,

\item[(ii)] there exists a bijection between the set of partitions of $m$ ($m\geq
0$) and the set of maximal vectors in ${\mathcal Z}(\Lambda)$ with
weight $\Lambda-m\delta$.
\end{itemize}

Let $Y=(y_k)_{k=0}^{\infty}\in{\mathcal Z}(\Lambda)$ be a maximal
vector, that is, $\tilde{E}_iY=0$ for all $i\in I$. Suppose that
$Y$ is the ground state wall $Y_{\Lambda}$. Since ${\rm
wt}(Y_{\Lambda})=\Lambda$ and ${\mathcal
Z}(\Lambda)_{\Lambda}=\{\,Y_{\Lambda}\,\}$, the multiplicity of
$V(\Lambda)$ in ${\mathcal F}(\Lambda)$ is $1$. Now, we assume
that $Y\neq Y_{\Lambda}$. Let $l$ be the largest integer such that
$|y_l|\neq 0$. Suppose that $Y'=(y_k)_{k=l+1}^{\infty}\in{\mathcal
Y}(\Lambda_j)$ for some $j\in I$.\vskip 3mm

{\bf Case 1.} $j=1$

We see from the pattern for ${\mathcal Z}(\Lambda)$ that
$\Lambda=\Lambda_1$. By the maximality of $Y$, the $1$-signature
of $y_{l+1}$ is $+$ and the $1$-signature of $y_l$ is $-$ or $-+$.
Hence $y_l$ is obtained by adding some $\delta$-columns to the
ground state wall. If $l\neq 0$, let $l'$ be the largest integer
such that $l'<l$ and $|y_{l'}|>|y_l|$. Again by the maximality of
$Y$, the $1$-signature of $y_{l'}$ is $-$ or $-+$, which means
that $y_{l'}$ is also obtained by adding some $\delta$-columns to
the ground state wall. Repeating the above argument from left to
right, we conclude that for $0\leq k\leq l$, the total volume of
the blocks added on the $k$th column is $2m_k$ for some $m_k\geq
1$. Hence, $(m_0,\cdots,m_l)$ forms a partition, and ${\rm
wt}(Y)=\Lambda-(\sum_{k=0}^lm_k)\delta$.\vskip 3mm

{\bf Case 2.} $j=0,2$

We may assume that $j=0$. By the maximality of $Y$, we observe
that the $y_{l+1}$ is $0$-admissible and $y_{l}$ is $0$-removable
but not $2$-removable. Hence $y_l$ is obtained by adding some
$\delta$-columns to the ground-state wall. If $l\neq 0$, let $l'$
be the largest integer such that $l'<l$ and $|y_{l'}|>|y_l|$. If
$y_{l'+1}$ is $0$-admissible, then by the maximality of $Y$,
$y_{l'}$ is $0$-removable but not $2$-removable. On the other
hand, if $y_{l'+1}$ is $2$-admissible, then by the maximality of
$Y$, $y_{l'}$ is $2$-removable but not $0$-removable. As in Case
1, by repeating the above argument, we conclude that for $0\leq
k\leq l$, the total volume of the blocks added on the $k$th column
is $2m_k$ for some $m_k\geq 1$. Hence, $(m_0,m_1,\cdots,m_l)$
forms a partition and ${\rm
wt}(Y)=\Lambda-(\sum_{k=0}^lm_k)\delta$.\vskip 3mm

Conversely, for a given partition $(m_k)_{k=0}^{\infty}$ of $m\geq
0 $, we can find a unique $Y=(y_k)_{k=0}^{\infty}\in {\mathcal
Z}(\Lambda)$ such that $y_k$ is obtained by adding $m_k$ many
$\delta$-columns to the $k$th column of the ground state wall
$Y_{\Lambda}$ (hence the total volume of the blocks added to the
$k$th column is $2m_k$). It is easy to check that $Y$ is a maximal
vector with ${\rm wt}(Y)=\Lambda-m\delta$.  \qed

\begin{ex}{\rm In ${\mathcal Y}(\Lambda_1)$, the maximal vector
of weight $\Lambda-4\delta$ corresponding to the partition
$(2,1,1)$ is \vskip 5mm

\begin{center}
{\begin{texdraw} \drawdim mm \setunitscale 0.5
\fontsize{7}{7}\selectfont \textref h:C v:C

\move(0 0)\recs  \move(10 0)\recs \move(-10 0)\recs

\move(-10 0)\bsegment

\move(0 0)\lvec(10 0)\lvec(10 10)\lvec(0 10)\lvec(0 0) \move(0
5)\lvec(10 5)\htext(5 2.5){$1$}\htext(5 7.5){$1$}

\move(0 10)\lvec(10 10)\lvec(10 20)\lvec(0 20)\lvec(0 10)\move(0
10)\lvec(10 20)\htext(2.5 17.5){$0$}\htext(7.5 12.5){$2$}

\move(0 20)\lvec(10 20)\lvec(10 25)\lvec(0 25)\lvec(0 20)\htext(5
22.5){$1$} \esegment

\move(0 0)\bsegment

\move(0 0)\lvec(10 0)\lvec(10 10)\lvec(0 10)\lvec(0 0) \move(0
5)\lvec(10 5)\htext(5 2.5){$1$}\htext(5 7.5){$1$}

\move(0 10)\lvec(10 10)\lvec(10 20)\lvec(0 20)\lvec(0 10)\move(0
10)\lvec(10 20)\htext(2.5 17.5){$2$}\htext(7.5 12.5){$0$}

\move(0 20)\lvec(10 20)\lvec(10 25)\lvec(0 25)\lvec(0 20)\htext(5
22.5){$1$}

\esegment

\move(10 0)\bsegment

\move(0 0)\lvec(10 0)\lvec(10 10)\lvec(0 10)\lvec(0 0) \move(0
5)\lvec(10 5)\htext(5 2.5){$1$}\htext(5 7.5){$1$}

\move(0 10)\lvec(10 10)\lvec(10 20)\lvec(0 20)\lvec(0 10)\move(0
10)\lvec(10 20)\htext(2.5 17.5){$0$}\htext(7.5 12.5){$2$}

\move(0 20)\lvec(10 20)\lvec(10 30)\lvec(0 30)\lvec(0 20)\move(0
25)\lvec(10 25)\htext(5 22.5){$1$}\htext(5 27.5){$1$}

\move(0 30)\lvec(10 30)\lvec(10 40)\lvec(0 40)\lvec(0 30)\move(0
30)\lvec(10 40)\htext(2.5 37.5){$0$}\htext(7.5 32.5){$2$}

\move(0 40)\lvec(10 40)\lvec(10 45)\lvec(0 45)\lvec(0 40)\htext(5
42.5){$1$}

\esegment

\end{texdraw}}\hskip 5mm .
\end{center}}
\end{ex}

\section{Global bases}
In this section, we will describe an algorithm for computing the
global basis for the basic representation $V(\Lambda)$ where
$\Lambda=\Lambda_i$ ($i\in I$). We have seen that there exists an
embedding $V(\Lambda)\hookrightarrow {\mathcal F}(\Lambda)$, and
that ${\mathcal Y}(\Lambda)$ is a crystal of $V(\Lambda)$. By
Theorem \ref{global basis}, there exists an $\mathbb{A}$-basis
$G(\Lambda)=\{\,G(Y)\,|\,Y\in {\mathcal Y}(\Lambda)\,\}$ of
$V(\Lambda)^{\mathbb{A}}$. For each $Y\in {\mathcal Y}(\Lambda)$,
the global basis element $G(Y)$ can be written as an
$\mathbb{A}$-linear combination of proper Young walls in
${\mathcal Z}(\Lambda)$. Hence, our algorithm is to compute the
coefficients of proper Young walls in each global basis element
$G(Y)$. We will follow the arguments in \cite{KK}.\vskip 3mm

We start with certain orderings. For $Y=(y_k)_{k=0}^{\infty}$ and
$Z=(z_k)_{k=0}^{\infty}$ in ${\mathcal Z}(\Lambda)$, consider
their associated partitions $|Y|$ and $|Z|$. We define
$|Y|\unrhd|Z|$ if and only if
$\sum_{k=l}^{\infty}|y_k|\geq\sum_{k=l}^{\infty}|z_k|$ for all
$l\geq 0$. Note that it is not a partial ordering on ${\mathcal
Z}(\Lambda)$ since there exist $Y\neq Z$ in ${\mathcal
Z}(\Lambda)$ such that $|Y|=|Z|$. We also define $|Y|> |Z|$ if
$|y_k|> |z_k|$ where $k$ is the largest integer such that
$|y_k|\neq|z_k|$. Note that $|Y|\unrhd|Z|$ implies $|Y|\geq |Z|$.
Next, on the set of the proper Young walls with the same
associated partitions, we fix an arbitrary total ordering $\succ$.
Then we define a total ordering $>$ on ${\mathcal Z}(\Lambda)$ as
follows:
\begin{equation}
Y>Z\quad \Leftrightarrow  \quad (\,\,|Y|>|Z|\,\,) \,\,\text{ or
}\,\, (\,\,|Y|=|Z| \text{ and } Y \succ Z\,\,).
\end{equation}

For $Y\in {\mathcal Z}(\Lambda)$, we write
\begin{equation}
f_i^{(r)}Y=\sum_{\substack{Z\in {\mathcal Z}(\Lambda) \\ {\rm
wt}(Z)={\rm wt}(Y)-r\alpha_i }}Q_{Y,Z}(q)Z,
\end{equation}
where $i\in I$, $r\geq 1$, and $Q_{Y,Z}(q)\in \mathbb{Q}(q)$. If
$Z=(z_k)_{k=0}^{\infty}\in {\mathcal Z}(\Lambda)$ satisfies
$Q_{Y,Z}(q)\neq 0$, we can find a unique sequence of proper Young
walls $Y=Y_{0},Y_{1},\cdots,Y_{r}=Z$ such that

\begin{itemize}
\item[(i)] $\lambda_{k+1}Y_{k+1}=Y_{k}\swarrow b_{k+1}$ for a (virtually)
admissible $i$-slot $b_{k+1}$ of $Y_{k}$ and
$\lambda_{k+1}\in\mathbb{Z}[q,q^{-1}]$,
\item[(ii)] $b_{k+1}$ is placed on $b_k$ or to the right of $b_k$.
\end{itemize}
For each $k$, let $Q_{Y_k,Y_{k+1}}(q)$ be the coefficient of
$Y_{k+1}$ in the expression of $f_iY_{k}$. We define
\begin{equation}
Q^{\circ}_{Y,Z}(q)=\prod_{k=0}^{r-1}Q_{Y_{k},Y_{k+1}}(q)\in\mathbb{Z}[q,q^{-1}].
\end{equation}
Suppose that $i=0,2$ (that is, $i$-blocks are of type II). Then by
induction on $r$, we have
\begin{equation}\label{QYZ1}
Q_{Y,Z}(q)=
Q^{\circ}_{Y,Z}(q)q^{2\binom{r}{2}}\in\mathbb{Z}[q,q^{-1}].
\end{equation}

Suppose that $i=1$ (that is, $i$-blocks are of type I). Let us
assume that $b_k$ is located in the $i_k$th column of $Y_{k-1}$
($1\leq k\leq r$). Note that each $b_k$ ($1\leq k\leq r$) can be
viewed as a block (not necessarily removable) in $Z$. Set
\begin{equation}\label {JST}
\begin{split}
J_1&=\{\,k\,|\,b_{k-1} \text{ is beneath } b_k\,\}, \\
J_2&=\{\,k\,|\,\text{there exists a $1$-block {\rm(}$\neq
b_{k-1}${\rm )} beneath $b_k$}\,\}, \\
J_3&=\{\,k\,|\,\text{there exists no $1$-block on and beneath
$b_k$ in $Z$}\,\}, \\
S&=\{\,k\in J_2\,|\,k-1\in J_3 \text{ and }
|z_{i_{k-1}}|=|z_{i_k}|-1\,\}.
\end{split}
\end{equation}
Put $n_j=|J_j|$ ($j=1,2,3$), and $\mu_k=q\lambda_k$ ($k\in S$).
Note that $2n_1+n_2+n_3=r$. By induction on $r$, we have
\begin{equation}\label{QYZ2}
Q_{Y,Z}(q)=Q^{\circ}_{Y,Z}(q)
\dfrac{q^{\sigma(n_1,n_2,n_3)}}{[2]^{n_1}\prod_{k\in S}\mu_k},
\end{equation}
where
$\sigma(n_1,n_2,n_3)=4\binom{n_1}{2}+\binom{n_2}{2}+\binom{n_3}{2}+2n_1(n_2+n_3)+n_2n_3$
(see \cite{KK}). Note that $[2]^{n_1}\prod_{k\in S}\mu_k$ divides
$Q^{\circ}_{Y,Z}(q)$, which implies that
$Q_{Y,Z}(q)\in\mathbb{Z}[q,q^{-1}]$.

Therefore, in both cases, $Q_{Y,Z}(q)$ is a Laurent polynomial
with integral coefficients.

For $Y\in {\mathcal Z}(\Lambda)$, let $b$ be a block in $Y$. We
define the {\it coordinate of $b$} to be the pair $(k,l)$ if $b$
is located in the $k$th column of $Y$ and the maximal number of
unit cubes lying below $b$ is $l$. Note that a block in $Y$ is not
uniquely determined by its coordinates since two different blocks,
which form a unit cube, have the same coordinate.

For a given coordinate $c=(k,l)$ ($k,l\geq 0$), we define the {\it
ladder at $c$} to be the finite sequence of coordinates as
follows;
\begin{equation*}
c=(k,l), (k-1, l+2), (k-2, l+4), \cdots , (0,l+2k).
\end{equation*}

For $Y \in {\mathcal Y}(\Lambda)$, let $y_k$ be the left-most
column in $Y$ such that $|y_k|\neq 0$. Choose an $i$-block $b$
placed on top of $y_k$ with a coordinate $c=(k,l)$. If the
$i$-block is of type II and there is another block of type II on
top of $y_l$, we choose the block at the front. Let $L_c$ be the
ladder at $c$. Then it is the left-most ladder having nontrivial
intersection with $Y$. We define $\overline{Y}$ to be the proper
Young wall which is obtained by removing all the blocks in $Y$,
which are contained $L_c$. It is easy to see that $\overline{Y}$
is also reduced.
\begin{ex}{\rm Let $Y$ be a reduced proper Young wall given in the
following figure. Then $\overline{Y}$ is obtained by removing all
the blocks in the ladder $L_{(2,0)} : (2,0),(1,2),(0,4)$.\vskip
5mm

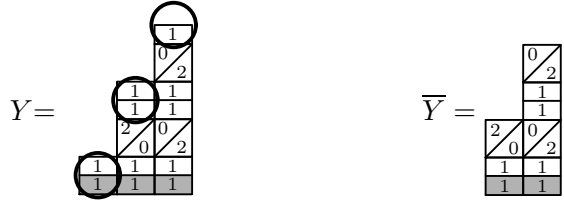
\begin{center}
$Y$=\hskip 2.5mm\raisebox{-0.4\height}{\begin{texdraw} \drawdim mm
\setunitscale 0.5 \fontsize{7}{7}\selectfont \textref h:C v:C

\move(0 0)\recs  \move(10 0)\recs \move(-10 0)\recs

\move(-10 0)\bsegment

\move(0 0)\lvec(10 0)\lvec(10 10)\lvec(0 10)\lvec(0 0) \move(0
5)\lvec(10 5)\htext(5 2.5){$1$}\htext(5 7.5){$1$}

\esegment

\move(0 0)\bsegment

\move(0 0)\lvec(10 0)\lvec(10 10)\lvec(0 10)\lvec(0 0) \move(0
5)\lvec(10 5)\htext(5 2.5){$1$}\htext(5 7.5){$1$}

\move(0 10)\lvec(10 10)\lvec(10 20)\lvec(0 20)\lvec(0 10)\move(0
10)\lvec(10 20)\htext(2.5 17.5){$2$}\htext(7.5 12.5){$0$}

\move(0 20)\lvec(10 20)\lvec(10 30)\lvec(0 30)\lvec(0 20)\move(0
25)\lvec(10 25)\htext(5 22.5){$1$}\htext(5 27.5){$1$}

\esegment

\move(10 0)\bsegment

\move(0 0)\lvec(10 0)\lvec(10 10)\lvec(0 10)\lvec(0 0) \move(0
5)\lvec(10 5)\htext(5 2.5){$1$}\htext(5 7.5){$1$}

\move(0 10)\lvec(10 10)\lvec(10 20)\lvec(0 20)\lvec(0 10)\move(0
10)\lvec(10 20)\htext(2.5 17.5){$0$}\htext(7.5 12.5){$2$}

\move(0 20)\lvec(10 20)\lvec(10 30)\lvec(0 30)\lvec(0 20)\move(0
25)\lvec(10 25)\htext(5 22.5){$1$}\htext(5 27.5){$1$}

\move(0 30)\lvec(10 30)\lvec(10 40)\lvec(0 40)\lvec(0 30)\move(0
30)\lvec(10 40)\htext(2.5 37.5){$0$}\htext(7.5 32.5){$2$}

\move(0 40)\lvec(10 40)\lvec(10 45)\lvec(0 45)\lvec(0 40)\htext(5
42.5){$1$}

\move(-15 5)\linewd 1 \lcir r:6 \move(-5 25)\linewd 1 \lcir r:6
\move(5 45)\linewd 1 \lcir r:6

\esegment
\end{texdraw}}\hskip 3cm $\overline{Y}=$
\raisebox{-0.5\height}{\begin{texdraw} \drawdim mm \setunitscale
0.5 \fontsize{7}{7}\selectfont \textref h:C v:C

\move(0 0)\recs  \move(10 0)\recs

\move(0 0)\bsegment

\move(0 0)\lvec(10 0)\lvec(10 10)\lvec(0 10)\lvec(0 0) \move(0
5)\lvec(10 5)\htext(5 2.5){$1$}\htext(5 7.5){$1$}

\move(0 10)\lvec(10 10)\lvec(10 20)\lvec(0 20)\lvec(0 10)\move(0
10)\lvec(10 20)\htext(2.5 17.5){$2$}\htext(7.5 12.5){$0$}

\esegment

\move(10 0)\bsegment

\move(0 0)\lvec(10 0)\lvec(10 10)\lvec(0 10)\lvec(0 0) \move(0
5)\lvec(10 5)\htext(5 2.5){$1$}\htext(5 7.5){$1$}

\move(0 10)\lvec(10 10)\lvec(10 20)\lvec(0 20)\lvec(0 10)\move(0
10)\lvec(10 20)\htext(2.5 17.5){$0$}\htext(7.5 12.5){$2$}

\move(0 20)\lvec(10 20)\lvec(10 30)\lvec(0 30)\lvec(0 20)\move(0
25)\lvec(10 25)\htext(5 22.5){$1$}\htext(5 27.5){$1$}

\move(0 30)\lvec(10 30)\lvec(10 40)\lvec(0 40)\lvec(0 30)\move(0
30)\lvec(10 40)\htext(2.5 37.5){$0$}\htext(7.5 32.5){$2$}

\esegment

\end{texdraw}}

\end{center}\vskip 5mm}
\end{ex}

Furthermore,  if $Y$ is a reduced proper Young wall and
$\overline{Y}$ is obtained by removing $r$ many $i$-blocks from
$Y$, then the coefficient of $Y$ in $f_i^{(r)}\overline{Y}$ is 1
by \eqref{QYZ1} and \eqref{QYZ2}.

Let $Y$ be a proper Young wall in ${\mathcal Z}(\Lambda)$. Let $L$
be a ladder such that there exists at least one block in $Y$ whose
coordinate is in $L$. We denote $Y\cap L$ by the set of all the
blocks in $Y$ whose coordinates are in $L$. Suppose that there are
$r$ many $i$-blocks in $Y\cap L$ for some $r\geq 0$ and $i\in I$.
Move these $i$-blocks to the first $r$ many $i$-slots in $L$ from
the bottom. Repeat this procedure ladder by ladder until no block
can be moved downward along a ladder. Then we obtain another
proper Young wall $Y^{R}$, which we call the {\it reduced form of
$Y$}. By definition, $Y^{R}$ is a reduced proper Young wall and
$|Y^R|\unrhd|Y|$, where the equality holds if and only if $Y$ is
reduced.

\begin{ex}{\rm \mbox{}
\vskip 5mm
\begin{center}
If \hskip 3mm $Y=$ \raisebox{-0.4\height}{\begin{texdraw} \drawdim
mm \setunitscale 0.5 \fontsize{7}{7}\selectfont \textref h:C v:C

\move(0 0)\recs  \move(10 0)\recs \move(-10 0)\recs

\move(-10 0)\bsegment

\move(0 0)\lvec(10 0)\lvec(10 10)\lvec(0 10)\lvec(0 0) \move(0
5)\lvec(10 5)\htext(5 2.5){$1$}\htext(5 7.5){$1$}

\esegment

\move(0 0)\bsegment

\move(0 0)\lvec(10 0)\lvec(10 10)\lvec(0 10)\lvec(0 0) \move(0
5)\lvec(10 5)\htext(5 2.5){$1$}\htext(5 7.5){$1$}

\move(0 10)\lvec(10 10)\lvec(10 20)\lvec(0 20)\lvec(0 10)\move(0
10)\lvec(10 20)\htext(2.5 17.5){$2$}\htext(7.5 12.5){$0$}

\move(0 20)\lvec(10 20)\lvec(10 30)\lvec(0 30)\lvec(0 20)\move(0
25)\lvec(10 25)\htext(5 22.5){$1$}\htext(5 27.5){$1$}

\move(0 30)\lvec(10 30)\lvec(10 40)\lvec(0 40)\lvec(0 30)\move(0
30)\lvec(10 40)\htext(2.5 37.5){$2$}\htext(7.5 32.5){$0$}

\esegment

\move(10 0)\bsegment

\move(0 0)\lvec(10 0)\lvec(10 10)\lvec(0 10)\lvec(0 0) \move(0
5)\lvec(10 5)\htext(5 2.5){$1$}\htext(5 7.5){$1$}

\move(0 10)\lvec(10 10)\lvec(10 20)\lvec(0 20)\lvec(0 10)\move(0
10)\lvec(10 20)\htext(2.5 17.5){$0$}\htext(7.5 12.5){$2$}

\move(0 20)\lvec(10 20)\lvec(10 30)\lvec(0 30)\lvec(0 20)\move(0
25)\lvec(10 25)\htext(5 22.5){$1$}\htext(5 27.5){$1$}

\move(0 30)\lvec(10 30)\lvec(10 40)\lvec(0 40)\lvec(0 30)\move(0
30)\lvec(10 40)\htext(2.5 37.5){$0$}\htext(7.5 32.5){$2$}

\move(0 40)\lvec(10 40)\lvec(10 50)\lvec(0 50)\lvec(0 40)\move(0
45)\lvec(10 45)\htext(5 42.5){$1$}\htext(5 47.5){$1$}

\move(0 50)\lvec(10 50)\lvec(10 60)\lvec(0 60)\lvec(0 50)\move(0
50)\lvec(10 60)\htext(2.5 57.5){$0$}\htext(7.5 52.5){$2$}

\move(0 60)\lvec(10 60)\lvec(10 70)\lvec(0 70)\lvec(0 60)\move(0
65)\lvec(10 65)\htext(5 62.5){$1$}\htext(5 67.5){$1$}

\move(0 70)\lvec(10 70)\lvec(10 80)\lvec(0 80)\lvec(0 70)\move(0
70)\lvec(10 80)\htext(2.5 77.5){$0$}\htext(7.5 72.5){$2$}

\esegment

\end{texdraw}}\hskip 3mm ,
\hskip .5cm then \hskip 5mm$Y^R=$
\raisebox{-0.6\height}{\begin{texdraw} \drawdim mm \setunitscale
0.5 \fontsize{7}{7}\selectfont \textref h:C v:C

\move(0 0)\recs  \move(10 0)\recs \move(-10 0)\recs \move(-20
0)\recs

\move(-20 0)\bsegment

\move(0 0)\lvec(10 0)\lvec(10 10)\lvec(0 10)\lvec(0 0) \move(0
5)\lvec(10 5)\htext(5 2.5){$1$}\htext(5 7.5){$1$}

\move(0 10)\lvec(10 10)\lvec(10 20)\lvec(0 20)\lvec(0 10)\move(0
10)\lvec(10 20)\htext(2.5 17.5){$0$}\htext(7.5 12.5){$2$}

\esegment

\move(-10 0)\bsegment

\move(0 0)\lvec(10 0)\lvec(10 10)\lvec(0 10)\lvec(0 0) \move(0
5)\lvec(10 5)\htext(5 2.5){$1$}\htext(5 7.5){$1$}

\move(0 10)\lvec(10 10)\lvec(10 20)\lvec(0 20)\lvec(0 10)\move(0
10)\lvec(10 20)\htext(2.5 17.5){$0$}\htext(7.5 12.5){$2$}

\move(0 20)\lvec(10 20)\lvec(10 25)\lvec(0 25)\lvec(0 20)\htext(5
22.5){$1$}

\esegment

\move(0 0)\bsegment

\move(0 0)\lvec(10 0)\lvec(10 10)\lvec(0 10)\lvec(0 0) \move(0
5)\lvec(10 5)\htext(5 2.5){$1$}\htext(5 7.5){$1$}

\move(0 10)\lvec(10 10)\lvec(10 20)\lvec(0 20)\lvec(0 10)\move(0
10)\lvec(10 20)\htext(2.5 17.5){$2$}\htext(7.5 12.5){$0$}

\move(0 20)\lvec(10 20)\lvec(10 30)\lvec(0 30)\lvec(0 20)\move(0
25)\lvec(10 25)\htext(5 22.5){$1$}\htext(5 27.5){$1$}

\move(0 30)\lvec(10 30)\lvec(10 40)\lvec(0 40)\lvec(0 30)\move(0
30)\lvec(10 40)\htext(2.5 37.5){$2$}\htext(7.5 32.5){$0$}

\esegment

\move(10 0)\bsegment

\move(0 0)\lvec(10 0)\lvec(10 10)\lvec(0 10)\lvec(0 0) \move(0
5)\lvec(10 5)\htext(5 2.5){$1$}\htext(5 7.5){$1$}

\move(0 10)\lvec(10 10)\lvec(10 20)\lvec(0 20)\lvec(0 10)\move(0
10)\lvec(10 20)\htext(2.5 17.5){$0$}\htext(7.5 12.5){$2$}

\move(0 20)\lvec(10 20)\lvec(10 30)\lvec(0 30)\lvec(0 20)\move(0
25)\lvec(10 25)\htext(5 22.5){$1$}\htext(5 27.5){$1$}

\move(0 30)\lvec(10 30)\lvec(10 40)\lvec(0 40)\lvec(0 30)\move(0
30)\lvec(10 40)\htext(2.5 37.5){$0$}\htext(7.5 32.5){$2$}

\move(0 40)\lvec(10 40)\lvec(10 50)\lvec(0 50)\lvec(0 40)\move(0
45)\lvec(10 45)\htext(5 42.5){$1$}\htext(5 47.5){$1$}

\esegment

\end{texdraw}}\hskip 3mm .

\end{center}
\vskip 5mm}
\end{ex}

Now, we are in a position to describe the algorithm. First, we
will construct an $\mathbb{A}$-basis $A(\Lambda)=\{\,A(Y)\,|\,Y\in
{\mathcal Y}(\Lambda)\,\}$ for $V(\Lambda)^{\mathbb{A}}$, which is
invariant under the involution $-$.

For $Y\in{\mathcal Y}(\Lambda)$, there exists a unique sequence of
reduced proper Young walls $\{Y_{k}\}_{k=0}^N$ such that
$Y_{0}=Y$, $Y_{k+1}=\overline{Y_{k}}$ ($0\leq k < N$),
$Y_{N}=Y_{\Lambda}$. Suppose that ${\rm wt}(Y_{k})={\rm
wt}(Y_{k-1})+r_k\alpha_{i_k}$ ($1\leq k\leq N$). We define
\begin{equation}
A(Y)=f_{i_1}^{(r_1)}\cdots f_{i_N}^{(r_N)}Y_{\Lambda}\in
V(\Lambda)^{\mathbb{A}}.
\end{equation}
It is clear that $\overline{A(Y)}=A(Y)$. We write
\begin{equation}
A(Y)=\sum_{Z\in{\mathcal Z}(\Lambda)} A_{Y,Z}(q)Z,
\end{equation}
where $A_{Y,Z}(q)\in\mathbb{Z}[q,q^{-1}]$ (see \eqref{QYZ1} and
\eqref{QYZ2}).
\begin{prop}\label {AYZ} Let $Y$ be a reduced proper Young wall.
For $Z\in{\mathcal Z}(\Lambda)$, we have
\begin{itemize}
\item[{\rm (a)}] if $A_{Y,Z}(q)\neq 0$, then $|Y|\unrhd |Z^R|$ and
${\rm wt}(Y)={\rm wt}(Z)$;

\item[{\rm (b)}] if $A_{Y,Z}(q)\neq 0$ and $|Y|=|Z|$, then $Y=Z$ and
$A_{Y,Y}(q)=1$;

\item[{\rm (c)}] $A(\Lambda)$ is a $\mathbb{Q}(q)$-basis of
$V(\Lambda)$.
\end{itemize}
\end{prop}
To prove this, we need the following technical lemma.

\begin{lem}\label {ord}{\rm \cite{KK}}
Let $Y\in{\mathcal Y}(\Lambda)$ and $Z\in{\mathcal Z}(\Lambda)$ be
such that $|\overline{Y}|\unrhd |Z^{R}|$. Suppose that ${\rm
wt}(\overline{Y})={\rm wt}(Y)+r\alpha_i$ for some $i\in I$ and
$r\geq 1$. Then, for each $W\in{\mathcal Z}(\Lambda)$ occurring in
the expansion of $f_i^{(r)}Z$, we have

{\rm (a)} $|Y|\unrhd|W^{R}|$;

{\rm (b)} if $|Y|=|W|$, then $|\overline{Y}|=|Z|$ and $Z$ is
reduced;

{\rm (c)} if $|Y|=|W|$ and $\overline{Y}=Z$, then $Y=W$.\qed
\end{lem}

{\bf Proof of Proposition \ref{AYZ}}. We will use induction on
$l$, the number of blocks in $Y$ which have been added to
$Y_{\Lambda}$. If $l=1$, it is clear. Suppose that $l>1$, and (a)
and (b) hold for $l'<l$. If $A(Y)=f_{i_1}^{(r_1)}\cdots
f_{i_N}^{(r_N)}Y_{\Lambda}$ for some $N\geq 1$, then we have
\begin{equation}
\begin{split}
A(Y)&=f_{i_1}^{(r_1)}A(\overline{Y})
=\sum_{|\overline{Y}|\unrhd|Z^R|}A_{\overline{Y},Z}(q)f_{i_1}^{(r_1)}Z
\\
&=\sum_{|\overline{Y}|\unrhd|Z^R|}A_{\overline{Y},Z}(q) \left (
\sum_{|Y|\unrhd|W^R|}Q_{Z,W}(q)W\right ) \hskip 1cm\text{by Lemma
\ref{ord} (a)}
\\
&=\sum_{|Y|\unrhd|W^R|}\left
(\sum_{|\overline{Y}|\unrhd|Z^R|}A_{\overline{Y},Z}(q)Q_{Z,W}(q)
\right )W.
\end{split}
\end{equation}
We have
\begin{equation}\label{AYW}
A_{Y,W}(q)=\sum_{|\overline{Y}|\unrhd|Z^R|}A_{\overline{Y},Z}(q)Q_{Z,W}(q)
\in\mathbb{Z}[q,q^{-1}],
\end{equation}
and $A_{Y,W}(q)=0$ unless $|Y|\unrhd |W^R|$ and ${\rm wt}(Y)={\rm
wt}(W)$. If $A_{Y,W}(q)\neq 0$ and $|Y|=|W|$, then Lemma \ref{ord}
(b) and \eqref{AYW} imply that $|Z|=|\overline{Y}|$ and
$A_{\overline{Y},Z}(q)\neq 0$. Hence, $Z=\overline{Y}$ by
induction hypothesis. Finally, we have $Y=W$ by Lemma \ref{ord}
(c), and then
$A_{Y,Y}(q)=A_{\overline{Y},\overline{Y}}(q)Q_{\overline{Y},Y}(q)=1$,
which completes the induction argument.

By (a) and (b), $A(\Lambda)$ is linearly independent over
$\mathbb{Q}(q)$. Since $\dim V(\Lambda)_{\lambda}=|{\mathcal
Y}(\Lambda)_{\lambda}|$ for all $\lambda\leq \Lambda$,
$A(\Lambda)$ is a $\mathbb{Q}(q)$-basis of $V(\Lambda)$. This
proves (c). \qed\vskip 3mm

For $Y\in {\mathcal Y}(\Lambda)$, we write
\begin{equation}
G(Y)=\sum_{Z\in{\mathcal Z}(\Lambda)} G_{Y,Z}(q)Z\in
V(\Lambda)^{\mathbb{A}}
\end{equation}
where $G_{Y,Z}\in \mathbb{Q}(q)$. Note that the coefficient
$G_{Y,Z}(q)$ satisfies
\begin{itemize}
\item[(i)] $G_{Y,Z}(q)\in\mathbb{Q}[q]$,
\item[(ii)] $G_{Y,Z}(q)\in q\mathbb{Q}[q]$ unless $Y=Z$,
\item[(iii)] $G_{Y,Y}(q)=1$.
\end{itemize}
Consider the following three matrices indexed by ${\mathcal
Y}(\Lambda)$ where the indices are decreasing with respect to the
total ordering $>$;
\begin{equation}
G=(G_{Y,Z}(q)),\quad H=(H_{Y,W}(q)), \quad A=(A_{W,Z}(q)),
\end{equation}
where $H$ is the transition matrix from $A(\Lambda)$ to
$G(\Lambda)$ as $\mathbb{Q}(q)$-bases of $V(\Lambda)$. We have
$G=HA$. By the $-$ invariance of $G(\Lambda)$ and $A(\Lambda)$,
$H$ is also invariant under the involution $-$. Since $A$ is a
unipotent matrix with entries in $\mathbb{A}$, $H$ is a matrix
with entries in $\mathbb{A}$, which yields:
\begin{prop}
$A(\Lambda)$ is an $\mathbb{A}$-basis of
$V(\Lambda)^{\mathbb{A}}$.\qed
\end{prop}

Also by the argument in \cite{LLT}, $H$ is a unipotent matrix.
Hence for each $Y\in{\mathcal Y}(\Lambda)_{\lambda}$ ($\lambda\leq
\Lambda$), $A(Y)$ can be expressed uniquely as follows;
\begin{equation}\label {H'YZ}
A(Y)=G(Y)+\sum_{\substack{Z\in {\mathcal Y}(\Lambda)_{\lambda}
\\ Y>Z }}\gamma_{Y,Z}(q)G(Z),
\end{equation}
for some $\gamma_{Y,Z}(q)\in \mathbb{Q}[q,q^{-1}]$ such that
$\gamma_{Y,Z}(q)=\gamma_{Y,Z}(q^{-1})$. If $Y$ is the minimal
element in ${\mathcal Y}(\Lambda)_{\lambda}$, then $A(Y)=G(Y)$.
Suppose that $Y$ is not minimal and $G(Y')$ are given for
$Y'\in{\mathcal Y}(\Lambda)_{\lambda}$ such that $Y'<Y$. Then
$\gamma_{Y,Y'}(q)$ are determined inductively as follows;
\begin{itemize}
\item[(1)] if $Y'$ is the maximal one such that $Y>Y'$ and
$A_{Y,Y'}(q)=\sum_{i=-r'}^ra_iq^{-i}$, then
$\gamma_{Y,Y'}(q)=\sum_{i=1}^{r}a_i(q^i +q^{-i}) + a_0$.

\item[(2)] if the coefficient of $Y'$ in $A(Y)-\sum_{Y>Z>Y'}
\gamma_{Y,Z}(q)G(Z)$ is given by $\sum_{i=-r'}^ra_iq^{-i}$, then
$\gamma_{Y,Y'}(q)=\sum_{i=1}^{r}a_i(q^i +q^{-i}) + a_0$.
\end{itemize}
To summarize, we have
\begin{thm} For a reduced proper Young wall $Y\in{\mathcal Y}(\Lambda)_{\lambda}$ {\rm ($\lambda\leq \Lambda$)}, the
corresponding global basis element is of the following form;
\begin{equation}\label {G(Y)}
G(Y)=Y+\sum_{\substack{Z\in{\mathcal Z}(\Lambda)_{\lambda} \\
|Y|\rhd|Z^R|}}G_{Y,Z}(q)Z,
\end{equation}
where  $G_{Y,Z}(q)\in q\mathbb{Z}[q]$ for $Y\neq Z$.\qed
\end{thm}

\begin{ex}{\rm
In the following, we list $G(Y)$, where $Y$ is the reduced proper
Young wall in ${\mathcal Y}(\Lambda_1)$ (Example \ref{graph}). Set
$k=\sum_{i\in I}k_i$ where ${\rm wt}(Y)=\Lambda_1-\sum_{i\in
I}k_i\alpha_i$ \vskip 5mm

\noindent(1) $k=1$\vskip 3mm

$G(~\raisebox{-0.3\height}{\begin{texdraw} \drawdim mm
\fontsize{7}{7}\selectfont \textref h:C v:C \setunitscale 0.5

\move(0 0)\recs \move(0 0) \bsegment \move(0 0)\lvec(10 0)\lvec(10
10)\lvec(0 10)\lvec(0 0) \move(10 5)\lvec(0 5) \htext(5
2.5){$1$}\htext(5 7.5){$1$} \esegment
\end{texdraw}}~)=f_1Y_{\Lambda_1} =\raisebox{-0.3\height}{\begin{texdraw} \drawdim mm
\fontsize{7}{7}\selectfont \textref h:C v:C \setunitscale 0.5

\move(0 0)\recs \move(0 0) \bsegment \move(0 0)\lvec(10 0)\lvec(10
10)\lvec(0 10)\lvec(0 0) \move(10 5)\lvec(0 5) \htext(5
2.5){$1$}\htext(5 7.5){$1$} \esegment
\end{texdraw}}$\vskip 5mm

\noindent(2) $k=2$\vskip 3mm

$G(~\raisebox{-0.3\height}{\begin{texdraw} \drawdim mm
\fontsize{7}{7}\selectfont \textref h:C v:C \setunitscale 0.5

\move(0 0)\recs \move(0 0)  \bsegment \move(0 0)\lvec(10
0)\lvec(10 20)\lvec(0 20)\lvec(0 0) \move(0 5)\lvec(10 5) \move(10
10)\lvec(0 10)\lvec(10 20) \htext(5 2.5){$1$}\htext(5
7.5){$1$}\htext(2.5 17.5){$0$} \esegment
\end{texdraw}}~)=f_0f_1Y_{\Lambda_1} =
\raisebox{-0.3\height}{\begin{texdraw} \drawdim mm
\fontsize{7}{7}\selectfont \textref h:C v:C \setunitscale 0.5

\move(0 0)\recs \move(0 0)  \bsegment \move(0 0)\lvec(10
0)\lvec(10 20)\lvec(0 20)\lvec(0 0) \move(0 5)\lvec(10 5) \move(10
10)\lvec(0 10)\lvec(10 20) \htext(5 2.5){$1$}\htext(5
7.5){$1$}\htext(2.5 17.5){$0$} \esegment
\end{texdraw}}$\vskip 5mm

$G(\raisebox{-0.3\height}{
\begin{texdraw}
\drawdim mm \fontsize{7}{7}\selectfont \textref h:C v:C
\setunitscale 0.5 \move(0 0)\recs \move(0 0) \bsegment \move(0
0)\lvec(10 0)\lvec(10 20)\lvec(0 10)\lvec(0 0) \move(0 5)\lvec(10
5) \move(10 10)\lvec(0 10)\lvec(10 20) \htext(5 2.5){$1$}\htext(5
7.5){$1$}\htext(7.5 12.5){$2$} \esegment
\end{texdraw}}~)=f_2f_1Y_{\Lambda_1} =
\raisebox{-0.3\height}{\begin{texdraw} \drawdim mm
\fontsize{7}{7}\selectfont \textref h:C v:C \setunitscale 0.5
\move(0 0)\recs \move(0 0) \bsegment \move(0 0)\lvec(10 0)\lvec(10
20)\lvec(0 10)\lvec(0 0) \move(0 5)\lvec(10 5) \move(10 10)\lvec(0
10)\lvec(10 20) \htext(5 2.5){$1$}\htext(5 7.5){$1$}\htext(7.5
12.5){$2$} \esegment
\end{texdraw}}$ \vskip 5mm

\noindent(3) $k=3$\vskip 3mm

$G(\raisebox{-0.3\height}{
\begin{texdraw}
\drawdim mm \fontsize{7}{7}\selectfont \textref h:C v:C
\setunitscale 0.5 \move(0 0)\recs \move(0 0) \bsegment \move(0
0)\lvec(10 0)\lvec(10 20)\lvec(0 20)\lvec(0 0) \move(0 5)\lvec(10
5) \move(10 10)\lvec(0 10)\lvec(10 20) \htext(5 2.5){$1$}\htext(5
7.5){$1$}\htext(7.5 12.5){$2$}\htext(2.5 17.5){$0$} \esegment
\end{texdraw}}~)=f_0f_2f_1Y_{\Lambda_1} =
\raisebox{-0.3\height}{\begin{texdraw} \drawdim mm
\fontsize{7}{7}\selectfont \textref h:C v:C \setunitscale 0.5
\move(0 0)\recs \move(0 0) \bsegment \move(0 0)\lvec(10 0)\lvec(10
20)\lvec(0 20)\lvec(0 0) \move(0 5)\lvec(10 5) \move(10 10)\lvec(0
10)\lvec(10 20) \htext(5 2.5){$1$}\htext(5 7.5){$1$}\htext(7.5
12.5){$2$}\htext(2.5 17.5){$0$} \esegment
\end{texdraw}}$ \vskip 5mm

$G(\raisebox{-0.3\height}{
\begin{texdraw}
\drawdim mm \fontsize{7}{7}\selectfont \textref h:C v:C
\setunitscale 0.5 \move(0 0)\recs \move(-10 0)\recs \move(0 0)
\bsegment \move(10 10)\lvec(-10 10)\lvec(-10 0)\lvec(10 0)\lvec(10
20)\lvec(0 20) \lvec(0 0) \move(-10 5)\lvec(10 5) \move(0
10)\lvec(10 20) \htext(5 2.5){$1$}\htext(5 7.5){$1$}\htext(2.5
17.5){$0$} \htext(-5 2.5){$1$}\htext(-5 7.5){$1$} \esegment
\end{texdraw}}~)=f_1f_0f_1Y_{\Lambda_1} =
\raisebox{-0.3\height}{\begin{texdraw} \drawdim mm
\fontsize{7}{7}\selectfont \textref h:C v:C \setunitscale 0.5
\move(0 0)\recs \move(-10 0)\recs  \move(0 0) \bsegment \move(10
10)\lvec(-10 10)\lvec(-10 0)\lvec(10 0)\lvec(10 20)\lvec(0 20)
\lvec(0 0) \move(-10 5)\lvec(10 5) \move(0 10)\lvec(10 20)
\htext(5 2.5){$1$}\htext(5 7.5){$1$}\htext(2.5 17.5){$0$}
\htext(-5 2.5){$1$}\htext(-5 7.5){$1$} \esegment
\end{texdraw}}$ \vskip 5mm

$G(\raisebox{-0.3\height}{
\begin{texdraw}
\drawdim mm \fontsize{7}{7}\selectfont \textref h:C v:C
\setunitscale 0.5 \move(0 0)\recs \move(-10 0)\recs \move(0 0)
\bsegment \move(10 10)\lvec(-10 10)\lvec(-10 0)\lvec(10 0)\lvec(10
20)\lvec(0 10) \lvec(0 0) \move(-10 5)\lvec(10 5) \move(0
10)\lvec(10 20) \htext(5 2.5){$1$}\htext(5 7.5){$1$}\htext(7.5
12.5){$2$} \htext(-5 2.5){$1$}\htext(-5 7.5){$1$} \esegment
\end{texdraw}}~)=f_1f_2f_1Y_{\Lambda_1} =
\raisebox{-0.3\height}{\begin{texdraw} \drawdim mm
\fontsize{7}{7}\selectfont \textref h:C v:C \setunitscale 0.5
\move(0 0)\recs \move(-10 0)\recs  \move(0 0) \bsegment \move(10
10)\lvec(-10 10)\lvec(-10 0)\lvec(10 0)\lvec(10 20)\lvec(0 10)
\lvec(0 0) \move(-10 5)\lvec(10 5) \move(0 10)\lvec(10 20)
\htext(5 2.5){$1$}\htext(5 7.5){$1$}\htext(7.5 12.5){$2$}
\htext(-5 2.5){$1$}\htext(-5 7.5){$1$} \esegment
\end{texdraw}}$ \vskip 5mm

\noindent (4) $k=4$\vskip 3mm

$G(\raisebox{-0.3\height}{
\begin{texdraw}
\drawdim mm \fontsize{7}{7}\selectfont \textref h:C v:C
\setunitscale 0.5 \move(0 0)\recs \move(-10 0)\recs \move(0 0)
\bsegment \move(10 10)\lvec(-10 10)\lvec(-10 0)\lvec(10 0)\lvec(10
20)\lvec(0 20) \lvec(0 0) \move(-10 5)\lvec(10 5) \move(0
10)\lvec(10 20) \htext(5 2.5){$1$}\htext(5 7.5){$1$}\htext(2.5
17.5){$0$}\htext(7.5 12.5){$2$} \htext(-5 2.5){$1$}\htext(-5
7.5){$1$} \esegment
\end{texdraw}}~)=f_1f_0f_2f_1Y_{\Lambda_1} =
\raisebox{-0.3\height}{\begin{texdraw} \drawdim mm
\fontsize{7}{7}\selectfont \textref h:C v:C \setunitscale 0.5
\move(0 0)\recs \move(-10 0)\recs  \move(0 0) \bsegment \move(10
10)\lvec(-10 10)\lvec(-10 0)\lvec(10 0)\lvec(10 20)\lvec(0 20)
\lvec(0 0) \move(-10 5)\lvec(10 5) \move(0 10)\lvec(10 20)
\htext(5 2.5){$1$}\htext(5 7.5){$1$}\htext(2.5
17.5){$0$}\htext(7.5 12.5){$2$} \htext(-5 2.5){$1$}\htext(-5
7.5){$1$} \esegment
\end{texdraw}}\hskip 3mm + \hskip 3mm q\hskip 2mm
\raisebox{-0.25\height}{\begin{texdraw} \drawdim mm
\fontsize{7}{7}\selectfont \textref h:C v:C \setunitscale 0.5
\move(0 0)\recs \move(0 0) \bsegment \move(10 10)\lvec(0
10)\lvec(0 0)\lvec(10 0)\lvec(10 20)\lvec(0 20) \lvec(0 0) \move(0
5)\lvec(10 5) \move(0 10)\lvec(10 20)\move(0 20)\lvec(0 25)\lvec
(10 25)\lvec(10 20) \htext(5 2.5){$1$}\htext(5 7.5){$1$}\htext(5
22.5){$1$}\htext(2.5 17.5){$0$}\htext(7.5 12.5){$2$}  \esegment
\end{texdraw}}
$ \vskip 5mm

$G(\raisebox{-0.3\height}{
\begin{texdraw}
\drawdim mm \fontsize{7}{7}\selectfont \textref h:C v:C
\setunitscale 0.5 \move(0 0)\recs \move(-10 0)\recs \move(0 0)
\bsegment \move(0 20)\lvec(0 0)\lvec(10 0)\lvec(10 20)\lvec(-10
20)\lvec(-10 0)\lvec(0 0) \move(-10 5)\lvec(10 5) \move(10
10)\lvec(0 10)\lvec(10 20) \move(0 10)\lvec(-10 10)\lvec(0 20)
\htext(5 2.5){$1$}\htext(5 7.5){$1$} \htext(-5 2.5){$1$}\htext(-5
7.5){$1$} \htext(-7.5 17.5){$2$}\htext(2.5 17.5){$0$} \esegment
\end{texdraw}}~)=f_2f_1f_0f_1Y_{\Lambda_1} =
\raisebox{-0.3\height}{\begin{texdraw} \drawdim mm
\fontsize{7}{7}\selectfont \textref h:C v:C \setunitscale 0.5
\move(0 0)\recs \move(-10 0)\recs  \move(0 0) \bsegment \move(0
20)\lvec(0 0)\lvec(10 0)\lvec(10 20)\lvec(-10 20)\lvec(-10
0)\lvec(0 0) \move(-10 5)\lvec(10 5) \move(10 10)\lvec(0
10)\lvec(10 20) \move(0 10)\lvec(-10 10)\lvec(0 20) \htext(5
2.5){$1$}\htext(5 7.5){$1$} \htext(-5 2.5){$1$}\htext(-5 7.5){$1$}
\htext(-7.5 17.5){$2$}\htext(2.5 17.5){$0$} \esegment
\end{texdraw}}\hskip 3mm + \hskip 3mm q^2 \hskip 2mm
\raisebox{-0.25\height}{\begin{texdraw} \drawdim mm
\fontsize{7}{7}\selectfont \textref h:C v:C \setunitscale 0.5
\move(0 0)\recs \move(-10 0)\recs \move(0 0) \bsegment \move(10
10)\lvec(-10 10)\lvec(-10 0)\lvec(10 0)\lvec(10 20)\lvec(0 20)
\lvec(0 0) \move(-10 5)\lvec(10 5) \move(0 10)\lvec(10 20)
\htext(5 2.5){$1$}\htext(5 7.5){$1$}\htext(2.5
17.5){$0$}\htext(7.5 12.5){$2$} \htext(-5 2.5){$1$}\htext(-5
7.5){$1$} \esegment
\end{texdraw}}
$ \vskip 5mm

$G(\raisebox{-0.3\height}{
\begin{texdraw}
\drawdim mm \fontsize{7}{7}\selectfont \textref h:C v:C
\setunitscale 0.5 \move(0 0)\recs \move(-10 0)\recs \move(0 0)
\bsegment \move(0 20)\lvec(0 0)\lvec(10 0)\lvec(10 20)\move(-10
10)\lvec(-10 0)\lvec(0 0) \move(-10 5)\lvec(10 5) \move(10
10)\lvec(0 10)\lvec(10 20) \move(0 10)\lvec(-10 10)\lvec(0 20)
\htext(5 2.5){$1$}\htext(5 7.5){$1$} \htext(-5 2.5){$1$}\htext(-5
7.5){$1$} \htext(-2.5 12.5){$0$}\htext(7.5 12.5){$2$} \esegment
\end{texdraw}}~)=f_0f_1f_2f_1Y_{\Lambda_1} =
\raisebox{-0.3\height}{\begin{texdraw} \drawdim mm
\fontsize{7}{7}\selectfont \textref h:C v:C \setunitscale 0.5
\move(0 0)\recs \move(-10 0)\recs  \move(0 0) \bsegment \move(0
20)\lvec(0 0)\lvec(10 0)\lvec(10 20)\move(-10 10)\lvec(-10
0)\lvec(0 0) \move(-10 5)\lvec(10 5) \move(10 10)\lvec(0
10)\lvec(10 20) \move(0 10)\lvec(-10 10)\lvec(0 20) \htext(5
2.5){$1$}\htext(5 7.5){$1$} \htext(-5 2.5){$1$}\htext(-5 7.5){$1$}
\htext(-2.5 12.5){$0$}\htext(7.5 12.5){$2$} \esegment
\end{texdraw}}\hskip 3mm + \hskip 3mm q^2 \hskip 2mm
\raisebox{-0.25\height}{\begin{texdraw} \drawdim mm
\fontsize{7}{7}\selectfont \textref h:C v:C \setunitscale 0.5
\move(0 0)\recs \move(-10 0)\recs \move(0 0) \bsegment \move(10
10)\lvec(-10 10)\lvec(-10 0)\lvec(10 0)\lvec(10 20)\lvec(0 20)
\lvec(0 0) \move(-10 5)\lvec(10 5) \move(0 10)\lvec(10 20)
\htext(5 2.5){$1$}\htext(5 7.5){$1$}\htext(2.5
17.5){$0$}\htext(7.5 12.5){$2$} \htext(-5 2.5){$1$}\htext(-5
7.5){$1$} \esegment
\end{texdraw}}
$ \vskip 5mm

\noindent(5) $k=5$\vskip 3mm
$G(\raisebox{-0.3\height}{
\begin{texdraw}
\drawdim mm \fontsize{7}{7}\selectfont \textref h:C v:C
\setunitscale 0.5 \move(0 0)\recs \move(-10 0)\recs \move(0 0)
\bsegment \move(10 10)\lvec(-10 10)\lvec(-10 0)\lvec(10 0)\lvec(10
20)\lvec(0 20) \lvec(0 0) \move(-10 5)\lvec(10 5) \move(0
10)\lvec(10 20) \htext(5 2.5){$1$}\htext(5 7.5){$1$}\htext(2.5
17.5){$0$}\htext(7.5 12.5){$2$} \htext(-5 2.5){$1$}\htext(-5
7.5){$1$} \move(10 20)\lvec(10 25)\lvec(0 25)\lvec(0 20) \htext(5
22.5){$1$}\esegment
\end{texdraw}}~)=f_1^{(2)}f_0f_2f_1Y_{\Lambda_1} =
\raisebox{-0.3\height}{\begin{texdraw} \drawdim mm
\fontsize{7}{7}\selectfont \textref h:C v:C \setunitscale 0.5
\move(0 0)\recs \move(-10 0)\recs  \move(0 0) \bsegment \move(10
10)\lvec(-10 10)\lvec(-10 0)\lvec(10 0)\lvec(10 20)\lvec(0 20)
\lvec(0 0) \move(-10 5)\lvec(10 5) \move(0 10)\lvec(10 20)
\htext(5 2.5){$1$}\htext(5 7.5){$1$}\htext(2.5
17.5){$0$}\htext(7.5 12.5){$2$} \htext(-5 2.5){$1$}\htext(-5
7.5){$1$} \move(10 20)\lvec(10 25)\lvec(0 25)\lvec(0 20) \htext(5
22.5){$1$} \esegment
\end{texdraw}}\hskip 3mm + \hskip 3mm q^2 \hskip 2mm
\raisebox{-0.25\height}{\begin{texdraw} \drawdim mm
\fontsize{7}{7}\selectfont \textref h:C v:C \setunitscale 0.5
\move(0 0)\recs \move(0 0) \bsegment \move(10 10)\lvec(0
10)\lvec(0 0)\lvec(10 0)\lvec(10 20)\lvec(0 20) \lvec(0 0) \move(0
5)\lvec(10 5) \move(0 10)\lvec(10 20)\move(0 20)\lvec(0 25)\lvec
(10 25)\lvec(10 20)\move(0 25)\lvec(0 30)\lvec(10 30)\lvec(10 25)
\htext(5 2.5){$1$}\htext(5 7.5){$1$}\htext(5 22.5){$1$}\htext(5
27.5){$1$}\htext(2.5 17.5){$0$}\htext(7.5 12.5){$2$}  \esegment
\end{texdraw}}$\vskip 5mm

$G(\raisebox{-0.3\height}{
\begin{texdraw}
\drawdim mm \fontsize{7}{7}\selectfont \textref h:C v:C
\setunitscale 0.5 \move(0 0)\recs \move(-10 0)\recs \move(0 0)
\bsegment \move(0 20)\lvec(0 0)\lvec(10 0)\lvec(10 20)\lvec(0 20
)\move(-10 10)\lvec(-10 0)\lvec(0 0) \move(-10 5)\lvec(10 5)
\move(10 10)\lvec(0 10)\lvec(10 20) \move(0 10)\lvec(-10
10)\lvec(0 20) \htext(5 2.5){$1$}\htext(5 7.5){$1$} \htext(-5
2.5){$1$}\htext(-5 7.5){$1$} \htext(-2.5 12.5){$0$}\htext(7.5
12.5){$2$}\htext(2.5 17.5){$0$} \esegment
\end{texdraw}}~)=f_0f_1f_0f_2f_1Y_{\Lambda_1} =
\raisebox{-0.3\height}{\begin{texdraw} \drawdim mm
\fontsize{7}{7}\selectfont \textref h:C v:C \setunitscale 0.5
\move(0 0)\recs \move(-10 0)\recs  \move(0 0) \bsegment \move(0
20)\lvec(0 0)\lvec(10 0)\lvec(10 20)\lvec(0 20 )\move(-10
10)\lvec(-10 0)\lvec(0 0) \move(-10 5)\lvec(10 5) \move(10
10)\lvec(0 10)\lvec(10 20) \move(0 10)\lvec(-10 10)\lvec(0 20)
\htext(5 2.5){$1$}\htext(5 7.5){$1$} \htext(-5 2.5){$1$}\htext(-5
7.5){$1$} \htext(-2.5 12.5){$0$}\htext(7.5 12.5){$2$}\htext(2.5
17.5){$0$} \esegment
\end{texdraw}}$\vskip 5mm

$G(\raisebox{-0.3\height}{
\begin{texdraw}
\drawdim mm \fontsize{7}{7}\selectfont \textref h:C v:C
\setunitscale 0.5 \move(0 0)\recs \move(-10 0)\recs \move(0 0)
\bsegment \move(0 20)\lvec(0 0)\lvec(10 0)\lvec(10 20)\lvec(-10
20)\lvec(-10 0)\lvec(0 0) \move(-10 5)\lvec(10 5) \move(10
10)\lvec(0 10)\lvec(10 20) \move(0 10)\lvec(-10 10)\lvec(0 20)
\htext(5 2.5){$1$}\htext(5 7.5){$1$} \htext(-5 2.5){$1$}\htext(-5
7.5){$1$} \htext(-7.5 17.5){$2$}\htext(2.5 17.5){$0$} \htext(7.5
12.5){$2$} \esegment
\end{texdraw}}~)=f_2f_1f_0f_2f_1Y_{\Lambda_1} =
\raisebox{-0.3\height}{\begin{texdraw} \drawdim mm
\fontsize{7}{7}\selectfont \textref h:C v:C \setunitscale 0.5
\move(0 0)\recs \move(-10 0)\recs  \move(0 0) \bsegment \move(0
20)\lvec(0 0)\lvec(10 0)\lvec(10 20)\lvec(-10 20)\lvec(-10
0)\lvec(0 0) \move(-10 5)\lvec(10 5) \move(10 10)\lvec(0
10)\lvec(10 20) \move(0 10)\lvec(-10 10)\lvec(0 20) \htext(5
2.5){$1$}\htext(5 7.5){$1$} \htext(-5 2.5){$1$}\htext(-5 7.5){$1$}
\htext(-7.5 17.5){$2$}\htext(2.5 17.5){$0$} \htext(7.5 12.5){$2$}
\esegment
\end{texdraw}}$\vskip 5mm

$G(\raisebox{-0.3\height}{
\begin{texdraw}
\drawdim mm \fontsize{7}{7}\selectfont \textref h:C v:C
\setunitscale 0.5 \move(0 0)\recs \move(-10 0)\recs \move(-20
0)\recs \move(0 0) \bsegment \move(0 20)\lvec(0 0)\lvec(10
0)\lvec(10 20)\lvec(-10 20)\lvec(-10 0)\lvec(0 0) \move(-10
5)\lvec(10 5) \move(10 10)\lvec(0 10)\lvec(10 20) \move(0
10)\lvec(-10 10)\lvec(0 20) \htext(5 2.5){$1$}\htext(5 7.5){$1$}
\htext(-5 2.5){$1$}\htext(-5 7.5){$1$} \htext(-7.5
17.5){$2$}\htext(2.5 17.5){$0$} \move(-10 0)\lvec(-20 0)\lvec(-20
10)\lvec(-10 10) \move(-10 5)\lvec(-20 5) \htext(-15 2.5){$1$}
\htext(-15 7.5){$1$} \esegment
\end{texdraw}}~)=f_1f_2f_1f_0f_1Y_{\Lambda_1} =
\raisebox{-0.3\height}{\begin{texdraw} \drawdim mm
\fontsize{7}{7}\selectfont \textref h:C v:C \setunitscale 0.5
\move(0 0)\recs \move(-10 0)\recs \move(-20 0)\recs  \move(0 0)
\bsegment \move(0 20)\lvec(0 0)\lvec(10 0)\lvec(10 20)\lvec(-10
20)\lvec(-10 0)\lvec(0 0) \move(-10 5)\lvec(10 5) \move(10
10)\lvec(0 10)\lvec(10 20) \move(0 10)\lvec(-10 10)\lvec(0 20)
\htext(5 2.5){$1$}\htext(5 7.5){$1$} \htext(-5 2.5){$1$}\htext(-5
7.5){$1$} \htext(-7.5 17.5){$2$}\htext(2.5 17.5){$0$} \move(-10
0)\lvec(-20 0)\lvec(-20 10)\lvec(-10 10) \move(-10 5)\lvec(-20 5)
\htext(-15 2.5){$1$} \htext(-15 7.5){$1$}\esegment
\end{texdraw}}\hskip 3mm + \hskip 3mm q \hskip 2mm
\raisebox{-0.25\height}{\begin{texdraw} \drawdim mm
\fontsize{7}{7}\selectfont \textref h:C v:C \setunitscale 0.5
\move(0 0)\recs \move(-10 0)\recs \move(0 0)  \bsegment \move(10
10)\lvec(-10 10)\lvec(-10 0)\lvec(10 0)\lvec(10 20)\lvec(0 20)
\lvec(0 0) \move(-10 5)\lvec(10 5) \move(0 10)\lvec(10 20)\move(0
20)\lvec(0 25)\lvec(10 25)\lvec(10 20)\htext(5 2.5){$1$}\htext(5
7.5){$1$}\htext(2.5 17.5){$0$}\htext(7.5 12.5){$2$}\htext(5
22.5){$1$} \htext(-5 2.5){$1$}\htext(-5 7.5){$1$} \esegment
\end{texdraw}}$\vskip 5mm

$G(\raisebox{-0.3\height}{
\begin{texdraw}
\drawdim mm \fontsize{7}{7}\selectfont \textref h:C v:C
\setunitscale 0.5 \move(0 0)\recs \move(-10 0)\recs \move(-20
0)\recs \move(0 0) \bsegment \move(0 20)\lvec(0 0)\lvec(10
0)\lvec(10 20)\move(-10 10)\lvec(-10 0)\lvec(0 0) \move(-10
5)\lvec(10 5) \move(10 10)\lvec(0 10)\lvec(10 20) \move(0
10)\lvec(-10 10)\lvec(0 20) \htext(5 2.5){$1$}\htext(5 7.5){$1$}
\htext(-5 2.5){$1$}\htext(-5 7.5){$1$} \htext(-2.5
12.5){$0$}\htext(7.5 12.5){$2$} \move(-10 0)\lvec(-20 0)\lvec(-20
10)\lvec(-10 10) \move(-10 5)\lvec(-20 5) \htext(-15 2.5){$1$}
\htext(-15 7.5){$1$} \esegment
\end{texdraw}}~)=f_1f_0f_1f_2f_1Y_{\Lambda_1} =
\raisebox{-0.3\height}{\begin{texdraw} \drawdim mm
\fontsize{7}{7}\selectfont \textref h:C v:C \setunitscale 0.5
\move(0 0)\recs \move(-10 0)\recs \move(-20 0)\recs  \move(0 0)
\bsegment \move(0 20)\lvec(0 0)\lvec(10 0)\lvec(10 20)\move(-10
10)\lvec(-10 0)\lvec(0 0) \move(-10 5)\lvec(10 5) \move(10
10)\lvec(0 10)\lvec(10 20) \move(0 10)\lvec(-10 10)\lvec(0 20)
\htext(5 2.5){$1$}\htext(5 7.5){$1$} \htext(-5 2.5){$1$}\htext(-5
7.5){$1$} \htext(-2.5 12.5){$0$}\htext(7.5 12.5){$2$} \move(-10
0)\lvec(-20 0)\lvec(-20 10)\lvec(-10 10) \move(-10 5)\lvec(-20 5)
\htext(-15 2.5){$1$} \htext(-15 7.5){$1$} \esegment
\end{texdraw}}\hskip 3mm + \hskip 3mm q \hskip 2mm
\raisebox{-0.25\height}{\begin{texdraw} \drawdim mm
\fontsize{7}{7}\selectfont \textref h:C v:C \setunitscale 0.5
\move(0 0)\recs \move(-10 0)\recs \move(0 0)  \bsegment \move(10
10)\lvec(-10 10)\lvec(-10 0)\lvec(10 0)\lvec(10 20)\lvec(0 20)
\lvec(0 0) \move(-10 5)\lvec(10 5) \move(0 10)\lvec(10 20)\move(0
20)\lvec(0 25)\lvec(10 25)\lvec(10 20)\htext(5 2.5){$1$}\htext(5
7.5){$1$}\htext(2.5 17.5){$0$}\htext(7.5 12.5){$2$}\htext(5
22.5){$1$} \htext(-5 2.5){$1$}\htext(-5 7.5){$1$} \esegment
\end{texdraw}}$\vskip 5mm
}
\end{ex}

\begin{ex}{\rm In the previous example, we have seen that $G(Y)=A(Y)$.
But this does not always hold for all $Y\in {\mathcal
Y}(\Lambda)$. Furthermore, the coefficient polynomial $G_{Y,Z}(q)$
do not always  have non-negative integral coefficients. Observe
that \vskip 5mm

$A(~\raisebox{-0.6\height}{\begin{texdraw} \drawdim mm
\setunitscale 0.5 \fontsize{7}{7}\selectfont \textref h:C v:C
\move(0 10)\tris \move(10 10)\tris \move(-20 10)\tris \move(-10
10)\tris \move(-20 0)\bsegment \move(0 10)\lvec(10 10)\lvec(10
20)\lvec(0 20)\lvec(0 10)\move(0 10)\lvec(10 20)\htext(2.5
17.5){$0$}\htext(7.5 12.5){$2$} \move(0 20)\lvec(10 20)\lvec(10
25)\lvec(0 25)\lvec(0 20)\htext(5 22.5){$1$} \esegment

\move(-10 0)\bsegment \move(0 10)\lvec(10 10)\lvec(10 20)\lvec(0
20)\lvec(0 10)\move(0 10)\lvec(10 20)\htext(2.5
17.5){$2$}\htext(7.5 12.5){$0$} \move(0 20)\lvec(10 20)\lvec(10
25)\lvec(0 25)\lvec(0 20)\htext(5 22.5){$1$} \esegment

\move(0 0)\bsegment \move(0 10)\lvec(10 10)\lvec(10 20)\lvec(0
20)\lvec(0 10)\move(0 10)\lvec(10 20)\htext(2.5
17.5){$0$}\htext(7.5 12.5){$2$} \move(0 20)\lvec(10 20)\lvec(10
25)\lvec(0 25)\lvec(0 20)\htext(5 22.5){$1$} \esegment

\move(10 0)\bsegment \move(0 10)\lvec(10 10)\lvec(10 20)\lvec(0
20)\lvec(0 10)\move(0 10)\lvec(10 20)\htext(2.5
17.5){$2$}\htext(7.5 12.5){$0$} \move(0 20)\lvec(10 20)\lvec(10
25)\lvec(0 25)\lvec(0 20)\htext(5 22.5){$1$} \esegment
\end{texdraw}}~)=f_1f_0f_1f_2f_1f_0f_1f_2 Y_{\Lambda_2}$

$=\hskip 2mm\raisebox{-0.6\height}{\begin{texdraw} \drawdim mm
\setunitscale 0.5 \fontsize{7}{7}\selectfont \textref h:C v:C
\move(0 10)\tris \move(10 10)\tris \move(-20 10)\tris \move(-10
10)\tris \move(-20 0)\bsegment \move(0 10)\lvec(10 10)\lvec(10
20)\lvec(0 20)\lvec(0 10)\move(0 10)\lvec(10 20)\htext(2.5
17.5){$0$}\htext(7.5 12.5){$2$} \move(0 20)\lvec(10 20)\lvec(10
25)\lvec(0 25)\lvec(0 20)\htext(5 22.5){$1$} \esegment

\move(-10 0)\bsegment \move(0 10)\lvec(10 10)\lvec(10 20)\lvec(0
20)\lvec(0 10)\move(0 10)\lvec(10 20)\htext(2.5
17.5){$2$}\htext(7.5 12.5){$0$} \move(0 20)\lvec(10 20)\lvec(10
25)\lvec(0 25)\lvec(0 20)\htext(5 22.5){$1$} \esegment

\move(0 0)\bsegment \move(0 10)\lvec(10 10)\lvec(10 20)\lvec(0
20)\lvec(0 10)\move(0 10)\lvec(10 20)\htext(2.5
17.5){$0$}\htext(7.5 12.5){$2$} \move(0 20)\lvec(10 20)\lvec(10
25)\lvec(0 25)\lvec(0 20)\htext(5 22.5){$1$} \esegment

\move(10 0)\bsegment \move(0 10)\lvec(10 10)\lvec(10 20)\lvec(0
20)\lvec(0 10)\move(0 10)\lvec(10 20)\htext(2.5
17.5){$2$}\htext(7.5 12.5){$0$} \move(0 20)\lvec(10 20)\lvec(10
25)\lvec(0 25)\lvec(0 20)\htext(5 22.5){$1$} \esegment
\end{texdraw}} \hskip 2mm +\hskip 2mm ^{q(1+q^6)}
\hskip 2mm\raisebox{-0.6\height}{\begin{texdraw} \drawdim mm
\setunitscale 0.5 \fontsize{7}{7}\selectfont \textref h:C v:C
\move(0 10)\tris \move(10 10)\tris \move(-20 10)\tris \move(-10
10)\tris \move(-20 0)\bsegment \move(0 10)\lvec(10 10)\lvec(10
20)\lvec(0 20)\lvec(0 10)\move(0 10)\lvec(10 20)\htext(2.5
17.5){$0$}\htext(7.5 12.5){$2$} \esegment

\move(-10 0)\bsegment \move(0 10)\lvec(10 10)\lvec(10 20)\lvec(0
20)\lvec(0 10)\move(0 10)\lvec(10 20)\htext(2.5
17.5){$2$}\htext(7.5 12.5){$0$} \move(0 20)\lvec(10 20)\lvec(10
25)\lvec(0 25)\lvec(0 20)\htext(5 22.5){$1$} \esegment

\move(0 0)\bsegment \move(0 10)\lvec(10 10)\lvec(10 20)\lvec(0
20)\lvec(0 10)\move(0 10)\lvec(10 20)\htext(2.5
17.5){$0$}\htext(7.5 12.5){$2$} \move(0 20)\lvec(10 20)\lvec(10
25)\lvec(0 25)\lvec(0 20)\htext(5 22.5){$1$} \esegment

\move(10 0)\bsegment \move(0 10)\lvec(10 10)\lvec(10 20)\lvec(0
20)\lvec(0 10)\move(0 10)\lvec(10 20)\htext(2.5
17.5){$2$}\htext(7.5 12.5){$0$} \move(0 20)\lvec(10 20)\lvec(10
25)\lvec(0 25)\lvec(0 20)\htext(5 22.5){$1$}

\move(0 25)\lvec(0 30)\lvec(10 30)\lvec(10 25)\htext(5 27.5){$1$}
 \esegment
\end{texdraw}} \hskip 2mm +\hskip 2mm ^{q(1-q^4)}
\hskip 2mm\raisebox{-0.4\height}{\begin{texdraw} \drawdim mm
\setunitscale 0.5 \fontsize{7}{7}\selectfont \textref h:C v:C
\move(0 10)\tris \move(10 10)\tris \move (-10 10)\tris \move(-10
0)\bsegment \move(0 10)\lvec(10 10)\lvec(10 20)\lvec(0 20)\lvec(0
10)\move(0 10)\lvec(10 20)\htext(2.5 17.5){$2$}\htext(7.5
12.5){$0$} \move(0 20)\lvec(10 20)\lvec(10 25)\lvec(0 25)\lvec(0
20)\htext(5 22.5){$1$} \esegment

\move(0 0)\bsegment \move(0 10)\lvec(10 10)\lvec(10 20)\lvec(0
20)\lvec(0 10)\move(0 10)\lvec(10 20)\htext(2.5
17.5){$0$}\htext(7.5 12.5){$2$} \move(0 20)\lvec(10 20)\lvec(10
25)\lvec(0 25)\lvec(0 20)\htext(5 22.5){$1$} \esegment

\move(10 0)\bsegment \move(0 10)\lvec(10 10)\lvec(10 20)\lvec(0
20)\lvec(0 10)\move(0 10)\lvec(10 20)\htext(2.5
17.5){$2$}\htext(7.5 12.5){$0$} \move(0 20)\lvec(10 20)\lvec(10
25)\lvec(0 25)\lvec(0 20)\htext(5 22.5){$1$} \move(0 25)\lvec(0
30)\lvec(10 30)\lvec(10 25)\htext(5 27.5){$1$}

\move(0 30)\lvec(10 30)\lvec(10 40)\lvec(0 30)\move(0 30)\lvec(10
40)\htext(7.5 32.5){$0$}
 \esegment
\end{texdraw}}
$

$+\hskip 2mm ^{q^2(1+q^2)(1-q^4)} \hskip
2mm\raisebox{-0.4\height}{\begin{texdraw} \drawdim mm
\setunitscale 0.5 \fontsize{7}{7}\selectfont \textref h:C v:C
\move(0 10)\tris \move(10 10)\tris \move (-10 10)\tris

\move(-10 0)\bsegment \move(0 10)\lvec(10 10)\lvec(10 20)\lvec(0
20)\lvec(0 10)\move(0 10)\lvec(10 20)\htext(2.5
17.5){$2$}\htext(7.5 12.5){$0$} \esegment

\move(0 0)\bsegment \move(0 10)\lvec(10 10)\lvec(10 20)\lvec(0
20)\lvec(0 10)\move(0 10)\lvec(10 20)\htext(2.5
17.5){$0$}\htext(7.5 12.5){$2$} \move(0 20)\lvec(10 20)\lvec(10
25)\lvec(0 25)\lvec(0 20)\htext(5 22.5){$1$} \move(0 25)\lvec(0
30)\lvec(10 30)\lvec(10 25)\htext(5 27.5){$1$} \esegment

\move(10 0)\bsegment \move(0 10)\lvec(10 10)\lvec(10 20)\lvec(0
20)\lvec(0 10)\move(0 10)\lvec(10 20)\htext(2.5
17.5){$2$}\htext(7.5 12.5){$0$} \move(0 20)\lvec(10 20)\lvec(10
25)\lvec(0 25)\lvec(0 20)\htext(5 22.5){$1$} \move(0 25)\lvec(0
30)\lvec(10 30)\lvec(10 25)\htext(5 27.5){$1$}

\move(0 30)\lvec(10 30)\lvec(10 40)\lvec(0 30)\move(0 30)\lvec(10
40)\htext(7.5 32.5){$0$} \esegment
\end{texdraw}}\hskip 2mm
+\hskip 2mm ^{(1+q^2)^2} \hskip
2mm\raisebox{-0.4\height}{\begin{texdraw} \drawdim mm
\setunitscale 0.5 \fontsize{7}{7}\selectfont \textref h:C v:C
\move(0 10)\tris \move(10 10)\tris

\move(0 0)\bsegment \move(0 10)\lvec(10 10)\lvec(10 20)\lvec(0
20)\lvec(0 10)\move(0 10)\lvec(10 20)\htext(2.5
17.5){$0$}\htext(7.5 12.5){$2$} \move(0 20)\lvec(10 20)\lvec(10
25)\lvec(0 25)\lvec(0 20)\htext(5 22.5){$1$} \move(0 25)\lvec(0
30)\lvec(10 30)\lvec(10 25)\htext(5 27.5){$1$} \esegment

\move(10 0)\bsegment \move(0 10)\lvec(10 10)\lvec(10 20)\lvec(0
20)\lvec(0 10)\move(0 10)\lvec(10 20)\htext(2.5
17.5){$2$}\htext(7.5 12.5){$0$} \move(0 20)\lvec(10 20)\lvec(10
25)\lvec(0 25)\lvec(0 20)\htext(5 22.5){$1$} \move(0 25)\lvec(0
30)\lvec(10 30)\lvec(10 25)\htext(5 27.5){$1$}

\move(0 30)\lvec(10 30)\lvec(10 40)\lvec(0 40)\lvec(0 30)\move(0
30)\lvec(10 40)\htext(2.5 37.5){$2$}\htext(7.5 32.5){$0$}\esegment
\end{texdraw}}
\hskip 2mm +\hskip 2mm ^{q(1+q^2)} \hskip
2mm\raisebox{-0.35\height}{\begin{texdraw} \drawdim mm
\setunitscale 0.5 \fontsize{7}{7}\selectfont \textref h:C v:C
\move(0 10)\tris \move(10 10)\tris

\move(0 0)\bsegment \move(0 10)\lvec(10 10)\lvec(10 20)\lvec(0
20)\lvec(0 10)\move(0 10)\lvec(10 20)\htext(2.5
17.5){$0$}\htext(7.5 12.5){$2$} \move(0 20)\lvec(10 20)\lvec(10
25)\lvec(0 25)\lvec(0 20)\htext(5 22.5){$1$}  \esegment

\move(10 0)\bsegment \move(0 10)\lvec(10 10)\lvec(10 20)\lvec(0
20)\lvec(0 10)\move(0 10)\lvec(10 20)\htext(2.5
17.5){$2$}\htext(7.5 12.5){$0$} \move(0 20)\lvec(10 20)\lvec(10
25)\lvec(0 25)\lvec(0 20)\htext(5 22.5){$1$} \move(0 25)\lvec(0
30)\lvec(10 30)\lvec(10 25)\htext(5 27.5){$1$}

\move(0 30)\lvec(10 30)\lvec(10 40)\lvec(0 40)\lvec(0 30)\move(0
30)\lvec(10 40)\htext(2.5 37.5){$2$}\htext(7.5 32.5){$0$}\move(0
40)\lvec(0 45)\lvec(10 45)\lvec(10 40)\htext(5 42.5){$1$}\esegment
\end{texdraw}}\hskip 3mm.$ \vskip 5mm

On the other hand, we have\vskip 5mm

$A(~\raisebox{-0.5\height}{\begin{texdraw} \drawdim mm
\setunitscale 0.5 \fontsize{7}{7}\selectfont \textref h:C v:C
\move(0 10)\tris \move(10 10)\tris

\move(0 0)\bsegment \move(0 10)\lvec(10 10)\lvec(10 20)\lvec(0
20)\lvec(0 10)\move(0 10)\lvec(10 20)\htext(2.5
17.5){$0$}\htext(7.5 12.5){$2$} \move(0 20)\lvec(10 20)\lvec(10
25)\lvec(0 25)\lvec(0 20)\htext(5 22.5){$1$} \move(0 25)\lvec(0
30)\lvec(10 30)\lvec(10 25)\htext(5 27.5){$1$} \esegment

\move(10 0)\bsegment \move(0 10)\lvec(10 10)\lvec(10 20)\lvec(0
20)\lvec(0 10)\move(0 10)\lvec(10 20)\htext(2.5
17.5){$2$}\htext(7.5 12.5){$0$} \move(0 20)\lvec(10 20)\lvec(10
25)\lvec(0 25)\lvec(0 20)\htext(5 22.5){$1$} \move(0 25)\lvec(0
30)\lvec(10 30)\lvec(10 25)\htext(5 27.5){$1$}

\move(0 30)\lvec(10 30)\lvec(10 40)\lvec(0 40)\lvec(0 30)\move(0
30)\lvec(10 40)\htext(2.5 37.5){$2$}\htext(7.5 32.5){$0$}\esegment
\end{texdraw}}~)=f_1^{(2)}f_2f_0^{(2)}f_1^{(2)}f_2 Y_{\Lambda_2}$

\hskip 2cm $=\hskip 2mm \raisebox{-0.5\height}{\begin{texdraw}
\drawdim mm \setunitscale 0.5 \fontsize{7}{7}\selectfont \textref
h:C v:C \move(0 10)\tris \move(10 10)\tris

\move(0 0)\bsegment \move(0 10)\lvec(10 10)\lvec(10 20)\lvec(0
20)\lvec(0 10)\move(0 10)\lvec(10 20)\htext(2.5
17.5){$0$}\htext(7.5 12.5){$2$} \move(0 20)\lvec(10 20)\lvec(10
25)\lvec(0 25)\lvec(0 20)\htext(5 22.5){$1$} \move(0 25)\lvec(0
30)\lvec(10 30)\lvec(10 25)\htext(5 27.5){$1$} \esegment

\move(10 0)\bsegment \move(0 10)\lvec(10 10)\lvec(10 20)\lvec(0
20)\lvec(0 10)\move(0 10)\lvec(10 20)\htext(2.5
17.5){$2$}\htext(7.5 12.5){$0$} \move(0 20)\lvec(10 20)\lvec(10
25)\lvec(0 25)\lvec(0 20)\htext(5 22.5){$1$} \move(0 25)\lvec(0
30)\lvec(10 30)\lvec(10 25)\htext(5 27.5){$1$}

\move(0 30)\lvec(10 30)\lvec(10 40)\lvec(0 40)\lvec(0 30)\move(0
30)\lvec(10 40)\htext(2.5 37.5){$2$}\htext(7.5 32.5){$0$}\esegment
\end{texdraw}} \hskip 2mm +
\hskip 2mm q \hskip 2mm \raisebox{-0.5\height}{\begin{texdraw}
\drawdim mm \setunitscale 0.5 \fontsize{7}{7}\selectfont \textref
h:C v:C \move(0 10)\tris \move(10 10)\tris

\move(0 0)\bsegment \move(0 10)\lvec(10 10)\lvec(10 20)\lvec(0
20)\lvec(0 10)\move(0 10)\lvec(10 20)\htext(2.5
17.5){$0$}\htext(7.5 12.5){$2$} \move(0 20)\lvec(10 20)\lvec(10
25)\lvec(0 25)\lvec(0 20)\htext(5 22.5){$1$}  \esegment

\move(10 0)\bsegment \move(0 10)\lvec(10 10)\lvec(10 20)\lvec(0
20)\lvec(0 10)\move(0 10)\lvec(10 20)\htext(2.5
17.5){$2$}\htext(7.5 12.5){$0$} \move(0 20)\lvec(10 20)\lvec(10
25)\lvec(0 25)\lvec(0 20)\htext(5 22.5){$1$} \move(0 25)\lvec(0
30)\lvec(10 30)\lvec(10 25)\htext(5 27.5){$1$}

\move(0 30)\lvec(10 30)\lvec(10 40)\lvec(0 40)\lvec(0 30)\move(0
30)\lvec(10 40)\htext(2.5 37.5){$2$}\htext(7.5 32.5){$0$}

\move(0 40)\lvec(0 45)\lvec(10 45)\lvec(10 40)\htext(5 42.5){$1$}
\esegment
\end{texdraw}}
\hskip 2mm + \hskip 2mm q^4 \hskip 2mm
\raisebox{-0.45\height}{\begin{texdraw} \drawdim mm \setunitscale
0.5 \fontsize{7}{7}\selectfont \textref h:C v:C \move(0 10)\tris
\move(10 10)\tris

\move(0 0)\bsegment \move(0 10)\lvec(10 10)\lvec(10 20)\lvec(0
20)\lvec(0 10)\move(0 10)\lvec(10 20)\htext(2.5
17.5){$0$}\htext(7.5 12.5){$2$} \esegment

\move(10 0)\bsegment \move(0 10)\lvec(10 10)\lvec(10 20)\lvec(0
20)\lvec(0 10)\move(0 10)\lvec(10 20)\htext(2.5
17.5){$2$}\htext(7.5 12.5){$0$} \move(0 20)\lvec(10 20)\lvec(10
25)\lvec(0 25)\lvec(0 20)\htext(5 22.5){$1$} \move(0 25)\lvec(0
30)\lvec(10 30)\lvec(10 25)\htext(5 27.5){$1$}

\move(0 30)\lvec(10 30)\lvec(10 40)\lvec(0 40)\lvec(0 30)\move(0
30)\lvec(10 40)\htext(2.5 37.5){$2$}\htext(7.5 32.5){$0$}

\move(0 40)\lvec(10 40)\lvec(10 45)\lvec(0 45)\lvec(0 40)\htext(5
42.5){$1$} \move(0 45)\lvec(0 50)\lvec(10 50)\lvec(10 45)\htext(5
47.5){$1$} \esegment
\end{texdraw}}
=G(~\raisebox{-0.5\height}{\begin{texdraw} \drawdim mm
\setunitscale 0.5 \fontsize{7}{7}\selectfont \textref h:C v:C
\move(0 10)\tris \move(10 10)\tris

\move(0 0)\bsegment \move(0 10)\lvec(10 10)\lvec(10 20)\lvec(0
20)\lvec(0 10)\move(0 10)\lvec(10 20)\htext(2.5
17.5){$0$}\htext(7.5 12.5){$2$} \move(0 20)\lvec(10 20)\lvec(10
25)\lvec(0 25)\lvec(0 20)\htext(5 22.5){$1$} \move(0 25)\lvec(0
30)\lvec(10 30)\lvec(10 25)\htext(5 27.5){$1$} \esegment

\move(10 0)\bsegment \move(0 10)\lvec(10 10)\lvec(10 20)\lvec(0
20)\lvec(0 10)\move(0 10)\lvec(10 20)\htext(2.5
17.5){$2$}\htext(7.5 12.5){$0$} \move(0 20)\lvec(10 20)\lvec(10
25)\lvec(0 25)\lvec(0 20)\htext(5 22.5){$1$} \move(0 25)\lvec(0
30)\lvec(10 30)\lvec(10 25)\htext(5 27.5){$1$}

\move(0 30)\lvec(10 30)\lvec(10 40)\lvec(0 40)\lvec(0 30)\move(0
30)\lvec(10 40)\htext(2.5 37.5){$2$}\htext(7.5 32.5){$0$}\esegment
\end{texdraw}}~).$\vskip 5mm

Therefore, we obtain \vskip 5mm

$G(~\raisebox{-0.6\height}{\begin{texdraw} \drawdim mm
\setunitscale 0.5 \fontsize{7}{7}\selectfont \textref h:C v:C
\move(0 10)\tris \move(10 10)\tris \move(-20 10)\tris \move(-10
10)\tris \move(-20 0)\bsegment \move(0 10)\lvec(10 10)\lvec(10
20)\lvec(0 20)\lvec(0 10)\move(0 10)\lvec(10 20)\htext(2.5
17.5){$0$}\htext(7.5 12.5){$2$} \move(0 20)\lvec(10 20)\lvec(10
25)\lvec(0 25)\lvec(0 20)\htext(5 22.5){$1$} \esegment

\move(-10 0)\bsegment \move(0 10)\lvec(10 10)\lvec(10 20)\lvec(0
20)\lvec(0 10)\move(0 10)\lvec(10 20)\htext(2.5
17.5){$2$}\htext(7.5 12.5){$0$} \move(0 20)\lvec(10 20)\lvec(10
25)\lvec(0 25)\lvec(0 20)\htext(5 22.5){$1$} \esegment

\move(0 0)\bsegment \move(0 10)\lvec(10 10)\lvec(10 20)\lvec(0
20)\lvec(0 10)\move(0 10)\lvec(10 20)\htext(2.5
17.5){$0$}\htext(7.5 12.5){$2$} \move(0 20)\lvec(10 20)\lvec(10
25)\lvec(0 25)\lvec(0 20)\htext(5 22.5){$1$} \esegment

\move(10 0)\bsegment \move(0 10)\lvec(10 10)\lvec(10 20)\lvec(0
20)\lvec(0 10)\move(0 10)\lvec(10 20)\htext(2.5
17.5){$2$}\htext(7.5 12.5){$0$} \move(0 20)\lvec(10 20)\lvec(10
25)\lvec(0 25)\lvec(0 20)\htext(5 22.5){$1$} \esegment
\end{texdraw}}~)=A(~\raisebox{-0.6\height}{\begin{texdraw} \drawdim mm
\setunitscale 0.5 \fontsize{7}{7}\selectfont \textref h:C v:C
\move(0 10)\tris \move(10 10)\tris \move(-20 10)\tris \move(-10
10)\tris \move(-20 0)\bsegment \move(0 10)\lvec(10 10)\lvec(10
20)\lvec(0 20)\lvec(0 10)\move(0 10)\lvec(10 20)\htext(2.5
17.5){$0$}\htext(7.5 12.5){$2$} \move(0 20)\lvec(10 20)\lvec(10
25)\lvec(0 25)\lvec(0 20)\htext(5 22.5){$1$} \esegment

\move(-10 0)\bsegment \move(0 10)\lvec(10 10)\lvec(10 20)\lvec(0
20)\lvec(0 10)\move(0 10)\lvec(10 20)\htext(2.5
17.5){$2$}\htext(7.5 12.5){$0$} \move(0 20)\lvec(10 20)\lvec(10
25)\lvec(0 25)\lvec(0 20)\htext(5 22.5){$1$} \esegment

\move(0 0)\bsegment \move(0 10)\lvec(10 10)\lvec(10 20)\lvec(0
20)\lvec(0 10)\move(0 10)\lvec(10 20)\htext(2.5
17.5){$0$}\htext(7.5 12.5){$2$} \move(0 20)\lvec(10 20)\lvec(10
25)\lvec(0 25)\lvec(0 20)\htext(5 22.5){$1$} \esegment

\move(10 0)\bsegment \move(0 10)\lvec(10 10)\lvec(10 20)\lvec(0
20)\lvec(0 10)\move(0 10)\lvec(10 20)\htext(2.5
17.5){$2$}\htext(7.5 12.5){$0$} \move(0 20)\lvec(10 20)\lvec(10
25)\lvec(0 25)\lvec(0 20)\htext(5 22.5){$1$} \esegment
\end{texdraw}}~)-G(~\raisebox{-0.5\height}{\begin{texdraw} \drawdim mm
\setunitscale 0.5 \fontsize{7}{7}\selectfont \textref h:C v:C
\move(0 10)\tris \move(10 10)\tris

\move(0 0)\bsegment \move(0 10)\lvec(10 10)\lvec(10 20)\lvec(0
20)\lvec(0 10)\move(0 10)\lvec(10 20)\htext(2.5
17.5){$0$}\htext(7.5 12.5){$2$} \move(0 20)\lvec(10 20)\lvec(10
25)\lvec(0 25)\lvec(0 20)\htext(5 22.5){$1$} \move(0 25)\lvec(0
30)\lvec(10 30)\lvec(10 25)\htext(5 27.5){$1$} \esegment

\move(10 0)\bsegment \move(0 10)\lvec(10 10)\lvec(10 20)\lvec(0
20)\lvec(0 10)\move(0 10)\lvec(10 20)\htext(2.5
17.5){$2$}\htext(7.5 12.5){$0$} \move(0 20)\lvec(10 20)\lvec(10
25)\lvec(0 25)\lvec(0 20)\htext(5 22.5){$1$} \move(0 25)\lvec(0
30)\lvec(10 30)\lvec(10 25)\htext(5 27.5){$1$}

\move(0 30)\lvec(10 30)\lvec(10 40)\lvec(0 40)\lvec(0 30)\move(0
30)\lvec(10 40)\htext(2.5 37.5){$2$}\htext(7.5 32.5){$0$}\esegment
\end{texdraw}}~)$

$=\hskip 2mm\raisebox{-0.6\height}{\begin{texdraw} \drawdim mm
\setunitscale 0.5 \fontsize{7}{7}\selectfont \textref h:C v:C
\move(0 10)\tris \move(10 10)\tris \move(-20 10)\tris \move(-10
10)\tris \move(-20 0)\bsegment \move(0 10)\lvec(10 10)\lvec(10
20)\lvec(0 20)\lvec(0 10)\move(0 10)\lvec(10 20)\htext(2.5
17.5){$0$}\htext(7.5 12.5){$2$} \move(0 20)\lvec(10 20)\lvec(10
25)\lvec(0 25)\lvec(0 20)\htext(5 22.5){$1$} \esegment

\move(-10 0)\bsegment \move(0 10)\lvec(10 10)\lvec(10 20)\lvec(0
20)\lvec(0 10)\move(0 10)\lvec(10 20)\htext(2.5
17.5){$2$}\htext(7.5 12.5){$0$} \move(0 20)\lvec(10 20)\lvec(10
25)\lvec(0 25)\lvec(0 20)\htext(5 22.5){$1$} \esegment

\move(0 0)\bsegment \move(0 10)\lvec(10 10)\lvec(10 20)\lvec(0
20)\lvec(0 10)\move(0 10)\lvec(10 20)\htext(2.5
17.5){$0$}\htext(7.5 12.5){$2$} \move(0 20)\lvec(10 20)\lvec(10
25)\lvec(0 25)\lvec(0 20)\htext(5 22.5){$1$} \esegment

\move(10 0)\bsegment \move(0 10)\lvec(10 10)\lvec(10 20)\lvec(0
20)\lvec(0 10)\move(0 10)\lvec(10 20)\htext(2.5
17.5){$2$}\htext(7.5 12.5){$0$} \move(0 20)\lvec(10 20)\lvec(10
25)\lvec(0 25)\lvec(0 20)\htext(5 22.5){$1$} \esegment
\end{texdraw}} \hskip 2mm +\hskip 2mm ^{q(1+q^6)}
\hskip 2mm\raisebox{-0.6\height}{\begin{texdraw} \drawdim mm
\setunitscale 0.5 \fontsize{7}{7}\selectfont \textref h:C v:C
\move(0 10)\tris \move(10 10)\tris \move(-20 10)\tris \move(-10
10)\tris \move(-20 0)\bsegment \move(0 10)\lvec(10 10)\lvec(10
20)\lvec(0 20)\lvec(0 10)\move(0 10)\lvec(10 20)\htext(2.5
17.5){$0$}\htext(7.5 12.5){$2$} \esegment

\move(-10 0)\bsegment \move(0 10)\lvec(10 10)\lvec(10 20)\lvec(0
20)\lvec(0 10)\move(0 10)\lvec(10 20)\htext(2.5
17.5){$2$}\htext(7.5 12.5){$0$} \move(0 20)\lvec(10 20)\lvec(10
25)\lvec(0 25)\lvec(0 20)\htext(5 22.5){$1$} \esegment

\move(0 0)\bsegment \move(0 10)\lvec(10 10)\lvec(10 20)\lvec(0
20)\lvec(0 10)\move(0 10)\lvec(10 20)\htext(2.5
17.5){$0$}\htext(7.5 12.5){$2$} \move(0 20)\lvec(10 20)\lvec(10
25)\lvec(0 25)\lvec(0 20)\htext(5 22.5){$1$} \esegment

\move(10 0)\bsegment \move(0 10)\lvec(10 10)\lvec(10 20)\lvec(0
20)\lvec(0 10)\move(0 10)\lvec(10 20)\htext(2.5
17.5){$2$}\htext(7.5 12.5){$0$} \move(0 20)\lvec(10 20)\lvec(10
25)\lvec(0 25)\lvec(0 20)\htext(5 22.5){$1$}

\move(0 25)\lvec(0 30)\lvec(10 30)\lvec(10 25)\htext(5 27.5){$1$}
 \esegment
\end{texdraw}} \hskip 2mm +\hskip 2mm ^{q(1-q^4)}
\hskip 2mm\raisebox{-0.4\height}{\begin{texdraw} \drawdim mm
\setunitscale 0.5 \fontsize{7}{7}\selectfont \textref h:C v:C
\move(0 10)\tris \move(10 10)\tris \move (-10 10)\tris \move(-10
0)\bsegment \move(0 10)\lvec(10 10)\lvec(10 20)\lvec(0 20)\lvec(0
10)\move(0 10)\lvec(10 20)\htext(2.5 17.5){$2$}\htext(7.5
12.5){$0$} \move(0 20)\lvec(10 20)\lvec(10 25)\lvec(0 25)\lvec(0
20)\htext(5 22.5){$1$} \esegment

\move(0 0)\bsegment \move(0 10)\lvec(10 10)\lvec(10 20)\lvec(0
20)\lvec(0 10)\move(0 10)\lvec(10 20)\htext(2.5
17.5){$0$}\htext(7.5 12.5){$2$} \move(0 20)\lvec(10 20)\lvec(10
25)\lvec(0 25)\lvec(0 20)\htext(5 22.5){$1$} \esegment

\move(10 0)\bsegment \move(0 10)\lvec(10 10)\lvec(10 20)\lvec(0
20)\lvec(0 10)\move(0 10)\lvec(10 20)\htext(2.5
17.5){$2$}\htext(7.5 12.5){$0$} \move(0 20)\lvec(10 20)\lvec(10
25)\lvec(0 25)\lvec(0 20)\htext(5 22.5){$1$} \move(0 25)\lvec(0
30)\lvec(10 30)\lvec(10 25)\htext(5 27.5){$1$}

\move(0 30)\lvec(10 30)\lvec(10 40)\lvec(0 30)\move(0 30)\lvec(10
40)\htext(7.5 32.5){$0$}
 \esegment
\end{texdraw}}
$

$+\hskip 2mm ^{q^2(1+q^2)(1-q^4)} \hskip
2mm\raisebox{-0.4\height}{\begin{texdraw} \drawdim mm
\setunitscale 0.5 \fontsize{7}{7}\selectfont \textref h:C v:C
\move(0 10)\tris \move(10 10)\tris \move (-10 10)\tris

\move(-10 0)\bsegment \move(0 10)\lvec(10 10)\lvec(10 20)\lvec(0
20)\lvec(0 10)\move(0 10)\lvec(10 20)\htext(2.5
17.5){$2$}\htext(7.5 12.5){$0$} \esegment

\move(0 0)\bsegment \move(0 10)\lvec(10 10)\lvec(10 20)\lvec(0
20)\lvec(0 10)\move(0 10)\lvec(10 20)\htext(2.5
17.5){$0$}\htext(7.5 12.5){$2$} \move(0 20)\lvec(10 20)\lvec(10
25)\lvec(0 25)\lvec(0 20)\htext(5 22.5){$1$} \move(0 25)\lvec(0
30)\lvec(10 30)\lvec(10 25)\htext(5 27.5){$1$} \esegment

\move(10 0)\bsegment \move(0 10)\lvec(10 10)\lvec(10 20)\lvec(0
20)\lvec(0 10)\move(0 10)\lvec(10 20)\htext(2.5
17.5){$2$}\htext(7.5 12.5){$0$} \move(0 20)\lvec(10 20)\lvec(10
25)\lvec(0 25)\lvec(0 20)\htext(5 22.5){$1$} \move(0 25)\lvec(0
30)\lvec(10 30)\lvec(10 25)\htext(5 27.5){$1$}

\move(0 30)\lvec(10 30)\lvec(10 40)\lvec(0 30)\move(0 30)\lvec(10
40)\htext(7.5 32.5){$0$} \esegment
\end{texdraw}}\hskip 2mm
+\hskip 2mm ^{q^2(2+q^2)} \hskip
2mm\raisebox{-0.4\height}{\begin{texdraw} \drawdim mm
\setunitscale 0.5 \fontsize{7}{7}\selectfont \textref h:C v:C
\move(0 10)\tris \move(10 10)\tris

\move(0 0)\bsegment \move(0 10)\lvec(10 10)\lvec(10 20)\lvec(0
20)\lvec(0 10)\move(0 10)\lvec(10 20)\htext(2.5
17.5){$0$}\htext(7.5 12.5){$2$} \move(0 20)\lvec(10 20)\lvec(10
25)\lvec(0 25)\lvec(0 20)\htext(5 22.5){$1$} \move(0 25)\lvec(0
30)\lvec(10 30)\lvec(10 25)\htext(5 27.5){$1$} \esegment

\move(10 0)\bsegment \move(0 10)\lvec(10 10)\lvec(10 20)\lvec(0
20)\lvec(0 10)\move(0 10)\lvec(10 20)\htext(2.5
17.5){$2$}\htext(7.5 12.5){$0$} \move(0 20)\lvec(10 20)\lvec(10
25)\lvec(0 25)\lvec(0 20)\htext(5 22.5){$1$} \move(0 25)\lvec(0
30)\lvec(10 30)\lvec(10 25)\htext(5 27.5){$1$}

\move(0 30)\lvec(10 30)\lvec(10 40)\lvec(0 40)\lvec(0 30)\move(0
30)\lvec(10 40)\htext(2.5 37.5){$2$}\htext(7.5 32.5){$0$}\esegment
\end{texdraw}}+\hskip 2mm ^{q^3} \hskip
2mm\raisebox{-0.35\height}{\begin{texdraw} \drawdim mm
\setunitscale 0.5 \fontsize{7}{7}\selectfont \textref h:C v:C
\move(0 10)\tris \move(10 10)\tris

\move(0 0)\bsegment \move(0 10)\lvec(10 10)\lvec(10 20)\lvec(0
20)\lvec(0 10)\move(0 10)\lvec(10 20)\htext(2.5
17.5){$0$}\htext(7.5 12.5){$2$} \move(0 20)\lvec(10 20)\lvec(10
25)\lvec(0 25)\lvec(0 20)\htext(5 22.5){$1$}  \esegment

\move(10 0)\bsegment \move(0 10)\lvec(10 10)\lvec(10 20)\lvec(0
20)\lvec(0 10)\move(0 10)\lvec(10 20)\htext(2.5
17.5){$2$}\htext(7.5 12.5){$0$} \move(0 20)\lvec(10 20)\lvec(10
25)\lvec(0 25)\lvec(0 20)\htext(5 22.5){$1$} \move(0 25)\lvec(0
30)\lvec(10 30)\lvec(10 25)\htext(5 27.5){$1$}

\move(0 30)\lvec(10 30)\lvec(10 40)\lvec(0 40)\lvec(0 30)\move(0
30)\lvec(10 40)\htext(2.5 37.5){$2$}\htext(7.5 32.5){$0$}\move(0
40)\lvec(0 45)\lvec(10 45)\lvec(10 40)\htext(5 42.5){$1$}\esegment
\end{texdraw}}$

$\hskip 2mm - \hskip 2mm ^{q^4} \hskip 2mm
\raisebox{-0.3\height}{\begin{texdraw} \drawdim mm \setunitscale
0.5 \fontsize{7}{7}\selectfont \textref h:C v:C \move(0 10)\tris
\move(10 10)\tris

\move(0 0)\bsegment \move(0 10)\lvec(10 10)\lvec(10 20)\lvec(0
20)\lvec(0 10)\move(0 10)\lvec(10 20)\htext(2.5
17.5){$0$}\htext(7.5 12.5){$2$} \esegment

\move(10 0)\bsegment \move(0 10)\lvec(10 10)\lvec(10 20)\lvec(0
20)\lvec(0 10)\move(0 10)\lvec(10 20)\htext(2.5
17.5){$2$}\htext(7.5 12.5){$0$} \move(0 20)\lvec(10 20)\lvec(10
25)\lvec(0 25)\lvec(0 20)\htext(5 22.5){$1$} \move(0 25)\lvec(0
30)\lvec(10 30)\lvec(10 25)\htext(5 27.5){$1$}

\move(0 30)\lvec(10 30)\lvec(10 40)\lvec(0 40)\lvec(0 30)\move(0
30)\lvec(10 40)\htext(2.5 37.5){$2$}\htext(7.5 32.5){$0$}

\move(0 40)\lvec(10 40)\lvec(10 45)\lvec(0 45)\lvec(0 40)\htext(5
42.5){$1$} \move(0 45)\lvec(0 50)\lvec(10 50)\lvec(10 45)\htext(5
47.5){$1$} \esegment
\end{texdraw}}\hskip 3mm.$ }
\end{ex}


{\small
\begin{thebibliography}{BK}
\bibitem{Ariki}
S.~Ariki, {\em On decomposition number of Hecke algebra of
$G(m,1,n)$}, J. Math. Kyoto Univ. \textbf{36} (1996), 789--808.


\bibitem{BK} J. Brundan, A. Kleshchev, {\em Hecke-Clifford superalgebras,
crystals of type $A^{(2)}_{2\ell}$ and modular branching rules for
$\hat{S}_n$}, Represent. Theory {\bf 5} (2001), 317--403
(electronic).

\bibitem{GL} I. Grojnowski, G. Lusztig,
{\em A comparison of bases of quantized enveloping algebras},
Contemp. Math. {\bf 153} (1993), 11-19.

\bibitem{HK} J. Hong, S.-J. Kang,  {\em Crystal graphs for basic representations
of the quantum affine algebra $U_q(C_2^{(1)})$}, Representations
and quantizations (Shanghai, 1998), 213--227, China High. Educ.
Press, Beijing, 2000.

\bibitem{Kang} S.-J.~Kang, {\em Crystal bases for quantum affine Lie algebras
and combinatorics of Young walls},  RIM-GARC preprint (2000) {\bf
00-2}, Seoul National University, to appear in Proc. London Math.
Soc.

\bibitem{KMN2} S.-J. Kang, M. Kashiwara, K. Misra, T. Miwa, T.
Nakashima, A. Nakayashiki, {\em Perfect crystals of quantum affine
Lie algebras}, Duke Math. J. {\bf 68} (1992), 499--607.

\bibitem{KK'}  S.-J. Kang, J.-H. Kwon, {\em Quantum affine algebras, combinatorics of Young
walls and global basis}, Electron. Res. Announc. Amer. Math. Soc. {\bf 8} (2002), 35-46.

\bibitem{KK} S.-J. Kang, J.-H. Kwon, {\em Fock space representations of quantum
affine algebras and generalized Lascoux-Leclerc-Thibon algorithm},
preprint (2002) math.QA/0208204.

\bibitem{Kash} M. Kashiwara, {\em On crystal bases of the $q$-analogue
of universal enveloping algebras}, Duke Math. J. {\bf 63} (1991),
465--516.

\bibitem{Kash93} M. Kashiwara, {\em Crystal bases and Littelmann's
refined Demazure character formula}, Duke Math. J. {\bf 71}
(1993), 839--858.

\bibitem{Kash94} M. Kashiwara, {\em Crystal bases of modified
quantized enveloping algebras}, Duke Math. J. {\bf 73} (1994),
383--413.

\bibitem{KMPY} M. Kashiwara, T. Miwa, J.-U. H. Petersen, C. M. Yung,
{\em Perfect crystals and $q$-deformed Fock spaces}, Selecta Math.
{\bf 2} (1996), 415--499.

\bibitem{KashNakash}
M. Kashiwara, T. Nakashima, {\em Crystal graphs for
representations of the $q$-analogue of classical Lie algebras}, J.
Algebra {\bf 165} (1994), 295--345.

\bibitem{KashSaito}
M. Kashiwara, Y. Saito, {\em Geometric construction of crystal
bases}, Duke Math. J. {\bf 89} (1997), 9--36.

\bibitem{LLT} A. Lascoux, B. Leclerc, J.-Y. Thibon, {\em Hecke algebras at
roots of unity and crystal bases of quantum affine algebras},
Comm. Math. Phys. {\bf 181} (1996), 205--263.

\bibitem{LT} B. Leclerc, P. Toffin,
{\em A simple algorithm for computing the global crystal basis of
an irreducible $U_q(\frak{sl}_n)$-module}, Int. J. Algebra
Computation, {\bf 10} (2000), 191-208.

\bibitem{L} C. Lecouvey, {\em An algorithm for computing the global basis of
an irreducible $U_q(\frak{sp}_{2n})$-module}, preprint (2002)
math.QA/0201143.

\bibitem{Lit}
P. Littelmann, {\em Paths and root operators in representation
theory}, Ann. of Math. (2) {\bf 142} (1995), 499--525.

\bibitem{Lusz} G. Lusztig, {\em Canonical bases arising from quantized
enveloping algebras}, J. Amer. Math. Soc. {\bf 3} (1990),
447--498.

\bibitem{Nakash1}
T. Nakashima, A. Zelevinsky, {\em Polyhedral realizations of
crystal bases for quantized Kac-Moody algebras}, Adv. Math. {\bf
131} (1997), 253--278.

\bibitem{Nakash2}
T. Nakashima, {\em Polyhedral realizations of crystal bases for
integrable highest weight modules}, J. Algebra {\bf 219} (1999),
571--597.


\end{thebibliography}
}
\end{document}